\documentclass[9pt,letterpaper,oneside,twocolumn,romanappendices]{IEEEtran}

\usepackage{amssymb}
\usepackage{amsthm, mathrsfs}
\usepackage{graphicx}
\usepackage{subfigure}
\usepackage{enumerate}
\usepackage{color}
\usepackage{arydshln}
\usepackage{multicol}
\usepackage{multirow}
\usepackage{hyperref}
\usepackage{algorithmic}
\usepackage[lined,ruled,linesnumbered]{algorithm2e}
\usepackage{microtype}
\usepackage{booktabs}
\usepackage{rotating}
\usepackage{booktabs}
\usepackage{float}

\usepackage{cite}

\newtheorem{theorem}{Theorem}[section]

\newtheorem{proposition}[theorem]{Proposition}

\ifCLASSINFOpdf
\else
\fi
\usepackage[cmex10]{amsmath}
\begin{document}
\title{High-Accuracy Total Variation for Compressed Video Sensing}

\author{Mahdi~S.~Hosseini,~\IEEEmembership{Student Member,~IEEE,}
        and~Konstantinos~N.~Plataniotis,~\IEEEmembership{Fellow,~IEEE}
\thanks{M. S. Hosseini and K. N. Plataniotis are with The Edward S. Rogers Sr. Department of Electrical and Computer Engineering, University of Toronto, Toronto, ON M5S 3G4, Canada e-mail: \url{mahdi.hosseini@mail.utoronto.ca}.}}%

\markboth{Submitted to IEEE Transaction on Image Processing}{Shell \MakeLowercase{\textit{et al.}}: Bare Demo of IEEEtran.cls for Journals}

\maketitle

\IEEEpeerreviewmaketitle

\begin{abstract}
Numerous total variation (TV) regularizers, engaged in image restoration problem, encode the gradients by means of simple $[-1,1]$ FIR filter. Despite its low computational processing, this filter severely deviates signal's high frequency components pertinent to edge/discontinuous information and cause several deficiency issues known as texture and geometric loss. This paper addresses this problem by proposing an alternative model to the TV regularization problem via high order accuracy differential FIR filters to preserve rapid transitions in signal recovery. A numerical encoding scheme is designed to extend the TV model into multidimensional representation (tensorial decomposition). We adopt this design to regulate the spatial and temporal redundancy in compressed video sensing problem to jointly recover frames from under-sampled measurements. We then seek the solution via alternating direction methods of multipliers and find a unique solution to quadratic minimization step with capability of handling different boundary conditions. The resulting algorithm uses much lower sampling rate and highly outperforms alternative state-of-the-art methods. This is evaluated both in terms of restoration accuracy and visual quality of the recovered frames.
\end{abstract}


\IEEEpeerreviewmaketitle

\section{Introduction}\label{sec:intro}
\IEEEPARstart{T}{otal} variation (TV) regularization, proposed by Rudin et al. \cite{RudinOsherFatemi:1992}, plays a vital role in image restoration problem to restore signals from their degraded measurements such as denoising \cite{RudinOsherFatemi:1992, chan2006total, ChanMarquinaMulet:2000, BrediesKunischPock:2010, BenningBruneBurgerMuller:2013}, compressed sensing (CS) \cite{CandesRombergTao1:2006, AfonsoBioucas-DiasFigueiredo2011} and deconvolution \cite{AfonsoBioucas-DiasFigueiredo2011, FigueiredoBioucasDias:2010, AlmeidaFigueiredo:2013, MatakosRamaniFessler:2013}. In such applications, TV regularizer impose stability to the recovery by preserving rapid transitions and discontinuities in the signal while prevents oscillation in the reconstruction level \cite{RudinOsherFatemi:1992, chan2006total}. The TV model in \cite{RudinOsherFatemi:1992} admits sparsity for piecewise constant signals and tends to underperform in smooth transitions recoveries causing blocky reconstructions: known as the stair-case effect \cite{chan2006total}.  Improvements via high order derivatives have been recently proposed in conjunction with a prior model in \cite{RudinOsherFatemi:1992} to reduce this artefact and provide more smooth transitions in the recovered signal \cite{ChanMarquinaMulet:2000, BrediesKunischPock:2010, BenningBruneBurgerMuller:2013}. The latter approaches, however, are case sensitive and suffer from deficiencies such as loss of geometry and textures. Such deficiencies still remain an open question to the field which needs to be properly addressed, please refer to \cite{chan2006total, ChanMarquinaMulet:2000, BrediesKunischPock:2010, BenningBruneBurgerMuller:2013} and reference therein for further explanation on deficiency issues.

TV encodes signal' features by means of first order differentiation. The most popular approach in discrete domain is to simply approximate the derivatives by $[-1,1]$ FIR filter with low computational complexity \cite{RudinOsherFatemi:1992}. Despite its popularity, this filter retains the lowest accuracy in numerical difference methods causing high frequency deviations in the filter response \cite{Lele:1992, fornberg:1988, Li:2005, Khan:1999, Khan:2000}. This property is in contrary with the notion of TV, which favours discontinuities referring to the high frequency components in the signal. In this article, we revisit the problem of TV regularization by means of high order (HO) accuracy differentiation in \cite{Holoborodko:2008}. The major benefit of HO accuracy FIR filter is to accurately transfer the high frequency components in estimation of the signal features. This is of great importance in video restoration problem, where the signal transitions occur either in spatial resolution: object boundaries and textural information, or temporal resolution: object movements/motion characteristics. Therefore, encoding such features via HO accuracy filters for approximating the TV norm is critical for better reconstructions.

Particular application of the proposed TV regularizer is discussed in this paper for the problem of compressed video sensing (CVS) an emerging technology for performing accelerated acquisitions from motion imagery. In this technology, under-sampled measurements are collected from multiple frames for reconstruction by exploiting significant temporal redundancy in video. Applications include, but are not limited to, commercial video processing cameras \cite{DuarteDavenportBaraniuk:2008, OikeGamal:2012, OikeGamal:2013}, video surveillance \cite{JiangLiHaimiCohenWilfordZhang:2012, Marcia_Willett_2008}, and dynamic MRI for medical imaging \cite{JungSungNayakKimYeChul:2009}. To date, CVS is explored in two main categories:
\begin{itemize}
\item motion estimation (ME) / motion compensation (MC): temporal correlation is encoded by an ME technique followed by compensating estimated motion MC in reconstruction level. ME/MC has been vastly used in CVS by altering between ME and MC for recovering frames from CS measurements, e.g. lifting-based multiscale frameworks \cite{ParkWakin:2009, ParkWakin:2013}, linear dynamic systems  \cite{AsifHamiltonBrummerRomberg:2012}, block-based approaches  \cite{JungSungNayakKimYeChul:2009, LiuLiPados:2013}, optical flow methods \cite{SankaranarayananStuderBaraniuk:2012}, and motion detection algorithms \cite{KashterLeviStern:2012}. 
\item tensorial decomposition: motion trajectories are jointly encoded with spatial features in each frame by means of tensorial decomposition \cite{DuarteBaraniuk:2011, CaiafaCichocki:2012}. The first idea is triggered in \cite{Marcia_Willett_2008} for super-resolution video problem, where the sparse decomposition is carried out by a temporal integration to {jointly recover the reference frame as well as the difference frames from under-sampled measurements.} {Alternative approaches in \cite{ShuAhuja:2011, JiangLiHaimiCohenWilfordZhang:2012} encode the frame differences in video by means TV regularizer.}
\end{itemize}

{ME/MC, adopted from video compression \cite{LeGall:1991}, modifies the CVS performance by compensating the motion residuals in the recovery and, hence, complements the low accuracy estimation of the temporal gradients in video. However, the motion between consecutive frames is again estimated by simple $[-1,1]$ filter and also requires a reference frame for recovery. Tensorial decomposition, compared to ME/MC, is more straight forward which is capable of separately encoding features in both space and time domains and avoid multistage recovery. Our aim here is to utilize the revised TV regularizer in CVS by means of tensorial decomposition for spatial-temporal feature encoding. The main benefit of such encoding in video is to preserve high frequency components in both space and time domains.}

In addition, this work offers a solution to a compound minimization problem that employs the proposed TV regularization and jointly reconstructs video frames from under-sampled measurements. {We find the solution to the minimization problem by utilizing the method of} split variable technique using the alternating direction methods of multipliers (ADMM) \cite{GabayMercier:1976, citeulike:10118819}. Conventional ADMM approaches use either periodic (circular) or reflective (Neumann) BC in TV regularizer for differentiation \cite{YangZhangYin:2009, AfonsoBioucas-DiasFigueiredo2011, FigueiredoBioucasDias:2010}. {Recent developments avoid the boundary effects in image deconvolution \cite{AlmeidaFigueiredo:2013, MatakosRamaniFessler:2013} by masking the effect of boundary pixels in 2D images. The mask size is directly defined by the size of the convoluting kernel. However, by extending this method to video application, if only limited frames are available, e.g. $N=32$, choosing high-length kernel, e.g. $L=13$, will eliminate most frames' information by downgrading factor of $(N-L+1)/N$. Also, when applied to CS problem, the notion of random distribution of observation will be violated since the outermost pixels are masked and hence the restricted isometric property (RIP) in \cite{CandesTao:2005} will be no longer satisfied. In particular, we derive the unique quadratic solution from ADMM by deploying different BCs and name the algorithm: BC-ADMM. In particular we suggest using anti-reflective BC to avoid reconstruction artefacts. 
}

Our preliminary approach to the CVS problem was published in \cite{HosseiniPlataniotis:2012}, where HO accuracy filter was used to differentiate only temporal resolution in the video embedded with zero BC to create differentiation matrix in conjunction with wavelet dictionaries to make the sparse tensorial decomposition invertible in $\ell_1$ regularizer. However, the method suffered from bias artifacts in recovering boundary frames. This article here is mainly different from \cite{HosseiniPlataniotis:2012} in several aspects by the following summarized contributions,
{
\begin{enumerate}
\item We have Introduced for the first time a HO accuracy differential FIR filter which can be effectively used within the TV regularization framework. The introduced solution offers high accuracy in approximating the $1^{\text{st}}$ order derivative in TV over $[-1,1]$ filter, is capable of encoding high range of frequency components in an underlying signal and affords high precision restoration results.
\item An implicit convolution matrix is designed where the  convoluting kernels are HO accuracy differential FIR filters embedded with various BCs, e.g. AR-BC. The designed matrices are used to decompose spatial-temporal features in video for TV encoding by means of tensorial representation.
\item The designed space-time TV is used as a new regularizer in (\ref{eq4_2}) to recover under-sampled frames in CVS application. The solution to (\ref{eq4_2}) is driven by BC-ADMM which accommodates tensorial decomposition for solving a quadratic minimization problem with capability of handling different BCs.
\item The quality of reconstructed videos are much accurate compared to regularization via $[-1,1]$ filter deployed by several competitive algorithms. In particular, we have addressed the deficiency problem of texture/geometric loss in TV where the proposed TV regularizer is capable of restoring almost all textures in each individual frame. 
\end{enumerate} }

{
The remainder of this article is organized as follows. Section \ref{HO_vs_conventional} studies the major role of HO accuracy differentiation in a signal processing context and propose a new TV model. Section \ref{sec_TV_encoding} elaborates on the encoding scheme of the proposed TV into tensorial decomposition format. CVS enters to the picture in Section \ref{Sec:ProposedCVSFramework} by adopting the proposed TV regularizer and derives the minimization solution to BC-ADMM as well as providing formal convergence analysis. The experimental validations are demonstrated in Section \ref{Sec:Experiments} and finally Section \ref{Sec:Conclusion} concludes the article.}\vspace{-.1in}

{
\section{First-Order TV Revisited: The Role of High-Order Accuracy Differentiation}\label{HO_vs_conventional}
The main benefit of using high order (HO) accuracy operators for approximating the first order derivatives is to accurately transfer high frequency components in the applied signals during differentiation. To convey the importance of this point, the rest of this section elaborates on the utility of first-order numerical differentiation and its use in the TV regularization framework commonly applied in signal processing.}\vspace{-.1in}

\begin{figure}[h]
\centerline{
\subfigure[]{\includegraphics[width=0.2\textwidth]{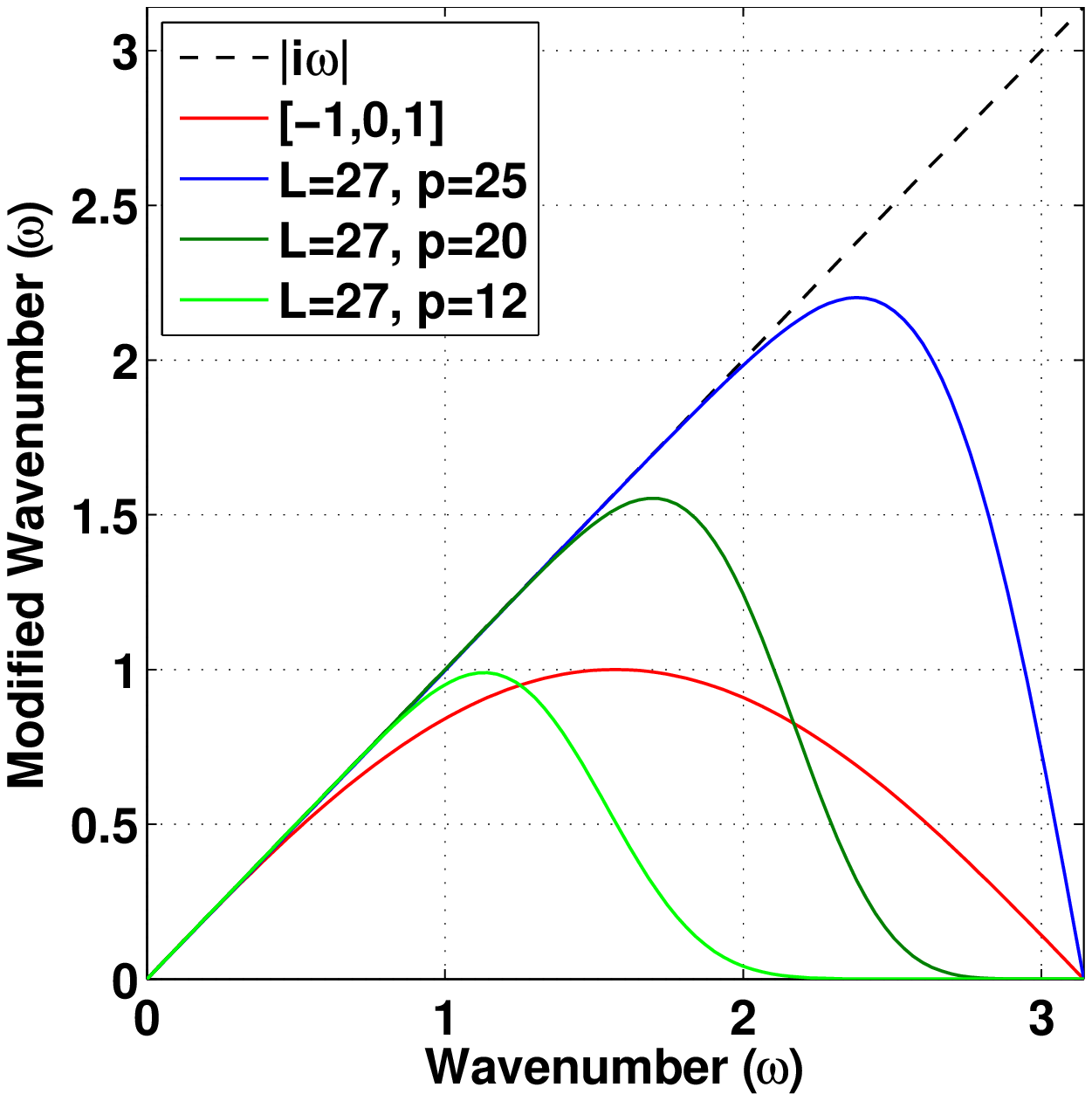}}
}\vspace{-.1in}
\caption{{Modified wave numbers (frequency response) of the first order numerical differentiation via different accuracy approximation}}
\label{fig_modified_wavenumber}
\end{figure}

{
\subsection{Features of HO accuracy differentiation}\label{HO_Accuracy_diff}
In general, the derivative of a continuous signal in time domain magnifies steep behaviours of the function and attenuates smooth variations. This is directly observable by taking the Fourier transform of the derivative function and yield the filter response as follows, $H_a(\omega)=\hat{f^{\prime}}(\omega)/\hat{f}(\omega)=i\omega$. The ideal derivative filter response $|i\omega|$, demonstrated in Figure \ref{fig_modified_wavenumber} in black-dashed-line, identifies the accurate amplification factor on different frequency range.}

{A common approach to approximate the ideal response in discrete domain is to reduce the model into digital filters either with $[-1,1]$ (forward), $[1,-1]$ (backward), or $[-1,0,1]/2$ (central) schemes \cite{Lele:1992, fornberg:1988, Li:2005, Khan:1999, Khan:2000}. However, all three cases exhibit the lowest accuracy in approximating the derivatives. To study this effect, the modified wavenumber (frequency response) of the central approach is demonstrated by the red plot in Figure \ref{fig_modified_wavenumber}. The discrete filter, compared to ideal response, \textit{deviates} frequencies beyond a certain range, e.g. $\omega>\pi/6$ in normalized frequency. Such deviation cause \textit{inappropriate transform} of high frequency components into lower magnitudes. Such high frequencies are not necessarily related to noise artefacts, but rather about the informative structures such as rapid transitions in a signal. It is important to note that this is different from the problem of aliasing in sampling theorem. Aliasing occurs due to insufficient (low resolution) sampling rate, while in the latter, we are provided with already-sampled function and required to estimate the derivatives from discrete samples. Therefore, the demand of HO accuracy filters is inevitable to accurately preserve high-frequency components in a discrete signal almost identical to ideal response and maintain accurate frequency transformation while approximating the derivatives.}

\begin{figure}[htp]
\centerline{
\subfigure{\includegraphics[width=0.4\textwidth]{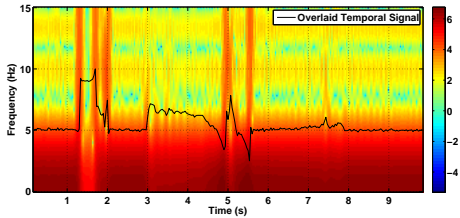}}
}\vspace{-.1in}
\caption{{Spectrogram and overlaid temporal video signal.}}
\label{fig_spectrogram_temporal_signal}
\end{figure}

{HO accuracy filters, in fact, has been vastly explored in the field of numerical methods for derivative approximation using Lagrangian interpolating polynomials \cite{fornberg:1988, Li:2005}, Taylor series expansion formulas \cite{Khan:1999, Khan:2000}, and Padde approximators \cite{Lele:1992}. Such methods are popular in numerical differentiation to achieve high accuracy rates as close as possible to the ideal response. Please, also refer to the references therein the papers \cite{Lele:1992, fornberg:1988, Li:2005, Khan:1999, Khan:2000} for further information. The aforementioned HO accuracy filters perform precisely when the signal meets two main conditions:
\begin{enumerate}
\item If the maximum frequency component of the signal is close to Nyquist rate (half of the sampling rate), then a high accuracy filter should be designed to avoid high frequency deviation in differentiating the signal.
\item The signal should be exact with noiseless measurements. If the signal is contaminated with noise, e.g. additive Gaussian, then the filter magnifies the noise-band/high frequency and hence cause error.
\end{enumerate}
Note the following: if the signal spectrum (band-limited) in continuous domain is close to the Nyquist rate, then the noise, e.g. additive white Gaussian, will alias into signal domain during the sampling. This is because the noise' frequency is beyond the Nyquist rate and related high-frequency artefacts will shift into the frequency domain that is representable by the sampled signal. For example, consider a temporal signal from HallMonitor video sketched in black line in Figure \ref{fig_spectrogram_temporal_signal}. The frame rate in this video is $30$ fps which is considered to be low for a class of natural video from a motion scene \cite{SukHwanApostolopoulosGamal:2005}. The corresponding Nyquist rate limits the maximum possible frequency in the temporal signal which should not exceed $15$ Hz. However, the temporal frequency components in HallMonitor video are beyond this limit. We show this feature by demonstrating the spectrogram analysis (time vs. frequency decomposition) in Figure \ref{fig_spectrogram_temporal_signal} for an overlaid signal in the graph which is fixed at spatial coordinates $(x,y)=(151, 174)$ of HallMonitor to study the temporal behaviours in the video. The signal is recorded for $10$ seconds, i.e. $300$ samples. To generate the spectrogram, we use Hamming window of length $10$ with $9$ overlapping samples in MATLAB. The colorbar is in logarithmic scale. Note the sharp signal transitions occupying the whole frequency range, i.e. zero to $15$ Hz. This class of temporal signal lies in the first condition discussed above. In this case, the high frequency band contains a mixture of signal and noise information. If the highest possible accuracy is set for differentiation, regardless of the noise effect, the approximated derivatives will much better preserve the sharp transitions in the signal. Therefore,  the signal derivative yields a better SNR condition, compared to the low accuracy estimation.}

\vspace{-.1in}

{In contrast, if the sampling rate is high enough (far away from maximum frequency component), then low-pass filtering the high-band is useful to attenuate the noise/aliasing effect. Smoothening techniques such as least square lanczos differentiators \cite{McDevitt:2012} known as Savitzky-Golay filters, low-pass filtering techniques combined with Padde approximation \cite{Kim:2010} and noise-robust (NR) approach \cite{Holoborodko:2008} are among existing methods to approximate the derivatives with controlled accuracy rate while low-passing the high-band frequency. In this study, we choose the noise-robust (NR) differentiation method from \cite{Holoborodko:2008} to encode frequency components in a signal to alter between two conditions discussed above. In particular, the user is provided with a tuning parameter $p$ to control the accuracy rate of the differentiation. The second parameter $q$, which has reverse relation to $p$, is also defined to control the smoothening factor of the low-pass filtering. The design of this filter is fully explained in detail in Appendix \ref{Sec:NRD}. The benefit of using NR filter is to satisfy both above conditions in two different scenarios: (a) in case of low sampling rate $p$ is chosen high to track all possible frequency components in the signal, (b) in case of high sampling rate, $p$ is chosen low. Consequently, this raise the $q$ value to avoid noise/aliasing effect.}

\vspace{-.1in}
{
\subsection{TV via HO accuracy filters}\label{TV_HO_Accuracy}
Bounded variation (BV) functionals refer to the functions with representable gradients distributed in a predefined domain set with finite Radon measure (vector-valued). Let $f\in L^1(\Omega)$ defined in an open subset $\Omega\in\mathbb{R}^n$. The TV of this function is defined by:
\begin{equation}\label{eq2_1}
\text{TV}(f,\Omega)=\int_{\Omega}{|\nabla f|dx},
\end{equation}
where the integration is evaluated on the function domain $\Omega$ \cite{AmbrosioFuscoPallara:2000}. TV encodes amount of oscillation in the function and is said to belongs to the BV space, i.e. $f\in\text{BV}(\Omega)$, if and only if its TV is finite that is $\text{TV}(f,\Omega)<\infty$. The celebrated work from Rudin et al. \cite{RudinOsherFatemi:1992} used BV as a regularization space in discrete model for image denoising problem. This model has been used over the last decades in variety of image restoration problems such as denoising, deblurring, inpainting, and CS. We refer the reader to \cite{chan2006total, Chambolle:2007, BrediesKunischPock:2010, BenningBruneBurgerMuller:2013} and references therein for recent advances on the topic.}


{
The TV norm of a function in BV space is equivalent to the surface area occupied by the functional in its domain which implies a direct relation between the geometry of the functional, e.g. an image in 2D space, and its TV measure \cite{chan2006total}. In discrete domain, failure to accurately approximate the derivatives for the TV measure leads to improper solution from a regularizer. Our claim in this paper is as follows: due to utilization of low accuracy filters in approximating the derivatives, either first order $[-1,1]$ or second order $[-1,2,-1]$, the conventional approaches fail to accurately encode the variational features (with high frequency component) in the signal and prone to errors in TV regularization problem. Variational features can be depicted by sharp transitions in the signal, i.e. discontinuities and edgal information. 
}

\begin{figure}[htp]
\centerline{
\subfigure[]{\includegraphics[width=0.3\textwidth]{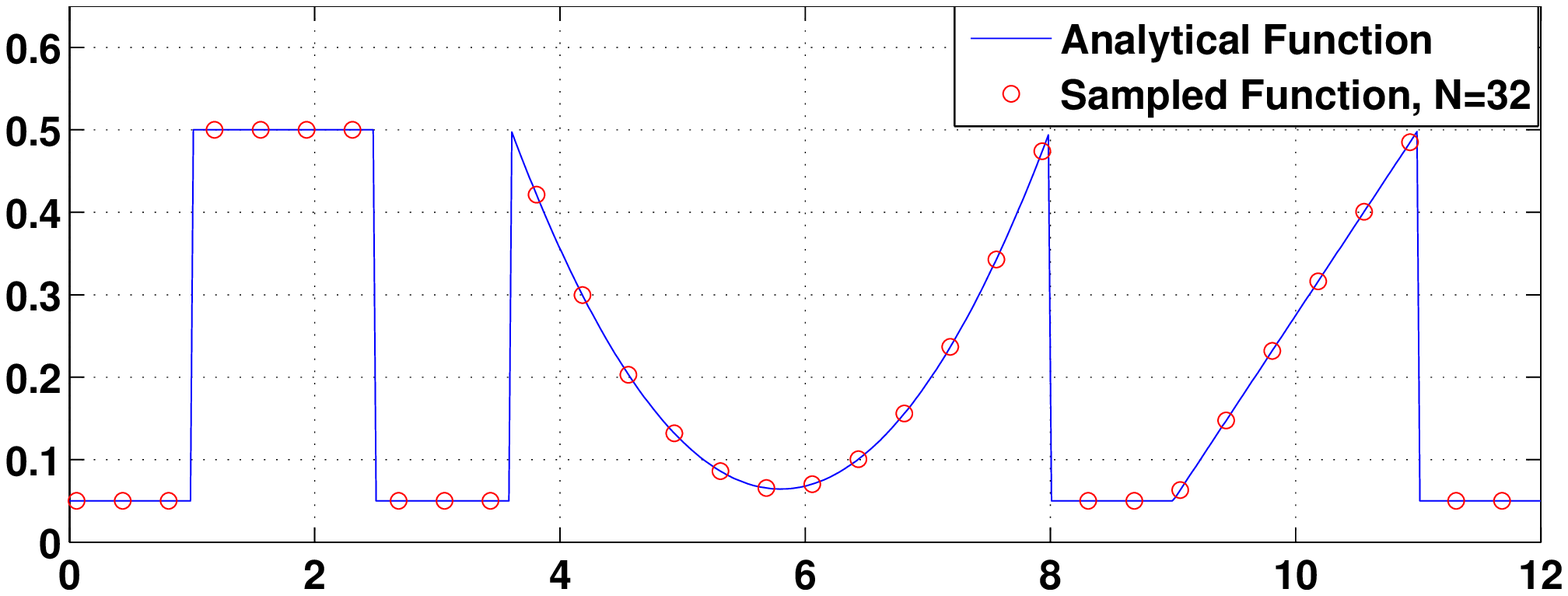}}
}
\centerline{
\subfigure[]{\includegraphics[width=0.3\textwidth]{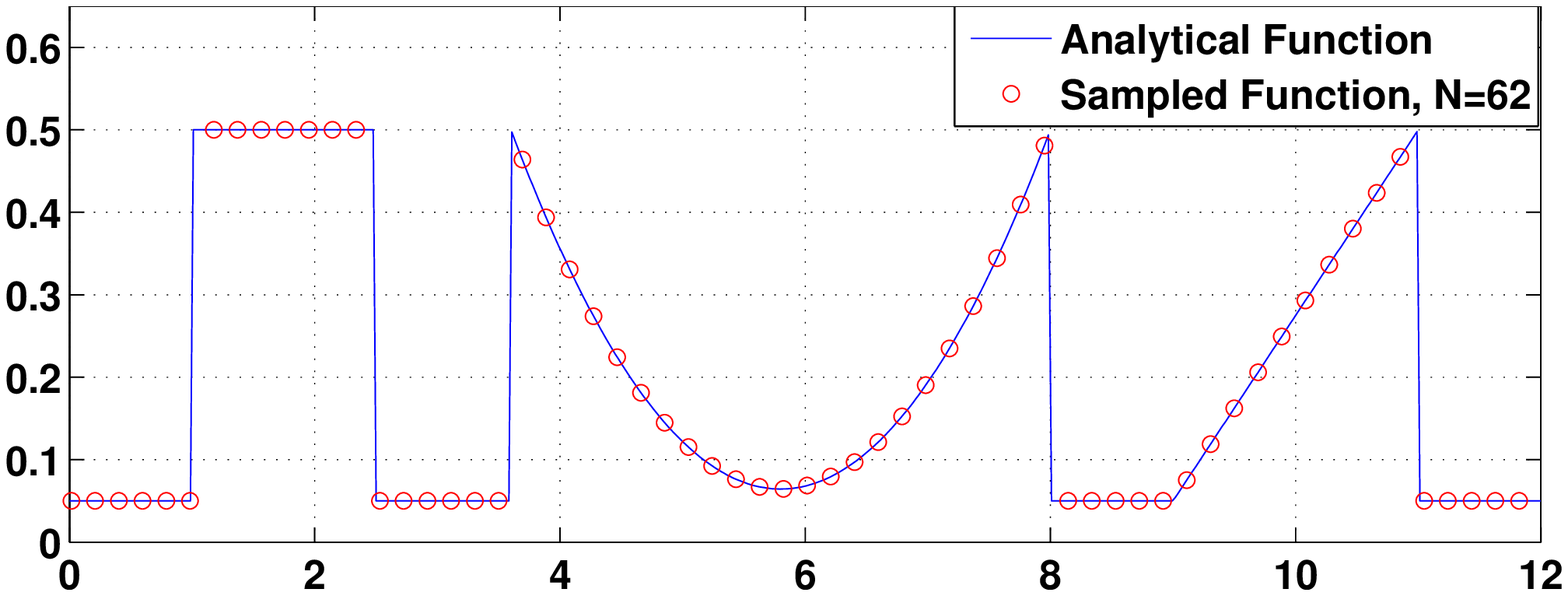}}
}\vspace{-.1in}
\caption{{Synthetic signal $u(x)$ with discrete samples: (a) $N=32$ (b) $N=62$}}
\label{fig_analytical_signal}
\end{figure}

{
We validate the above argument by numerically analyzing the effect of different derivative filters on approximating the TV norm value. We generate a synthetic signal similar to the one in \cite{ChanMarquinaMulet:2000} for testing benchmark. In particular, the signal contains several discontinuities with piecewise constant, linear, and second order polynomial variations to artificially simulate possible transitions as follows:
\begin{equation}\label{TV1}
u(x)=
\left\{
\begin{array}{lr}
\alpha, & x<a \\
\beta, & a\leq x<b \\
\alpha, & b\leq x<c \\
\frac{\beta-\mu}{(c-d)^2}(2x-c-d)^2+\mu, & c\leq x<d \\
\alpha, & d\leq x<e \\
\frac{\beta-\alpha}{f-e}(x-e)+\alpha, & e\leq x<f \\
\alpha, & f\leq x<g.
\end{array}
\right.
\end{equation}
Here, parameters are set to $a=1$, $b=2.5$, $c=3.6$, $d=8$, $e=9$, $f=11$, $e=12$ and $\{\alpha,\beta,\mu\}=\{0.05, 0.5, 0.09\}$ to separate transitions in the signal. The plot of synthetic signal $u(x)$ in (\ref{TV1}) is demonstrated in Figure \ref{fig_analytical_signal}. The calculated norm of the synthetic signal in (\ref{TV1}) using (\ref{eq2_1}) yields $TV(u)=\int^{\infty}_{-\infty}|u^{\prime}(x)|dx=8\beta-6\alpha-2\mu$.}

\begin{figure*}[htp]
\centerline{
\subfigure[]{\includegraphics[height=0.2\textwidth]{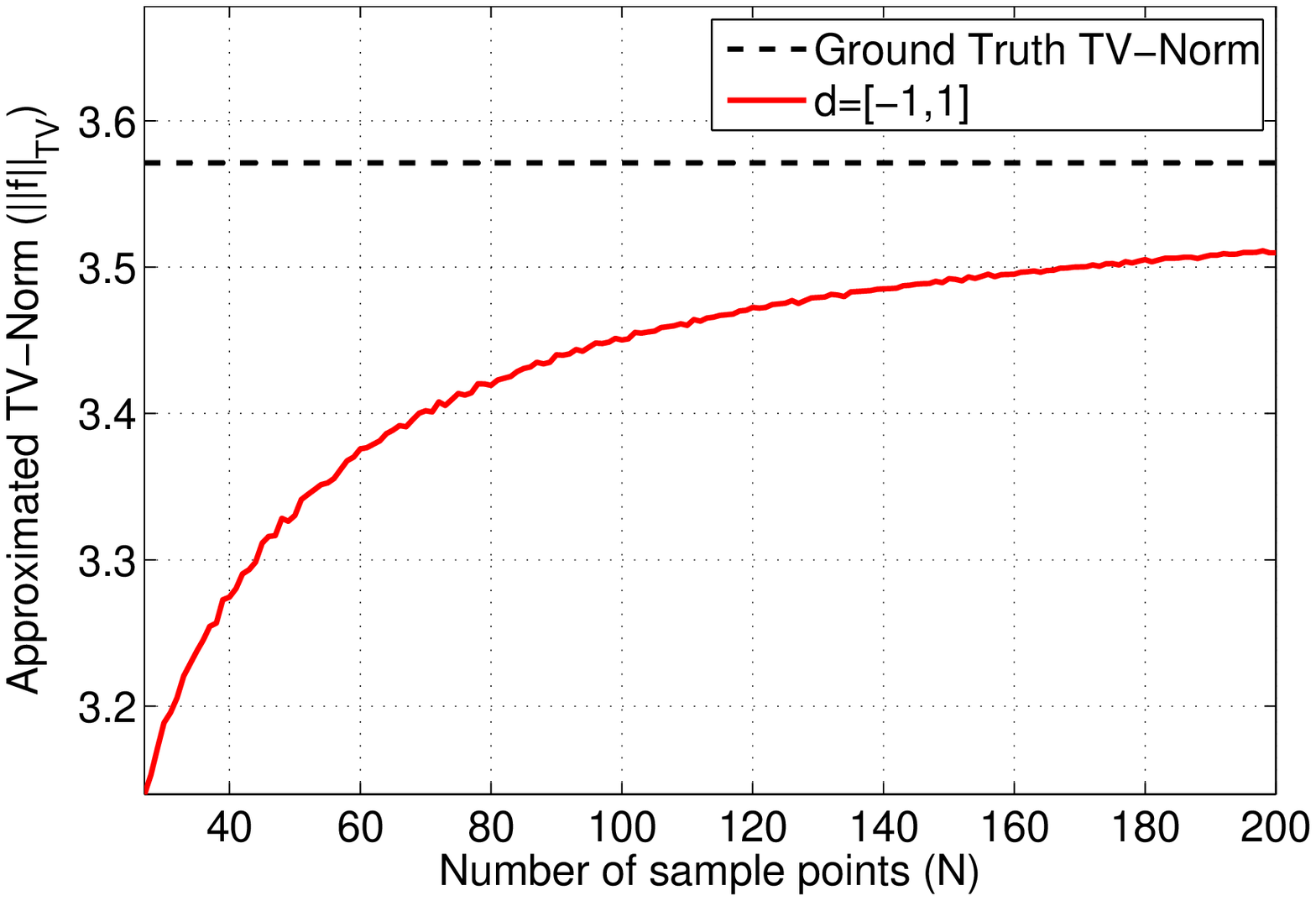}}
\subfigure[]{\includegraphics[height=0.2\textwidth]{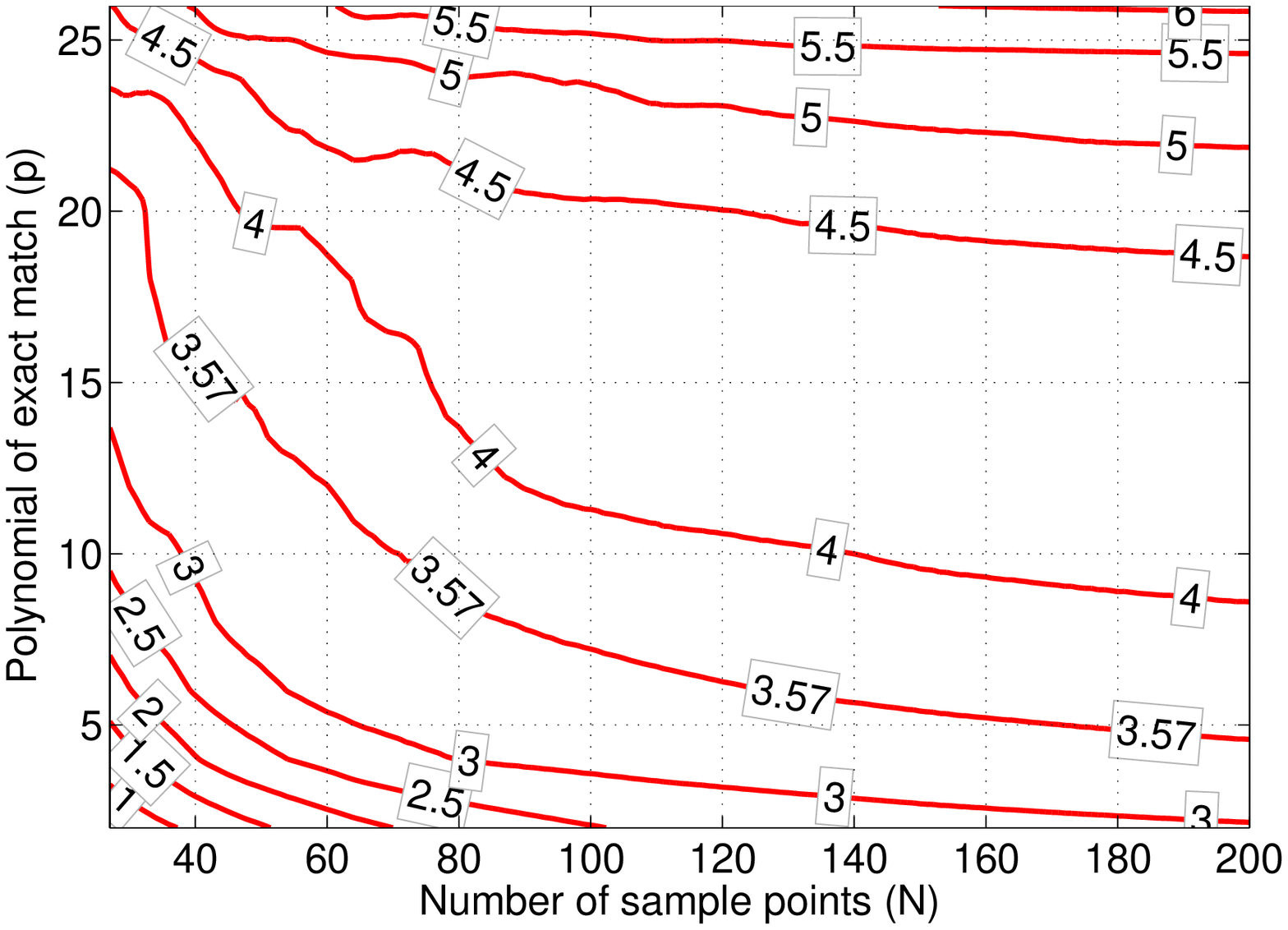}}
\subfigure[]{\includegraphics[height=0.205\textwidth]{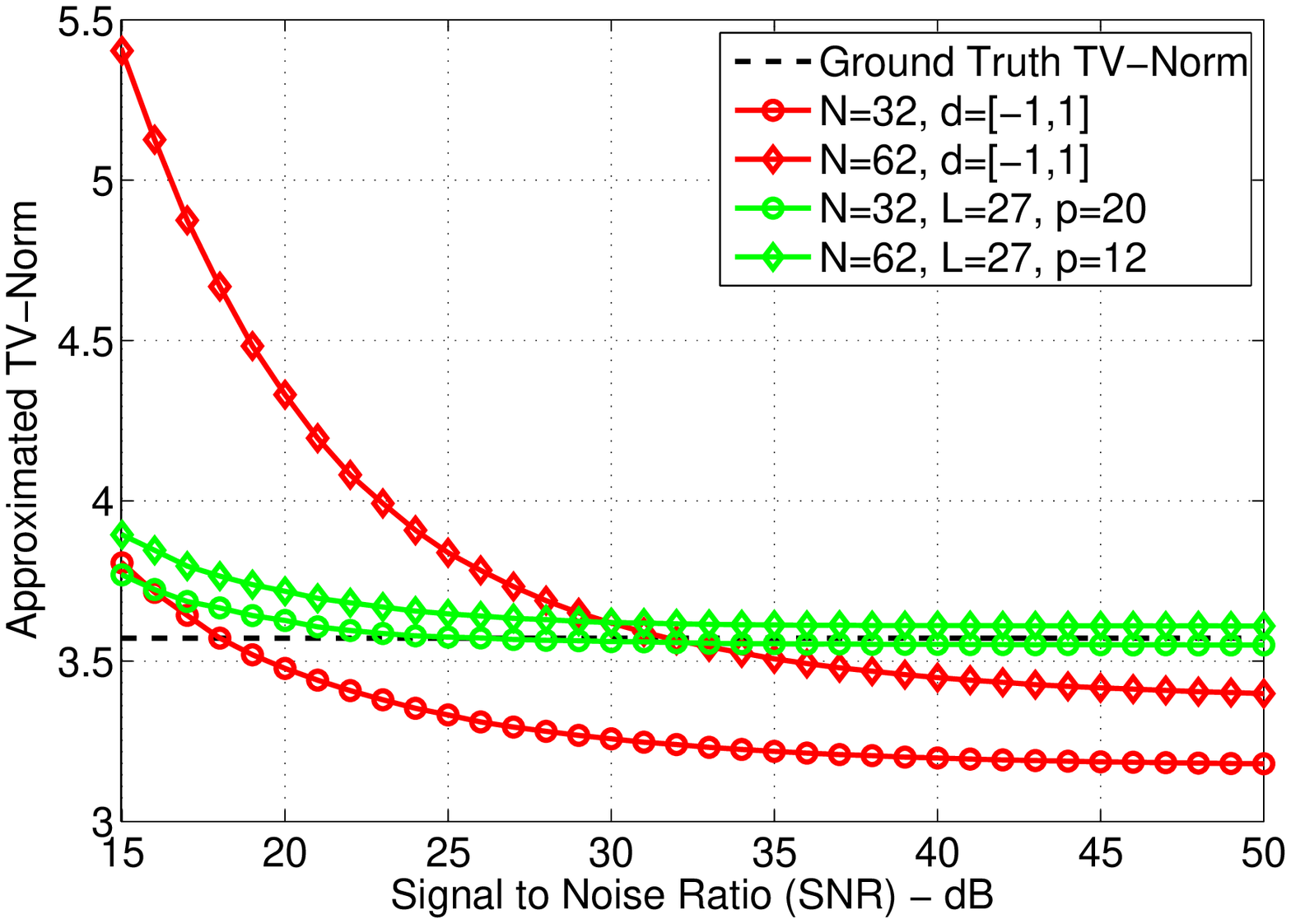}}
}\vspace{-.1in}
\caption{{(a) Approximated TV norm for different sampling points via conventional differentiation $[-1,1]$;} {(b) Approximated TV norm via different HO accuracy filters;} {(c) Noise sensitivity analysis of HO accuracy differentiation vs $[-1,1]$.}}
\label{fig_analytical_signal_HO_vs_L3_TV_performance}
\end{figure*}

{
Note in continuous domain, Dirac delta functions are generated in differentiating the signal. This is considered in calculating the TV norm, where the ground truth value for the aforementioned parameters is $TV(u)=3.5712$. In the following we generate discrete samples from the synthetic signal $u(x)$ using different sampling rates from $27$ to $200$ sample points. Increasing this rate consequently increases the Nyquist resolution and hence less aliasing will occurs for estimating the sharp edges in discrete signal. The minimum sample points is set to $27$. This is because the highest length NR filter is set to $L=27$. The discrete approximation of the TV norm in \cite{RudinOsherFatemi:1992} is evaluated for the generated discrete signals using: (a) conventional differentiation scheme ${\bf d} = [-1,1]$, and (b) HO accuracy differentiation using NR method introduced in Appendix \ref{Sec:NRD} for $L=27$ with different accuracy  factors $p=\{2,4,\hdots,26\}$. As a result, the conventional scheme deviates high frequency components in the signal during sharp transitions. This is the main reason why the approximated TV norm fails to accurately match with the ground truth value demonstrated in Figure \ref{fig_analytical_signal_HO_vs_L3_TV_performance}.a. We resolve this issue by deploying the proposed HO accuracy filter for TV approximation. Figure \ref{fig_analytical_signal_HO_vs_L3_TV_performance}.b demonstrates the phase transition of approximated TV values for different sampling $N$ and accuracy $p$ rates. The contour level sets of the approximated TV norm are shown. The empirical results suggest the accuracy parameter $p$ in HO accuracy filter should be set high when the Nyquist rate is low. This is to achieve close match between approximated and groundtruth values in order to accurately transfer high frequency components of the signal transitions in TV approximation. We further evaluate the performance of the HO accuracy differentiation via noise-sensitivity analysis. Two different sampling rates $N=32$ and $N=62$ are chosen to analyze different sampling rates. Both discrete signals are overlaid over the continuous signal in Figure \ref{fig_analytical_signal}.a and \ref{fig_analytical_signal}.b, respectively. The result of approximated TV through different noise levels is shown in Figure \ref{fig_analytical_signal_HO_vs_L3_TV_performance}.c. The results elaborate on the performance of the applied HO accuracy filter through different SNR values that provides robust and accurate approximation compared to conventional scheme. The type of simulated noise here is white Gaussian.}

{
To conclude this section, the issue of frequency deviation via conventional approach has been raised and compared to HO accuracy differentiation filters. The tuning of the NR filter parameters are discussed under different Nyquist rates and found that high accuracy parameters should be set when the maximum frequency component in a signal is close to the Nyquist rate. The direct impact of this phenomenon is further analyzed for TV problem where HO accuracy NR filter were able to accurately approximate the TV norm in (\ref{eq2_1}) compared to the conventional differentiator.}

\vspace{-.1in}
{
\section{Proposed TV Encoding Scheme}\label{sec_TV_encoding}
In this section we propose a numerical scheme to revisit the problem of TV regularization based on the benefits discussed in previous Section \ref{HO_vs_conventional}. We extend the TV encoder to space-time signals in multidimensional format, e.g. video (3D). In fact, the history of space-time TV backs to seminal works in \cite{WeickertJoachimSchnChristoph:2001, PapenbergBruhnBroxDidasWeickert:2006, werlberger2009anisotropic, ZachPockBischof:2007} for optical-flow estimation in video. Applications in deconvolution, denoising and CS problems can be also found in \cite{ShuAhuja:2011, ChanKhoshabehGibsonGillNguyen:2011}. We use similar concept to extend our results for space-time TV regularization.\vspace{-.1in}

\subsection{Numerical Design: Integrating HO accuracy Differentiation with Boundary Conditions}\label{sec_numerical_design}
The design of numerical difference scheme is presented in this section. As mentioned in Section \ref{HO_Accuracy_diff}, we adopt the HO of accuracy filters from \cite{Holoborodko:2008} to generate a digital convoluting kernel for differentiation. The complete design of this filter is presented in Appendix \ref{Sec:NRD}. The idea of integrating the so-called boundary condition (BC) to HO accuracy differentiation is also taken from the field of image deconvolution \cite{NgChanTang:1999,Capizzano:2003}. Similar technique is used to approximate the derivatives close to the signal boundaries. Handling BC is to appropriately preserve the derivative continuity in differentiation. This is of special importance in differentiating with HO accuracy filter where the length of the kernel is usually high.}

\begin{figure}[htp]
\centerline{
\subfigure{\includegraphics[width=0.3\textwidth]{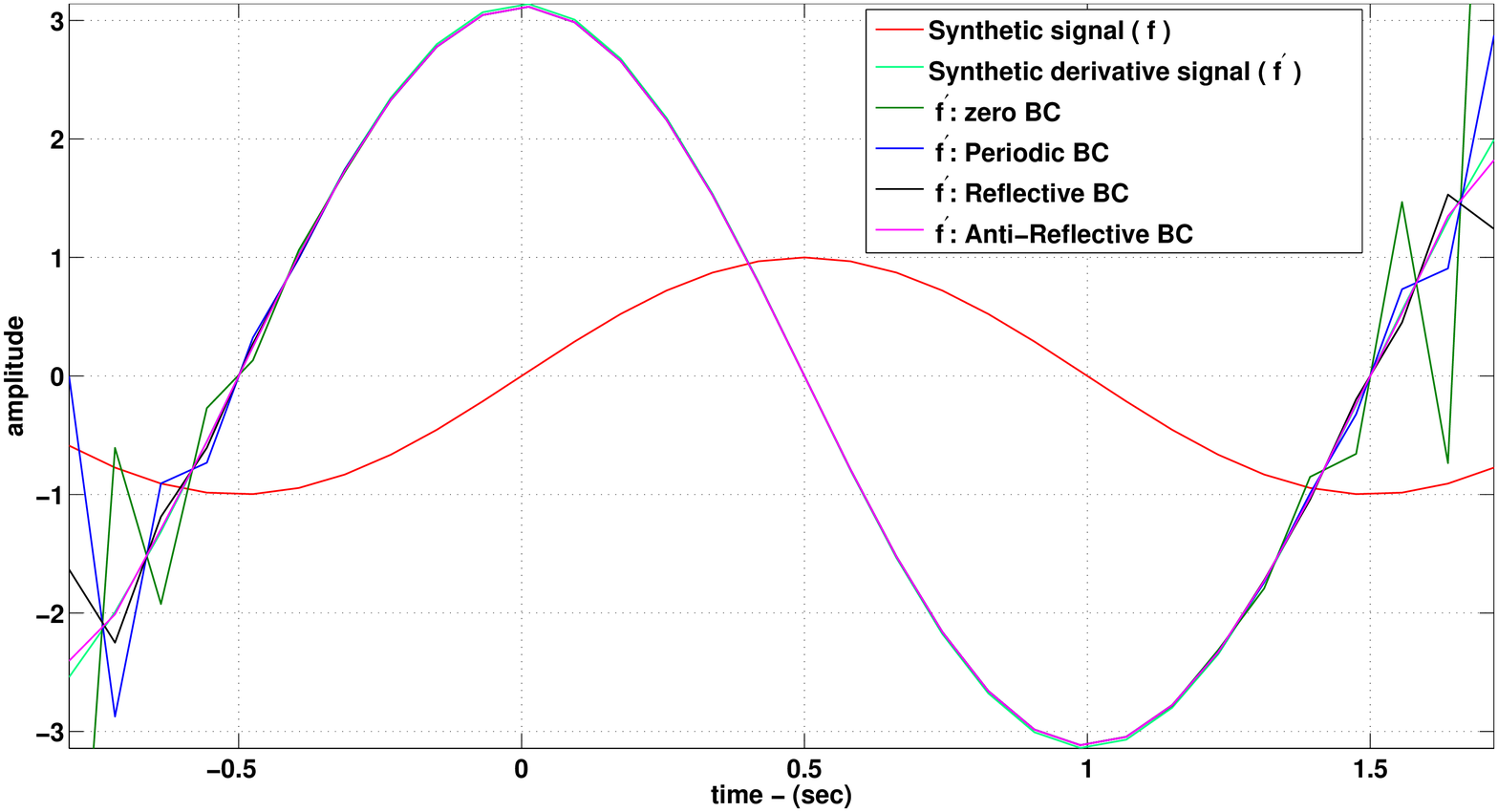}}
}\vspace{-.1in}
\caption{{Performance of different BCs on approximating sinusoidal signal derivatives}}
\label{fig_synthetic_signal_diff_BCs}
\end{figure}

In particular, we design a matrix ${\bf D}\in\mathbb{R}^{n\times n}$ which can operate on a discrete vector valued function $f\in\mathbb{R}^{n}$ and approximates the corresponding derivatives, i.e. $f^{\prime}\approx {\bf D} f$. {The derivative matrix ${\bf D}$ contains two key attributes: (a) adopts HO accuracy difference filters, and (b) is capable of embedding different BCs.} Consider the derivative kernels, designed in Appendix \ref{Sec:NRD}, and denote by
\begin{equation}\label{TD17}
{\bf d}=\left[d_{\frac{L-1}{2}},\hdots,d_1,0,-d_1,\hdots,-d_{\frac{L-1}{2}}\right]/T.
\end{equation}
Here, the kernel is assumed to be an anti-symmetric digital FIR filter with length $L$ (odd). The numerical differentiation using the kernel in (\ref{TD7}) can now be written as a convolution problem, i.e.
\begin{equation}\label{TD19}
{f}^{\prime}_j = \left({\bf d}\ast{\bf f}\right)(j)=\sum^{\frac{L-1}{2}}_{\ell=-\frac{L-1}{2}}{d_{\ell}\cdot f_{j-\ell}},
\end{equation}
where, $j\in\{0,1,\hdots,n-1\}$ corresponds to the vector index of the signal $f_j$. Here, the filter length should not exceed the number of the signal vector, i.e. $L\leq n$. The equivalent matrix-vector format of the Equation (\ref{TD19}) is
\begin{equation}\label{TD20}
{f}^{\prime}=
\left[
\begin{array}{c:c:c}
{\bf B_L}&{\bf D_T}&{\bf B_R}
\end{array}
\right]
\left[
\begin{array}{c}
f_{\bf L} \\
{f} \\
f_{\bf R}
\end{array}
\right],
\end{equation}
where ${\bf B}_L$, ${\bf D}_T$ and ${\bf B}_R$ are the Toeplitz matrices in (\ref{TD21}),
\begin{figure*}[htp]
\normalsize
\begin{equation}\label{TD21}
\begin{array}{l}
\hspace{.3cm}\overbrace{\hspace{2.75cm}}^{{\bf B_L}}
\hspace{.14cm}
\overbrace{\hspace{7.15cm}}^{{\bf D_T}}
\hspace{.1cm}
\overbrace{\hspace{3.25cm}}^{{\bf B_R}}
\\
\left[
{\begin{array}{ccc:ccccccc:ccc}
d_{\frac{L-1}{2}}&\hdots&d_1&0&\hdots&-d_{\frac{L-1}{2}}	&&&&&&\\
& \ddots&\vdots&\vdots&\ddots&\vdots&\ddots&&&&&\\
&& d_{\frac{L-1}{2}}&\vdots&\ddots&\ddots&\ddots&-d_{\frac{L-1}{2}}&&&&\\
&&& d_{\frac{L-1}{2}}&&&0&&&&&\\
&&&&\ddots&\ddots&\ddots&\ddots&\ddots&&&\\
&&&& &\ddots&0&\ddots&&-d_{\frac{L-1}{2}} &\\
&&&&&d_{\frac{L-1}{2}}&&\ddots&\ddots&	&-d_{\frac{L-1}{2}}&&\\
&&&&&&\ddots&&\ddots&\ddots&\vdots&\ddots&\\
&&&&&&& d_{\frac{L-1}{2} }&&0&-d_1&\hdots&-d_{\frac{L-1}{2}}
\end{array}}
\right]
\end{array}
\end{equation}
\vspace{-.2in}
\end{figure*}
and
\begin{equation}\label{TD22}
{f_{\bf L}}=\left[
\begin{array}{c}
f_{-\frac{L-1}{2}}\\
\vdots \\
f_{-1}
\end{array}
\right],
{f}=\left[
\begin{array}{c}
f_{0}\\
\vdots \\
f_{n-1}
\end{array}
\right],
{f_{\bf R}}=\left[
\begin{array}{c}
f_{n}\\
\vdots \\
f_{n+\frac{L-3}{2}}
\end{array}
\right].
\end{equation}
Here $f_{\bf L}$ and $f_{\bf R}$ are considered outside the field of view (left and right) for the original signal $f\in\mathbb{R}^{n}$, respectively. We used the same notation in \cite{NgChanTang:1999} to define appropriate BCs. The BCs can either defined a priori or from the information inside the field of view. In general, the signals contain no prior information from the boundaries and in order to perform the differentiation, the BCs are considered to be driven within the signal itself. So, $f_{\bf L}$ and $f_{\bf R}$ in (\ref{TD22}) are replaced by
\begin{equation}\label{TD23}
\begin{array}{l}
f_{\bf L}={\bf S_L}{f} \\
f_{\bf R}={\bf S_R}{f},
\end{array}
\end{equation}
where ${\bf S_L}$ and ${\bf S_R}$ define the linear transforms from true signal domain $f$ into the boundary values. Finally, by substituting (\ref{TD23}) in (\ref{TD20}), the following equation determines the derivatives with an appropriate boundary assignment
\begin{equation}\label{TD24}
{f}^{\prime}={\bf D}{f}=\left({\bf B_L}{\bf S_L}+{\bf D_T}+{\bf B_R}{\bf S_R}\right){f}
\end{equation}
where ${\bf D}$ contains appropriate BC for difference calculation. Dependent on the signal application, here, four types of BCs are summarized for the convolution purpose.
\begin{enumerate}
\item{Zero (Z)-BC (Dirichlet)}: the outside boundary of the signals is assumed to be zero, i.e. ${\bf S_L=0}$ and ${\bf S_R=0}$. So, ${f_{\bf L}=0}$ and ${f_{\bf R}=0}$ and the difference values from (\ref{TD24}) are
\begin{equation}\label{TD25}
{f}^{\prime}={\bf D_T}{f}.
\end{equation}
\item{Periodic (P)-BC}: signal is duplicated in periodic manner \cite{NgChanTang:1999}. For periodicity, $f_{\bf L}$ and $f_{\bf R}$ should satisfy
\begin{equation}\label{TD26}
f_{\bf L}=\left[
\begin{array}{c}
f_{n-\frac{L-1}{2}}\\
\vdots \\
f_{n-1}
\end{array}
\right],~~~f_{\bf R}=\left[
\begin{array}{c}
f_{0}\\
\vdots \\
f_{\frac{L-1}{2}-1}
\end{array}
\right].
\end{equation}
By substituting (\ref{TD26}) in (\ref{TD23}), the corresponding linear transforms are driven by
\begin{equation}\label{TD27}
\begin{array}{l}
{\bf S_L}=
\left[\begin{array}{cc}
{\bf 0} & {\bf I}
\end{array}\right]_{\frac{L-1}{2} \times n} \\
{\bf S_R}=\left[\begin{array}{cc}
{\bf I} & {\bf 0}
\end{array}\right]_{\frac{L-1}{2} \times n}
\end{array},
\end{equation}
where ${\bf I}$ is an identity matrix with $\frac{L-1}{2}$ dimension.
\item{Reflective (R)-BC}: the boundary values are assumed to be the reflected values from the true signal $f$ from the both ends \cite{NgChanTang:1999}, i.e.
\begin{equation}\label{TD28}
f_{\bf L}=\left[
\begin{array}{c}
f_{\frac{L-1}{2}-1+s}\\
\vdots \\
f_{s}
\end{array}
\right],~~~{f_{\bf L}}=\left[
\begin{array}{c}
f_{n-1-s}\\
\vdots \\
f_{n-\frac{L-1}{2}-s}
\end{array}
\right],
\end{equation}
where $s$ is the shifted index. The shift is some times useful to avoid sample repetition at the both ends, so $s\in\{0,1\}$. The corresponding linear transforms are provided by
\begin{equation}\label{TD29}
\begin{array}{l}
{\bf S_L}=
\left[\begin{array}{ccc}
{\bf 0} & {\bf I} & {\bf 0}_s
\end{array}\right]{\bf J} \\
{\bf S_R}=\left[\begin{array}{ccc}
{\bf 0}_s & {\bf I} & {\bf 0}
\end{array}\right]{\bf J}
\end{array},
\end{equation}
where $\bf J$ is the reversal permutation matrix with $n$ dimension and ${\bf 0}_s$ is a vector of zeros with $\frac{L-1}{2}$ elements.
\item{Anti-Reflective (AR)-BC}: the signal is reflected in anti-symmetric fashion towards the signal boundary values $f_0$ and $f_{n-1}$ \cite{Capizzano:2003}, i.e.
\begin{equation}\label{TD30}
f_{\bf L}=\left[
\begin{array}{c}
2f_{0}-f_{\frac{L-1}{2}-1+s}\\
\vdots \\
2f_{0}-f_s
\end{array}
\right],~~~f_{\bf R}=\left[
\begin{array}{c}
2f_{n-1}-f_{n-1-s}\\
\vdots \\
2f_{n-1}-f_{n-\frac{L-1}{2}-s}
\end{array}
\right].
\end{equation}
Here $s\in\{0,1\}$ is the same shift index used in R-BC and the corresponding linear transform is defined by
\begin{equation}\label{TD31}
\begin{array}{l}
{\bf S_L}=
\left(2\left[\begin{array}{cc}
{\bf 0} & {\bf 1}
\end{array}\right]-
\left[\begin{array}{ccc}
{\bf 0} & {\bf I} & {\bf 0}_s
\end{array}\right]\right){\bf J}
\\
{\bf S_R}=
\left(2\left[\begin{array}{cc}
{\bf 1} & {\bf 0}
\end{array}\right]-
\left[\begin{array}{ccc}
{\bf 0}_s & {\bf I} & {\bf 0}
\end{array}\right]\right){\bf J}
\end{array},
\end{equation}
where $\bf J$ and ${\bf 0}_s$ are defined the same in R-BC. Please refer to \cite{NgChanTang:1999,Capizzano:2003, Ng:2004} for comprehensive review on the BC.
\end{enumerate}
{
Among the existing approaches, AR-BC is the only one preserves first-order derivative continuity compared to Z-BC, P-BC, and R-BC. This is critical since the signal's tangency at both ends is not zero and behaves arbitrarily. To elaborate more on this argument, we generate a synthetic sinusoidal signal $f(t)=\sin(\pi t)$ and its derivatives, i.e. $f^{\prime}=\pi\cos(\pi t)$ for empirical analysis. We define the domain $t\in[-0.8,1.8]$ to avoid periodicity and reflectivity at the boundaries (arbitrary behaviour). We descretize the signal at $N=32$ points and approximate the derivatives using HO accuracy filter from Appendix \ref{Sec:NRD} combined with different BCs discussed above. The tap length is $L=27$ and accuracy of $p=25$. The results are demonstrated in Figure \ref{fig_synthetic_signal_diff_BCs}. All three BCs, i.e. Z-BC, R-BC,  and P-BC, deviate the signal at the boundaries, while AR-BC preserves mostly continuous behaviour compared to the synthetic derivatives.}\vspace{-.1in}

\vspace{-.1in}

{In general, the selection of BC is dependent to the signal application for better estimation of the boundaries. For instance, in case of periodic signal, we can directly use P-BC which can be efficiently implemented by means of fast Fourier transform (FFT). However, we do not control the dynamics of the signal at the boundaries which can end with arbitrary behaviour. Therefore, using AR-BC is more suitable to avoid possible irregular deviations. We provide appropriate solutions  from literature for fast manipulation of different BCs in next section that can be utilized in ADMM algorithm.}\vspace{-.1in}

\subsection{Tensorial Feature Encoding}\label{sec_tensorial_tv}
\begin{figure}[htp]
\centerline{
\includegraphics[width=0.225\textwidth]{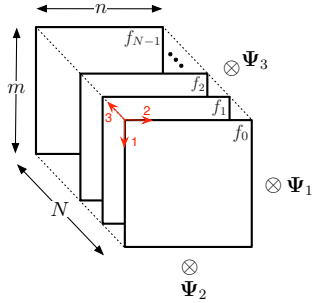}}\vspace{-.1in}
\caption{3D structure of $\bf F$ and relevant bases decomposition at each dimension}
\label{Fig:CubeStack}
\end{figure}

Video, in general, corresponds to consecutive snapshots from a scene to study the dynamics of the object movements. Frames are collected by a fixed sampling rate in time intervals known as the frame-rate and is measured by the number of collected frames per second (fps) \cite{Jahne:1993:SIP:528863}. Let us consider a sequence of $N$ video frames stacked during several time intervals, i.e.
\begin{equation}\label{eq1}
{\bf F} = \left[f_0,f_1,\hdots,f_{N-1}\right],
\end{equation}
where ${\bf F}\in\mathbb{R}^{m\times n\times N}$ is the tensorial representation of video in 3D. Every element in tensor  ${\bf F}$ can be accessed by three modes (dimension), i.e. ${\bf F}(x,y,t)$, where the first two modes $(x,y)$ correspond to the spatial coordinates and the third mode $t$ identifies the temporal variable at the incident time. {Space-time volume representation of the tensor ${\bf F}$ in (\ref{eq1}) is represented in Figure \ref{Fig:CubeStack}}, where each frame contains $m$-by-$n$ sample grids. {According to the Kronecker decomposition of multidimensional signals \cite{CaiafaCichocki:2012}, the tensor ${\bf F}$ can be decomposed in each direction by means of separate factor bases using the following vectorized tensor representation,}
\begin{equation}\label{TF2}
{{\bf c}}= \left({\bf\Psi}_3 \otimes {\bf\Psi}_2 \otimes {\bf\Psi}_1\right){{\bf f}},
\end{equation}
{where, ${\bf\Psi}_1\in\mathbb{R}^{m\times m}$, ${\bf\Psi}_2\in\mathbb{R}^{n\times n}$, and ${\bf\Psi}_3\in\mathbb{R}^{N\times N}$ are the factor bases. Also, ${\bf c}=\text{vec}\left({\bf C}\right)\in\mathbb{R}^{mnN}$ and ${\bf f}=\text{vec}\left({\bf F}\right)\in\mathbb{R}^{mnN}$ are the vectorized forms of mode-$1$ vector of tensors ${\bf C}$ and ${\bf F}$, respectively.} Here, the mode-$1$ vector of a tensor is obtained by fixing all modes in the tensor except the one on mode-$1$.

\vspace{-0in}

We adopt the vectorized tensor decomposition in (\ref{TF2}) to approximate the gradients ${\bf\nabla}{\bf F}$ along the spatio-temporal coordinates in (\ref{eq1}). The bases factors used for decomposition are the derivative matrices designed in Appendix \ref{Sec:NRD} and Section \ref{sec_numerical_design}. The following decomposition extracts the variational features calculated by the directional derivatives on each feature axis, i.e.
\begin{equation}\label{eq3}
\text{vec}\left({\bf\nabla}{\bf F}\right)=
\left[
\begin{array}{c}
\text{vec}\left({\partial{\bf F}}/{\partial x}\right) \\
\text{vec}\left({\partial{\bf F}}/{\partial y}\right) \\
\text{vec}\left({\partial{\bf F}}/{\partial t}\right)
\end{array}
\right] =
\left[
\begin{array}{c}
{\bf I}_N\otimes {\bf I}_n \otimes {\bf D}_m \\
{\bf I}_N\otimes {\bf D}_n \otimes {\bf I}_m \\
{\bf D}_N\otimes {\bf I}_n \otimes {\bf I}_m
\end{array}
\right]\text{vec}\left({\bf F}\right)
\end{equation}
where, ${\bf D}_m\in\mathbb{R}^{m\times m}$, ${\bf D}_n\in\mathbb{R}^{n\times n}$, and ${\bf D}_N\in\mathbb{R}^{N\times N}$ correspond to the designed derivative matrices decomposing the first two spatial and temporal coordinates, respectively. The estimated gradients in (\ref{eq3}) can now be substituted in (\ref{eq2_1}) to obtain the discrete estimation of the TV norm as follows:
\begin{eqnarray}\label{eq5}
\|{\bf F}\|_{\text{TV},\xi}=\|{\bf\nabla}{\bf F}\|_{\xi}=\sum\limits^{m,n,N}_{i,j,k=1}\sqrt[\xi]{\left|\frac{\partial{\bf F}}{\partial x}\right|^{\xi}_{ijk}+\left|\frac{\partial{\bf F}}{\partial y}\right|^{\xi}_{ijk}+\left|\frac{\partial{\bf F}}{\partial t}\right|^{\xi}_{ijk}}
\end{eqnarray}
where, $\xi=1$ and $\xi=2$ correspond to the anisotropic and isotropic TV, respectively.  Our preference in this paper is to use anisotropic norm to preserve spatio-temporal discontinuities. {This norm fully separates the gradient components and the regularization is done in $\ell_1$-ball which is more suitable for sparse recovery. The isotropic form however smoothens the edge/discontinuous information in the signal.} Furthermore, the summation of three gradient components (directional derivatives) in (\ref{eq5}) can be adjusted by three regularizing parameters associated with the sampling rates $T$ in each tensorial direction. The related sampling rates are incorporated to construct the derivative matrices in (\ref{eq3}) explained in Section \ref{sec_numerical_design}.\vspace{-.1in}

\section{Proposed CVS Framework}\label{Sec:ProposedCVSFramework}
In this section we propose a recovery framework for a CVS problem in order to exploit significant redundancies by means of the revised TV in \ref{sec_TV_encoding}. Multiple video frames are reconstructed from under-sampled measurements which are provided from individual frames at incident times. In particular, by deploying the proposed TV measure as a regularizer, associated with HO accuracy filters and proper BC, we expect to preserve sharp transitions and edgal information in the recovered video frames.\vspace{-.1in}

\subsection{Spatio-Temporal Compressive Sampling}\label{Sec:CVSsampling}
{
Single-pixel-camera is the first CI camera prototype introduced in \cite{DuarteDavenportBaraniuk:2008}. The samples were acquired using a digital programmable hardware for optical projection of the image scene. The main drawback of such conventional CS cameras is acquiring random measurements in a serial manner. This forces the application to maintain photography from fixed-scenes rather than scenes containing motion. Otherwise, the collected samples correspond to different time frames and cause confusion. This issue has gained more attention to develop a new hardware design to support single-frame-shot architecture in CI cameras. The revolutionary prototype of CMOS image sensor has been proposed in \cite{OikeGamal:2012, OikeGamal:2013} to perform single-shot CS imaging. This architecture supports high frame rate of $120$ fps to $1920$ fps corresponding to $1$ to $1/16$ compression ratio which mainly solves the issue of the time delay and makes it quite suitable for CVS applications. The CMOS image sensor is comprised with compress sensing multiplexer (CS-MUX) which applies random binary pattern $\phi_{ij}\in\{0,1\}$  from Bernoulli distribution at analog-to-digital-converter (ADC) level to outperform compressed measurements. The technology has also been integrated by Rademacher distribution to generate CS measurement matrix $\phi_{ij}\in\{-1,1\}$ in \cite{KaticHosseiniSchmidVandergheynstPierreLeblebici:2013, TomasMagli:2013}. Instant class of Rademacher distribution is the Walsh-Hadamard transform that can be easily implemented in CMOS hardware for CS application \cite{KhanWentzloff:2013}. For recent discussion and evaluation of the existing CS-CMOS imaging technologies and applications in video hardware, please refer to \cite{DadkhahDeenShirani:2013}. We restrict the design of our CVS framework to single-shot CS sampling from individual frames discussed above.} In general, the samples are collected by a linear transform where the image data is been observed. This transform identifies the connection between the acquisition and sampling domain. Existing transfer functions in CI devices are random Gaussian/Walsh-Hadamard multiplexers \cite{DuarteDavenportBaraniuk:2008,OikeGamal:2012, OikeGamal:2013}, Fourier bases in rapid MRI \cite{LustigDonohoPauly:2007}, circulant bases in coded apertures \cite{Marcia_Willett_2008}. Due to single-shot manner, the sampling operator is modelled by concatenating 2D sampling bases over the temporal dimension in a separable Kronecker format \cite{Marcia_Willett_2008, DuarteBaraniuk:2011, CaiafaCichocki:2012}, i.e. ${\bf b} = \overline{\bf\Phi}{\text{vec}\left({\bf F}\right)}=\left({\bf I}_N\otimes{\bf\Phi}_n\otimes{\bf\Phi}_m\right){\text{vec}\left({\bf F}\right)}$, where the ${\bf\Phi}_n\otimes{\bf\Phi}_m$ is the 2D-sampling bases discussed above. Furthermore, subsamples from individual frames are collected, i.e. ${\bf b}_{\bf\Omega} = \overline{\bf\Phi}_{\bf\Omega}{\text{vec}\left({\bf F}\right)}$, where ${\bf\Omega}=\{\Omega_0,\Omega_1,\hdots,\Omega_{N-1}\}$ corresponds to under-sampling indices assigned to $N$ individual frames and limited by total number of  available samples, i.e. $|{\bf\Omega}|<mnN$.\vspace{-.1in}

\subsection{Regularization Problem: $\text{TV}\ell_1/\ell_1$ Basis Pursuit De-Noising}\label{sec_BPDN}
We use the Proposed TV model in Section \ref{sec_TV_encoding} as the first regularizer to encodes the fine-scale information (rapid changes) pertinent to the spatio-temporal discontinuities/edgal information. Coarse-scale information, in contrary to the fine-scale, represents low-resolution components related to the image background with smooth variations \cite{LustigDonohoPauly:2007}. Therefore, we introduce the second regularizer for sparse linear transform using wavelet decomposition by means of vectorized tensor representation as follows,
\begin{equation}\label{eq2}
{\bf c} = \overline{\bf\Psi}{\text{vec}\left({\bf F}\right)}=\left({\bf I}_N\otimes{\bf\Psi}_n\otimes {\bf\Psi}_m\right){\text{vec}\left({\bf F}\right)},
\end{equation}
where, ${\bf I}_N\in\mathbb{R}^{N\times N}$ is an identity matrix, {and the decomposition bases ${\bf\Psi}_n\in\mathbb{R}^{n\times n}$ and ${\bf\Psi}_n\in\mathbb{R}^{n\times n}$} are performed on every individual frame at spatial axes, i.e. mode-$1$ and mode-$2$ of the tensor ${\bf F}$ {for sparse representation}. We augment these two regularizers to explore spatio-temporal redundancies for joint video frame reconstruction. We seek the solution $\bar{\bf F}$ from under-sampled measurements ${\bf b}_{\bf\Omega}$ in \ref{Sec:CVSsampling} via the following $\text{TV}\ell_1/\ell_1$ Basis Pursuit De-Noising (BPDN) problem
\begin{equation}\label{eq4_2}
\displaystyle \arg\min_{\bf F}\left\{\alpha\|{\bf F}\|_{\text{TV},\xi}+\beta\|\overline{\bf\Psi} {\text{vec}\left({\bf F}\right)}\|_1+{\gamma}\|\overline{\bf\Phi}_{\bf\Omega}{\text{vec}\left({\bf F}\right)}-{\bf b}_{\bf\Omega}\|^2_2\right\},
\end{equation}
where a strict equality is relaxed in (\ref{eq4_2}) to obtain a convex problem. A margin of error is allowed when the solution is either assumed to be compressible (non-exact sparsity) or the measurements are contaminated with Gaussian noise. The BPDN problem in (\ref{eq4_2}) is called compound minimization due to combination of both $\ell_1$ and TV norms for regularization \cite{GoldsteinOsher:2009}. The first compound minimization problem was first introduced in \cite{LustigDonohoPauly:2007} with many works reporting a variety of methods using ADMM approach \cite{YangZhangYin:2009}, split Bregman method \cite{GoldsteinOsher:2009}, and flexible conic formulations \cite{BeckerCandesGrant:2011}, where they used block-circulant sampling bases for image recovery based on the practical settings.\vspace{-.1in}

\subsection{Solution via BC-ADMM}\label{sec_bc_admmm}
We seek the unique solution to (\ref{eq4_2}) via ADMM which is a class of split variable technique \cite{GabayMercier:1976, citeulike:10118819} and introduce a new procedure, called BC-ADMM, to adopt various BCs that can be embedded in different CVS applications for video reconstruction such as CS-CMOS cameras. To solve the problem in (\ref{eq4_2}), the dual variables are decomposed; please refer to \cite{citeulike:10118819, AfonsoBioucas-DiasFigueiredo2011, AlmeidaFigueiredo:2013} for a comprehensive review. We use similar techniques in \cite{AfonsoBioucas-DiasFigueiredo2011, AlmeidaFigueiredo:2013, YangZhangYin:2009} to define auxiliary variables ${\bf y}_{(1)}$ and ${\bf y}_{(2)}$ to remove the transfer functions ${\bf\nabla}{\bf f}$ and $\overline{\bf\Psi}{\bf f}$ out of the non-differentiable norms $\|\cdot\|_{\xi}$ and $\|\cdot\|_1$, and define the third auxiliary variable ${\bf y}_{(3)}$ to separate the transfer function $\overline{\bf\Phi}{\bf f}$ from its sub-sampling scheme as follows,
\begin{eqnarray}\label{eq6}
\lefteqn{\bar{\bf F}=\arg\min_{\bf f}\left\{\alpha\|{\bf y}_{(1)}\|_{\xi}+\beta\|{\bf y}_{(2)}\|_1+\gamma\|{\bf P}{\bf y}_{(3)}-{\bf b}\|^2_2~~\right. } \\
& & \left. \text{s.t.}~{\bf y}_{(1)}={\bf H}_{(1)}{\bf f},{\bf y}_{(2)}={\bf H}_{(2)}{\bf f},{\bf y}_{(3)}={\bf H}_{(3)}{\bf f}\right\}\nonumber
\end{eqnarray}
where, ${\bf H}_{(1)}={\bf\nabla}$, ${\bf H}_{(2)}=\overline{\bf\Psi}$, ${\bf H}_{(3)}=\overline{\bf\Phi}$, and ${\bf f}=\text{vec}\left({\bf F}\right)\in\mathbb{R}^{mnN}$. Here, ${\bf P}:={\bf P}({\bf\Omega})$ is defined to be a function of ${\bf\Omega}$ referring to the selection operator for CS identified by the indices in ${\bf\Omega}$. From the method of multipliers, the augmented Lagrangian (AL) form of (\ref{eq6}) is formed by augmenting the quadratic penalty functions with Lagrangian multipliers to penalize the difference between auxiliary variables ${\bf y}_{(j)}$ and their corresponding transfer functions ${\bf H}_{(j)}{\bf f}$, yielding the following approximation model to (\ref{eq4_2}):
\begin{eqnarray}\label{eq7}
\lefteqn{\mathcal{L}_{\mathcal{A}}\left({\bf y}_{(1)}, {\bf y}_{(2)}, {\bf y}_{(3)}, {\bf f}, {\bf\lambda}_{(1)},{\bf\lambda}_{(2)},{\bf\lambda}_{(3)}\right) =} \\
& & \sum^{3}_{j=1} c_{(j)}g_{(j)}\left({\bf y}_{(j)}\right)+{\bf\lambda}^T_{(j)}{\bf r}_{(j)}+\mu_{(j)}/2\|{\bf r}_{(j)}\|^2_2, \nonumber
\end{eqnarray}
where,
\begin{equation}\label{eq8}
\begin{array}{lll}
c_{(1)}=\alpha & g_{(1)}(\cdot)=\|\cdot\|_{\xi} & {\bf r}_{(1)}={\bf H}_{(1)}{\bf f}-{\bf y}_{(1)} \\
c_{(2)}=\beta& g_{(2)}(\cdot)=\|\cdot\|_1& {\bf r}_{(2)}={\bf H}_{(2)}{\bf f}-{\bf y}_{(2)} \\
c_{(3)}=\gamma & g_{(3)}(\cdot)=\|{\bf P}{(\cdot)}-{\bf b}\|^2_2 & {\bf r}_{(3)}= {\bf H}_{(3)}{\bf f}-{\bf y}_{(3)},
\end{array}\nonumber
\end{equation}
and $\mu_{(j)}$ are the AL penalizing parameters. It is easy to notice that the AL variables in (\ref{eq7}) are separable with respect to ${\bf y}_{(1)}, {\bf y}_{(2)}$ and ${\bf y}_{(3)}$ for fixed Lagrangian multipliers ${\bf \lambda}_{(j)}$. Using the scaled dual variables ${\bf u}_{j}=(1/\mu_{(j)}){\bf \lambda}_{(j)}$ from \cite{citeulike:10118819}, the AL multipliers in (\ref{eq7}) can be expressed as follows
\begin{eqnarray}\label{eq9_1}
\lefteqn{\mathcal{L}_{\mathcal{A}}\left({\bf y}_{(1)}, {\bf y}_{(2)}, {\bf y}_{(3)}, {\bf f}, {\bf u}\right) =} \\
& & \sum^{3}_{j=1}{c_{(j)}g_{(j)}\left({\bf y}_{(j)}\right)+\mu_{(j)}/2\left[\|{\bf r}_{(j)}+{\bf u}_{(j)}\|^2_2-\|{\bf u}_{(j)}\|^2_2\right]}, \nonumber
\end{eqnarray}
The dual function is derived by calculating the infimum of the Lagrangian function in (\ref{eq9_1}), i.e.
\begin{equation}\label{eq9_2}
 q({\bf u})=\inf_{{\bf y}_{(1)}, {\bf y}_{(2)}, {\bf y}_{(3)}, {\bf f}}\mathcal{L}_{\mathcal{A}}\left({\bf y}_{(1)}, {\bf y}_{(2)}, {\bf y}_{(3)}, {\bf f}, {\bf u}\right),\nonumber
\end{equation}
and the following iterations are obtained for scaled version of ADMM approach using the gradient ascent method
\begin{IEEEeqnarray}{l} 
{\bf y}^{k+1}_{(j)}\longleftarrow\arg\min_{{\bf y}_{(j)}}\left\{c_{(j)}g_{(j)}\left({\bf y}_{(j)}\right)+\mu_{(j)}/2\|{\bf r}_{(j)}+{\bf u}_{(j)}\|^2_2\right\}\IEEEnosubnumber\label{eq10_1}\\
{\bf f}^{k+1}\longleftarrow\arg\min_{{\bf f}}\left\{
\sum^{3}_{j=1}{\mu_{(j)}/2\|{\bf r}_{(j)}+{\bf u}_{(j)}\|^2_2}\right\}\IEEEnosubnumber\label{eq10_2}\\
{\bf u}^{k+1}_{(j)}\longleftarrow{\bf u}^{k}_{(j)}+{\bf r}^{k+1}_{(j)},~j\in\{1,2,3\} \IEEEnosubnumber\label{eq10_3}
\end{IEEEeqnarray}
Problems in (\ref{eq10_1}-\ref{eq10_2}) consist of four sub-minimization problems with respect to ${\bf y}_{(1)}, {\bf y}_{(2)}, {\bf y}_{(3)}$ and ${\bf f}$. The so-called ${\bf y}_{(j)}$-sub minimization tasks including ${\bf f}$-sub problem can be carried out in parallel.

\subsubsection{${\bf y}_{(1)}, {\bf y}_{(2)}$-Sub Problems}\label{SubSec:Y12}
The first two sub-minimization problems from (\ref{eq10_1}) with respect to ${\bf y}_{(1)}$ and ${\bf y}_{(2)}$ are known as \textit{Moreau proximity operator} of $c_{(j)}g_{(j)}$ \cite{doi:10.1137/080724265} and the solution is given by
\begin{eqnarray}\label{eq11}
\lefteqn{{\bf y}^{k+1}_{(j)}\longleftarrow \text{prox}_{c_{(j)/\mu_{(j)}}}g_{(j)}\left({\bf s}_{(j)}\right)\equiv} \\
& & \arg\min_{{\bf y}_{(j)}}\left\{\frac{c_{(j)}}{\mu_{(j)}}g_{(j)}\left({\bf y}_{(j)}\right)+\frac{1}{2}\|{\bf y}_{(j)}-{\bf s}_{(j)}\|^2_2\right\}. \nonumber
\end{eqnarray}
where, ${\bf s}_{(j)}={\bf H}_{(j)}{\bf f}+{\bf u}_{(j)}$ for $j\in\{1,2\}$. The solution to the proximity operator in (\ref{eq11}) is provided by the splitting techniques in \cite{doi:10.1137/080724265, HaleYinZhang:2008} by means of soft-thresholding functions. The convergence of the aforementioned splitting methods are strongly guaranteed for sparse solutions with non-strict convexity assumptions compared to the existing methods. Alternative solutions can be found in \cite{HaleYinZhang:2008} and references therein. For the cases of anisotropic $g_{(j)}(\cdot)=\|\cdot\|_1$ and isotropic $g_{(j)}(\cdot)=\|\cdot\|_2$ norms, the proximity operators are defined by the soft and vector-soft thresholding, respectively, i.e.
\begin{IEEEeqnarray}{l}
\text{prox}^{\text{aniso}}_{T}g\left({\bf s}\right) = \text{soft}({\bf s},T) = \text{sgn}\left({\bf s}\right)\odot\max\left(|{\bf s}|-T,0\right) \IEEEnosubnumber\label{eq12_1} \\
\text{prox}^{\text{iso}}_{T}g\left({\bf s}\right) = \text{vec-soft}({\bf s},T) = {{\bf s}}/{\|{\bf s}\|_2}\cdot\max\left(\|{\bf s}\|_2-T,0\right) \IEEEnosubnumber\label{eq12_2}
\end{IEEEeqnarray}
where, $\odot$ stands for element-wise multiplications, the absolute function $|\cdot|$ in isotropic solution is element-wise operator, and finally with the convention $0/\|0\|_2=0$.

\subsubsection{${\bf y}_{(3)}$-Sub Problem}\label{SubSec:Y3}
Minimizations with respect to ${\bf y}_{(3)}$ is attained by finding the minimum argument in (\ref{eq10_1}) for $j=3$, 
\begin{eqnarray}\label{eq13}
\lefteqn{{\bf y}^{k+1}_{(3)}\longleftarrow\arg\min_{{\bf y}_{(3)}}\left\{c_{(3)}\|{\bf P}{\bf y}_{(3)}-{\bf b}\|^2_2+\right.} \\
& & ~~~~~~~~~~~~~~~~~~~~~~\left. \frac{\mu_{(3)}}{2}\|{\bf H}_{(3)}{\bf f}-{\bf y}_{(3)}+{\bf u}_{(3)}\|^2_2\right\}. \nonumber
\end{eqnarray}
which is quadratic and the closed-form solution is driven by the minimum least squares with respect to ${\bf y}_{(3)}$, i.e.
\begin{eqnarray}\label{eq14}
\lefteqn{{\bf y}^{*}_{(3)}=\left[2c_{(3)}{\bf P}^T{\bf P}+\mu_{(3)}{\bf I}\right]^{-1}} \\
&&~~~~~~~~~~~\left[\mu_{(3)}{\bf H}_{(3)}{\bf f}+\mu_{(3)}{\bf u}_{(3)}+2c_{(3)}{\bf P}^T{\bf b}\right]. \nonumber
\end{eqnarray}
{Here, ${\bf P}^T{\bf P}$ constitutes a diagonal operator since ${\bf P}$ is an underdetermined matrix who's rows are a subset of rows from an identity matrix selected at random. Therefore, ${\bf P}^T{\bf P}$ constitutes an square matrix, where the diagonal components $\textit{diag}({\bf P}^T{\bf P})$ are combination of ones and zeros in random distribution identifying the indices of under-sampling pattern.} So, ${\bf y}^{k+1}_{(3)}$ can be easily calculated.

\subsubsection{${\bf f}$-Sub Problem}\label{SubSec:f}
In the final step, the minimization with respect to ${\bf f}$ forms a quadratic minimization problem, i.e.
\begin{equation}\label{eq15}
{\bf f}^{k+1}\longleftarrow \arg\min_{{\bf f}}\left\{\sum^{3}_{j=1}{\mu_{(j)}/2\|{\bf H}_{(j)}{\bf f}-{\bf\zeta}_{(j)}\|^2_2}\right\},
\end{equation}
where ${\bf\zeta}_{(j)}={\bf y}_{(j)}-{\bf u}_{(j)}$. Calculating the gradient with respect to ${\bf f}$ and equating it with zero, the unique solution ${\bf f}^{*}$ is driven for iterative updating, i.e.
\begin{equation}\label{eq16}
{\bf f}^{*}=\left[\sum^{3}_{j=1}{\mu_{(j)}{\bf H}_{(j)}^T{\bf H}_{(j)}}\right]^{-1}\left[\sum^{3}_{j=1}{\mu_{(j)}{\bf H}_{(j)}^T{\bf\zeta}_{(j)}}\right]
\end{equation}

{As mentioned in Section \ref{Sec:CVSsampling}, the sampling basis is identified from the feature of the CI devices, e.g. CS-CMOS imaging and single-pixel-camera. Here, we assume two orthonormal bases for ${\bf H}_{(3)}=\overline{\bf\Phi}$ in our experiments. The first one is matrix formed by populating the entries with i.i.d. Gaussian variables followed by orthonormalization of the rows, and the second is Walsh-Hadamard matrix transform. Both transforms are proved to hold RIP condition in \cite{CandesTao:2005, HowardCalderbankSearle:2008}}. Also, the decomposition bases ${\bf H}_{(2)}=\overline{\bf\Psi}$ is orthonormal wavelet bases. So, the matrix to be inverted in (\ref{eq16}) can now be expressed as
\begin{equation}\label{eq17}
{\bf O} = \mu_1{\bf\nabla}^T{\bf\nabla}+(\mu_2+\mu_3){\bf I},
\end{equation}
where,
\begin{equation}\label{eq17_1}
{\bf\nabla}^T{\bf\nabla}=
{\bf I}_N\otimes{\bf I}_n\otimes{{\bf D}^T_m{\bf D}_m} +
{\bf I}_N\otimes{{\bf D}^T_n{\bf D}_n}\otimes{\bf I}_m +
{{\bf D}^T_N{\bf D}_N}\otimes{\bf I}_n\otimes{\bf I}_m.
\end{equation}
The inverse operator ${\bf O}^{-1}$ can be directly utilized if it is diagonalizable by means of a proper orthogonal basis. Here, {\bf D} forms differential operator associated with four possible BCs constituted with anti-symmetric derivative kernels. Consequently, ${\bf D}^T{\bf D}$ contains symmetric kernels. Without loss of generality, let us assume the following decomposition for diagonalizing three different matrices pertinent with the aforementioned directional derivatives,
\begin{equation}\label{eq18}
\begin{array}{lll}
{\bf D}^T_m{\bf D}_m & = & {\bf Q}_m{\bf \Lambda}_m{\bf Q}^T_m \\
{\bf D}^T_n{\bf D}_n & = & {\bf Q}_n{\bf \Lambda}_n{\bf Q}^T_n \\
{\bf D}^T_N{\bf D}_N & = & {\bf Q}_N{\bf \Lambda}_N{\bf Q}^T_N,
\end{array}
\end{equation}
where ${\bf \Lambda}$ is a diagonal eigenvalue matrix and ${\bf Q}$ is the orthogonal basis for decomposition related to different BCs. {The eigenvalues ${\bf \Lambda}$ are non-negative corresponding to the Fourier spectrum of the second derivatives i.e. $|(i\omega)^2|$ from rank-deficient matrix ${\bf D}{\bf D}^T$, where the rank deficiency is caused by different BCs.} The operator ${\bf O}$ in (\ref{eq17}) can now be decomposed by ${\bf O} = {\overline{\bf Q}}{\bf\Lambda}_0{\overline{\bf Q}}^T$, where
\begin{IEEEeqnarray}{l}
\overline{\bf Q} = ({\bf Q}_N\otimes{\bf Q}_n\otimes{\bf Q}_m), \label{eq19_1}\\
{\bf\Lambda}_0 =
(\mu_2+\mu_3){\bf I}_{mnN} +\mu_{(1)}\left({\bf I}_N\otimes{\bf I}_n\otimes{\bf \Lambda}_m + \right. \label{eq19_2}\\
\hspace{1.15in}\left.{\bf I}_N\otimes{\bf \Lambda}_n\otimes{\bf I}_m +{\bf \Lambda}_N\otimes{\bf I}_n\otimes{\bf I}_m\right).\nonumber
\end{IEEEeqnarray}
Here, ${\bf\Lambda}_0$ is strictly positive diagonal matrix and can be directly inverted. {For comprehensive treatment of computational complexities for different BCs please refer to \cite{NgChanTang:1999,Capizzano:2003, Ng:2004}. In particular, P-BC and R-BC can be diagonalized by FFT matrices and AR-BC via fast Sine transform (FST) all with computational complexity of orders of $\mathcal{O}(n\log n)$. Consequently, the inverse operator ${\bf O}^{-1}$ can be easily calculated to boost up the processing speed in computer.}

Following the update solutions from the sub-problems in Sections \ref{SubSec:Y12},\ref{SubSec:Y3}, and \ref{SubSec:f}, the pseudocode of the proposed BC-ADMM method is summarized in algorithm \ref{algo_BC_ADMM}. Sub-problem ${\bf f}^k$ is updated in line \ref{algo:f-sub} of the algorithm using (\ref{eq16}) which is the computational bottleneck of the Algorithm \ref{algo_BC_ADMM}, since it needs the inversion of the matrix ${\bf O}$ in (\ref{eq17}). This is utilized by the aforementioned BCs using the eigenvalue decomposition in (\ref{eq19_1}-\ref{eq19_2}). Sub-problems ${\bf y}^{k}_{(j)}$ are also updated in lines \ref{algo:y12-sub} and \ref{algo:y3-sub}. Finally, the algorithm proceeds to update the dual variables in line \ref{Alg:update}. The algorithm is terminated when the optimality condition is satisfied. This is when the relative error change (tolerance error) of updating ${\bf f}$-sub problem is small enough, i.e. $ \|{{\bf f}^{k+1}-{\bf f}^{k}}\|/{\|{\bf f}^{k}\|}\leq\epsilon$. {The convergence of BC-ADMM is addressed in Section \ref{sec_convergence}}. Fast convergence of the ADMM algorithm with the length $\rho\in(0,(\sqrt{5}+1)/2)$ was also guaranteed in \cite{GlowinskiTallec:1989} in the context of the variational inequality. Thus, in line \ref{Alg:update} of the Algorithm \ref{algo_BC_ADMM}, a steplength $\rho$ is allowed for updating the dual variables ${\bf u}^{k}_{(j)}$.\vspace{-.1in}

\begin{algorithm}[htp]
\SetKwData{Left}{left}\SetKwData{This}{this}\SetKwData{Up}{up} \SetKwFunction{Union}{Union}\SetKwFunction{FindCompress}{FindCompress}
\SetKwInOut{Input}{input}
\SetKwInOut{Output}{output}
\Input{Selection operator ${\bf P}$, set ${\bf H}_{(1)}={\bf\nabla}$, ${\bf H}_{(2)}=\bar{\bf \Psi}$ and ${\bf H}_{(3)}=\bar{\bf\Phi}$, model parameters $\{c_{(j)},\mu_{(j)}\}^3_{j=1}$, initialize $\{{\bf y}^0_{(j)}={\bf 0},{\bf u}^0_{(j)}={\bf 0}\}^3_{j=1}$, and set $k=0$.}
\Output{optimized ${\bf f}^{*}\approx{\bf f}^{k+1}$} \BlankLine
\Repeat{convergence condition is satisfied}
{\label{forins}
{${\bf\zeta}^{k}_{(j)}~\leftarrow{\bf y}^{k}_{(j)}-{\bf u}^{k}_{(j)}$ for $j\in\{1,2,3\}$} \\
{${\bf f}^{k+1}\leftarrow$ ${\overline{\bf Q}}{\bf\Lambda}^{-1}_0{\overline{\bf Q}^T}\left[\sum^{3}_{j=1}{\mu_{(j)}{\bf H}_{(j)}^T{\bf\zeta}^k_{(j)}}\right]$}\label{algo:f-sub} \\
\For{$j=1$ to $2$}
{
${\bf y}^{k+1}_{(j)}\leftarrow \text{prox}_{\frac{c_{(j)}}{\mu_{(j)}}}g_{(j)}\left({\bf H}_{(j)}{\bf f}^{k+1}+{\bf u}^{k}_{(j)}\right)$\label{algo:y12-sub}}
{${\bf y}^{k+1}_{(3)}\leftarrow\left[2c_{(3)}{\bf P}^T{\bf P}+\mu_{(3)}{\bf I}\right]^{-1}$\label{algo:y3-sub}\\
$\hspace{.7in}\left[\mu_{(3)}{\bf H}_{(3)}{\bf f}^{k+1}+\mu_{(3)}{\bf u}^{k+1}_{(3)}+2c_{(3)}{\bf P}^T{\bf b}\right]$}\\
\For{$j=1$ to $3$}
{
${\bf u}^{k+1}_{(j)}\leftarrow{\bf u}^{k}_{(j)}+\rho\left[{\bf H}_{(j)}{\bf f}^{k+1}-{\bf y}^{k+1}_{(j)}\right]$
}\label{Alg:update}
$k\leftarrow k+1$
}
\caption{{\bf BC-ADMM}}\label{algo_BC_ADMM}
\end{algorithm}\vspace{-.2in}

{
\subsection{Convergence Validation of BC-ADMM}\label{sec_convergence}
Generalized ADMM algorithm is a special case of the Douglas-Rachford splitting method used to find the zeros of the sum of two maximal monotone operators \cite{GlowinskiMarrocco:1975, GabayMercier:1976, EcksteinBertsekas:1992}. Consider the following unconstrained optimization problem
\begin{equation}\label{sec_convergence_eq1}
\min_{{\bf f}\in\mathbb{R}^d}~~~f_1({\bf f})+f_2({\bf G}{\bf f}),
\end{equation}
where ${\bf G}\in\mathbb{R}^{r\times d}$. $f_1\in\mathbb{R}^d\mapsto\bar{\mathbb{R}}$ and $f_2\in\mathbb{R}^r\mapsto\bar{\mathbb{R}}$ are closed proper convex functions which remains as an instant class of maximal monotone operators. If (\ref{sec_convergence_eq1}) has a solution ${\bf f}^{*}$, then the sequences $\{{\bf f}\}_{k}$ generated by the generalized ADMM converges to the solution ${\bf f}^{*}$ under the conditions of the Theorem $8$ in \cite{EcksteinBertsekas:1992} from Eckstein and Bertsekas. In particular, when ${\bf G}$ has full column rank the convergence of the outcome sequences from ADMM holds, please refer to \cite{EcksteinBertsekas:1992} for more information. Finding the zero of the sum of $J>2$ maximal monotone operators remains a broad topic in split variable technique. Seminal works on the convergence analysis of such problem can be found in \cite{Spingarn:1983, EcksteinSvaiter:2009, GoldfarbMa:2012}. Simplified version of the maximal monotone operators are convex feasibility problems. Figueiredo et al. \cite{FigueiredoBioucasDias:2010, AfonsoBioucas-DiasFigueiredo2011, AlmeidaFigueiredo:2013} extended the validity of convergence for genrelized ADMM by defining a linear mapping that maps an objective function with $J>2$ operators into (\ref{sec_convergence_eq1}), where the conditions of the Theorem in \cite{EcksteinBertsekas:1992} hold for the argument.}

\begin{figure}[htp]
\centerline{\includegraphics[width=0.3\textwidth]{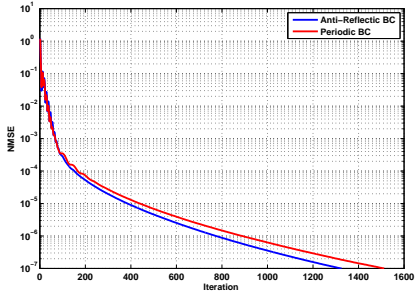}}\vspace{-.1in}
\caption{Evolution of updated sequence ${\bf f}_k$ for BC-ADMM.}
\label{Fig:convergence}
\end{figure}

{
The optimization problem presented by the authors fits exactly the variant model discussed in \cite{FigueiredoBioucasDias:2010}, Section III.B, by  mapping three ($J=3$) closed convex functions into (\ref{sec_convergence_eq1}) corresponding to the proposed $\text{TV}\ell_1/\ell_2$ BPDN problem in Section \ref{sec_BPDN} as follows,
\begin{equation}\label{sec_convergence_eq2}
\min_{{\bf f}\in\mathbb{R}^d}~~~\sum^{J=3}_{j=1}g_{(j)}\left({\bf H}_{(j)}{\bf f}\right),
\end{equation}
where $g_{(j)}$ are closed proper convex functions and ${\bf H}_{(j)}$ are matrices all defined in Section \ref{sec_bc_admmm}. Under the following assumption the minimization problem in (\ref{sec_convergence_eq2}) is equivalent to the problem in (\ref{sec_convergence_eq1}) ,
\begin{equation}\label{sec_convergence_eq3}
\begin{array}{ll}
f_1=0, & {\bf G}=\left[
\begin{array}{l}
{\bf H}_{(1)} \\
{\bf H}_{(2)} \\
{\bf H}_{(3)}
\end{array}
\right]\in\mathbb{R}^{r\times d}
\end{array},
\end{equation}
where, $r=5mnN$ and $d=mnN$ in this paper, and
\begin{equation}\label{sec_convergence_eq4}
f_2({\bf y})=\sum^{J=3}_{j=1}{g_{(j)}\left({\bf y}_{(j)}\right)}.
\end{equation}
Here, ${\bf y}_{(1)}\in\mathbb{R}^{3mnN}$, ${\bf y}_{(2)}\in\mathbb{R}^{mnN}$, and ${\bf y}_{(3)}\in\mathbb{R}^{mnN}$, where the following equivalency holds
\begin{equation}\label{sec_convergence_eq4_2}
{\bf y}=
\left[\begin{array}{l}
{\bf y}_{(1)} \\
{\bf y}_{(2)} \\
{\bf y}_{(3)}
\end{array}\right] =
\left[\begin{array}{l}
{\bf H}_{(1)}{\bf f} \\
{\bf H}_{(2)}{\bf f} \\
{\bf H}_{(3)}{\bf f}
\end{array}\right]={\bf G}{\bf f}
\in\mathbb{R}^{5mnN}
\end{equation}
This particular way of mapping expands the generalized ADMM into $J=3$ decoupled minimization steps that can be embedded in parallel and handled separately by three multicore processing computer. Also, $f_1=0$ provides a unique solution via a quadratic minimization step in the ADMM algorithm in (\ref{eq15}) which needs ${\bf G}$ to be full column rank. Please refer to \cite{FigueiredoBioucasDias:2010, AfonsoBioucas-DiasFigueiredo2011, AlmeidaFigueiredo:2013} for further assessments. By defining the above linear mapping, it can be claimed that the convergence conditions of the proposed BC-ADMM holds by means of Theorem $8$ in \cite{EcksteinBertsekas:1992}, where the updated sequences ${\bf y}^{k+1}_{(j)}$ to within some error $\varepsilon^k_{(j)}$ are absolutely summable and mainly requires ${\bf G}=[{\bf H}^T_{(1)}, {\bf H}^T_{(2)}, {\bf H}^T_{(3)}]^T$ to be full column rank by Proposition \ref{Proposition:Rank}. 
\begin{proposition}[Operator rank ${\bf G}$]\label{Proposition:Rank}
Consider the gradient operator ${\bf H}^T_{(1)}=\nabla$ in (\ref{eq3}). The rank of the operator is deficient $r(\nabla)=nmN-p^3$, where $p=\{0,2,1,1\}$ corresponds to Z-BC, P-BC, R-BC, and AR-BC, respectively. Furthermore, the stacked operator ${\bf G}=[{\bf H}^T_{(1)}, {\bf H}^T_{(2)}, {\bf H}^T_{(3)}]^T$ is of full-rank $r({\bf G})=nmN$.
\end{proposition}\vspace{-0in}
The proof of Propositions \ref{Proposition:Rank} is presented in Appendix \ref{Sec:Appendix2}. Furthermore, an optimal convergence example of BC-ADMM is demonstrated in Figure \ref{Fig:convergence} for the evolution of convergent sequence  ${\bf f}_k$ with minimum tolerance error of $\epsilon=10^{-7}$, where AR-BC-ADMM provides better convergence rate compared to P-BC-ADMM.}

\vspace{-.1in}

\section{Experimental Evaluations}\label{Sec:Experiments}
\begin{figure*}[htp]
\centerline{
\subfigure[]{\includegraphics[height=0.125\textwidth]{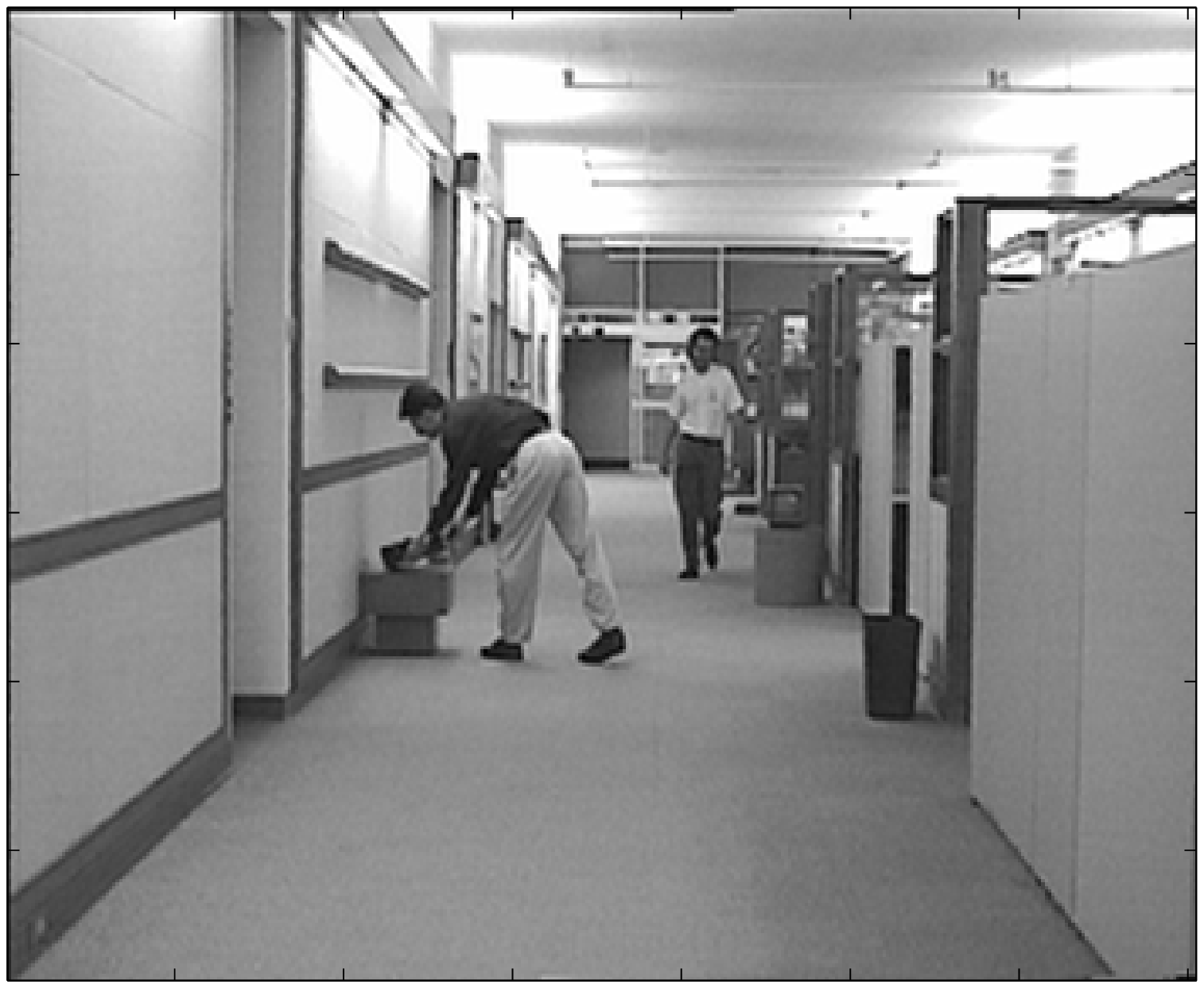}}
\subfigure[]{\includegraphics[height=0.125\textwidth]{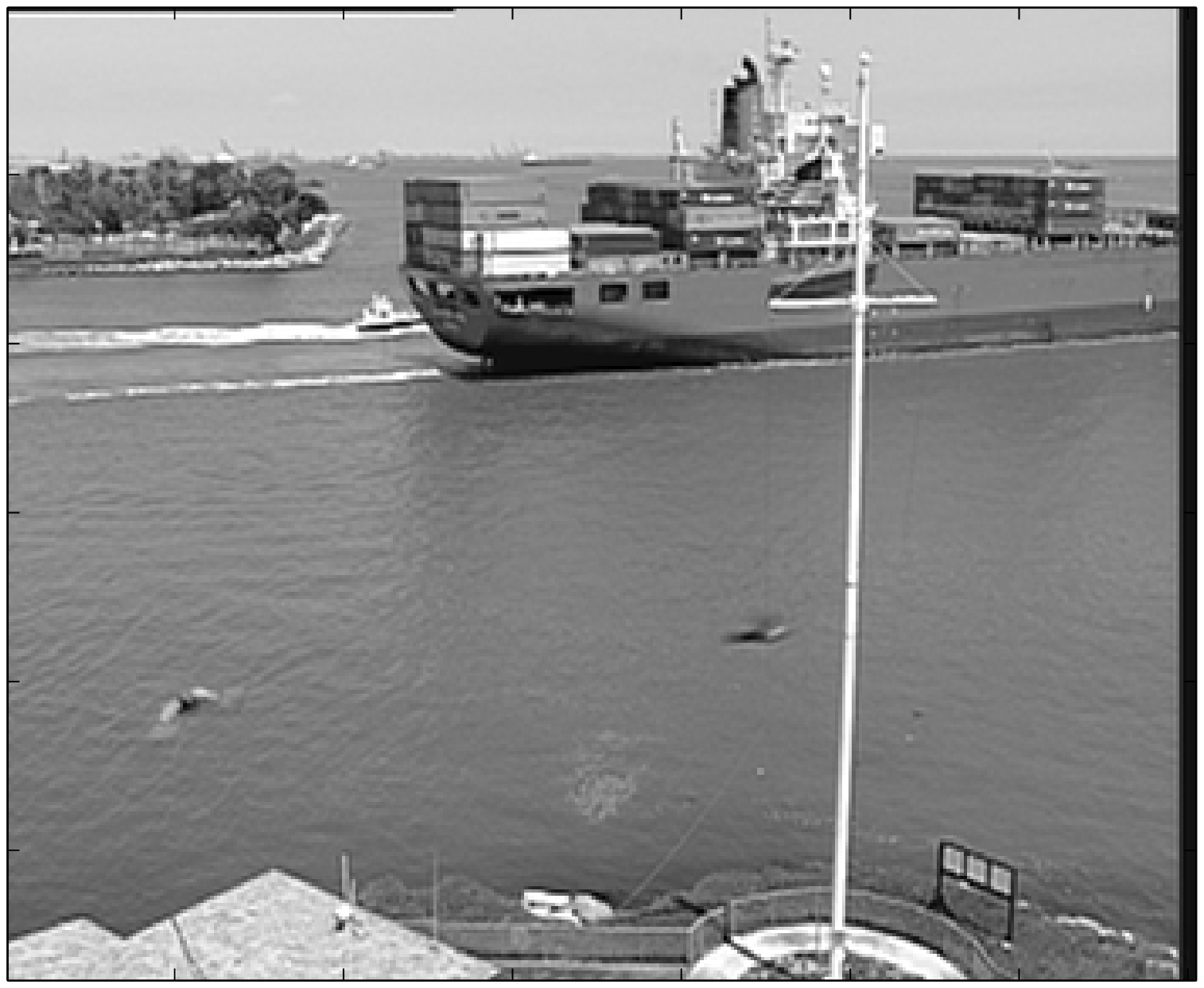}}
\subfigure[]{\includegraphics[height=0.125\textwidth]{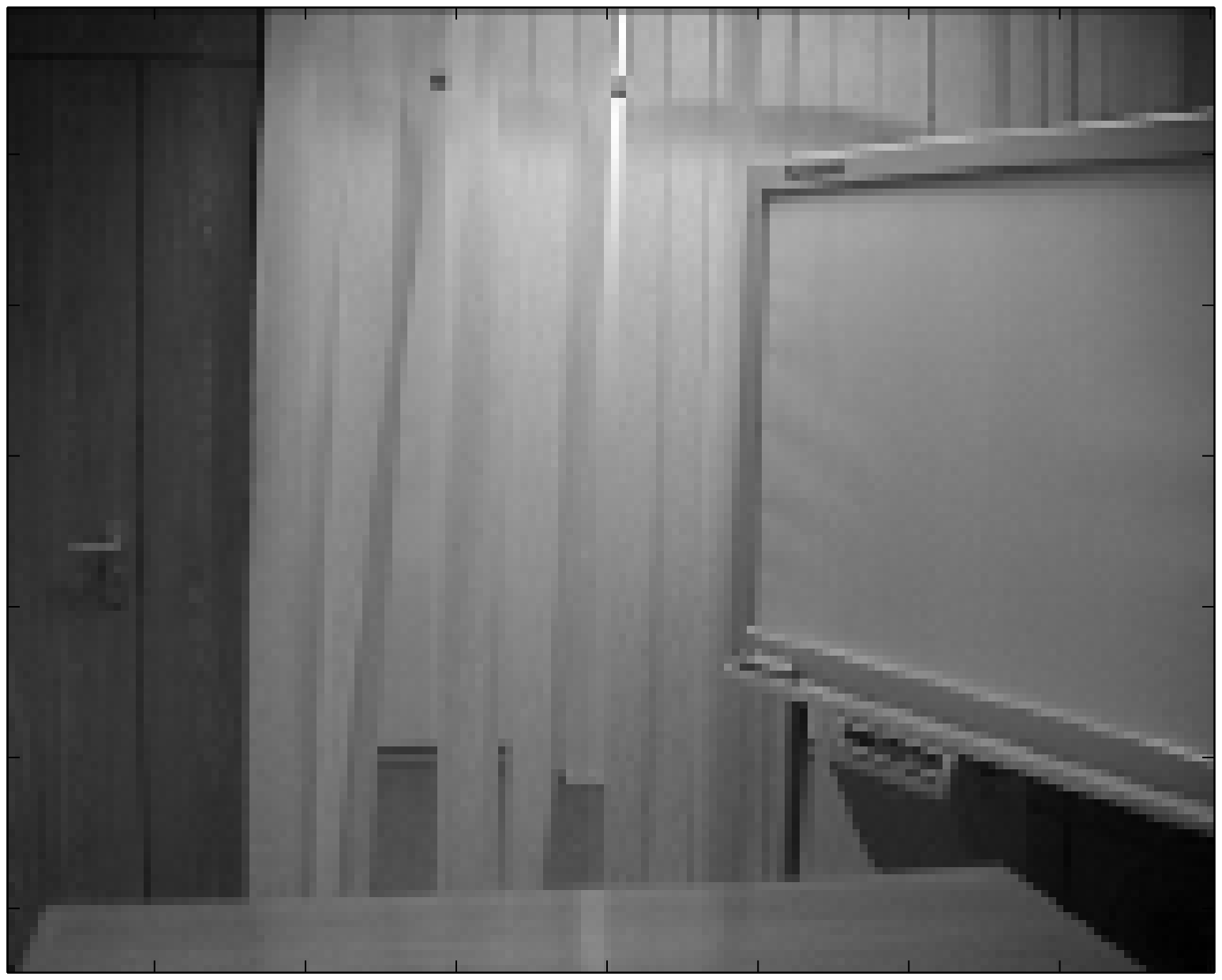}}
\subfigure[]{\includegraphics[height=0.125\textwidth]{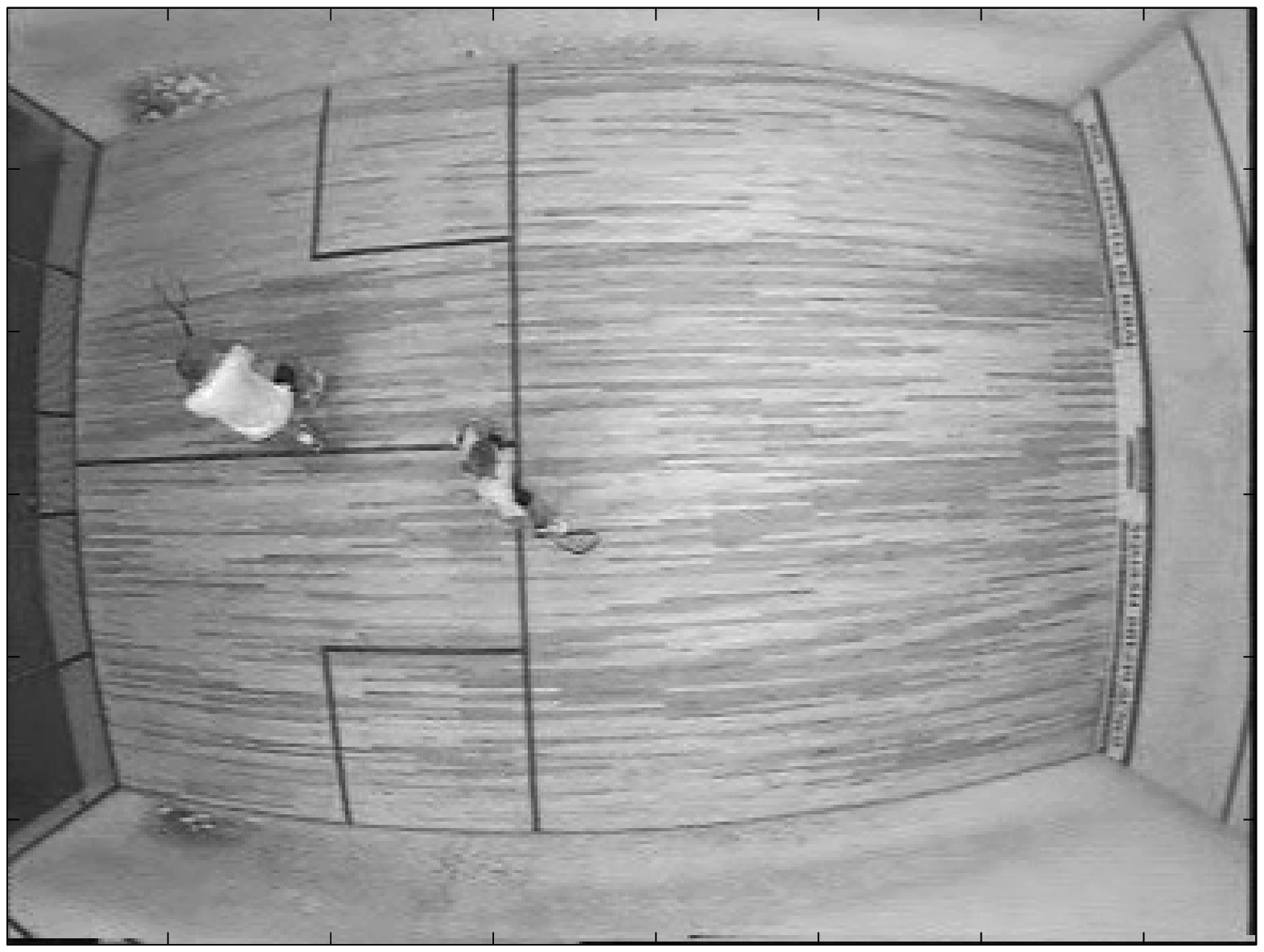}}
}\vspace{-.1in}
\caption{Test sequence frame No. $17$ out of $32$ from: (a) \textit{Hall-Monitor}; (b) \textit{Container}; (c) \textit{Office Environment}; and (d) \textit{Squash2}}
\label{Fig:TestSequences}
\end{figure*}

We evaluate BC-ADMM's performance to recover video frames from under-sampled measurements within the CVS paradigm. The proposed method was tested on four video sequences recorded in stationary background (no camera movements): \textit{Hall-Monitor} and \textit{Container} with CIF spatial resolution of $288\times 352$ pixels registered in 30 frames/s (fps), \textit{Office Environment} with $128\times 160$ pixels registered in 25 fps \cite{LiHuangTian:2004}, and \textit{Squash2} with $288\times 384$ pixels registered in 25 fps \cite{CVBASE06}. Based on the aforementioned frame rates, we fix the temporal resolution of the test data to $32$ frames to observe significant motion trajectories for analysis as follows. We select frames $101-132$ from \textit{Hall-Monitor} containing complex background with two men walking in the corridors and frames $240-271$ from \textit{Container} representing a container-ship moving slowly with sea/flag waving and two flying birds. Frames $1-32$ are segmented from \textit{Office Environment} with waving curtains and finally frames $1-32$ of  \textit{Squash2} are selected showing two men playing squash with fast motion trajectories. We select these sequences due to different levels of spatio-temporal complexities that is quantified by spatial information (SI) and temporal information (TI) indexes \cite{ITUTP910}. \textit{Hall-Monitor} and \textit{Container} contain high SI and mid TI indexes \cite{SimoneNaccariEbrahimi:2009}, while \textit{Office Environment} contain low SI and low TI indexes, and \textit{Squash2} contains high SI and high TI indexes. A test frame of four test clips is shown in Figure \ref{Fig:TestSequences}.
\vspace{-.1in}

{
\subsection{Parameter Selection for BC-ADMM}\label{sec_param_select}
The gradient operator ${\bf\nabla}$ in Algorithm \ref{algo_BC_ADMM} is constructed by the scheme design in Section \ref{sec_TV_encoding}. We select AR-BC and P-BC to construct the matrices ${\bf D}$ for every tensorial dimensions designed in Section \ref{sec_numerical_design}. The convoluting FIR filter ${\bf d}$ is defined by the HO accuracy NR method in Appendix \ref{Sec:NRD} to design the derivative matrices. The size of the filter length $L$ and order of exact polynomial $p$ identifies the maximum accuracy achievable for differentiation. Setting $L$ from low to high will increase the zooming range in frequency domain for differentiation. Provided with the second parameter $p$ defines the polynomial degree that the filter is exact on low frequencies. Setting $p$ low suppress high frequency band. In contrast, if we are keen to observe high frequency component in the signal, then $p$ should be close to $L-1$. This accuracy parameter can vary between $p\in\{2,4,\hdots,L-1\}$ in filter design for different purposes. Please refer to Appendix \ref{Sec:NRD} for the design procedure. The choice of $L$ is related to the closeness of the maximum frequency spectrum of the signal to the Nyquist rate. If the maximum frequency is close to the Nyquist rate then high accuracy parameter should be set for differentiation, i.e. $p=L-2$. In this study, for differentiating the class of videos such as Hall-Monitor, setting $L=27$ and $p=25$ satisfies the amount of zooming window in frequency range and able to observe the most frequency components in the signal either in spatial or temporal domain. For more information on the latter please refer to Section \ref{HO_Accuracy_diff}. Also, the number of collected samples for differentiation should be beyond the filter length, i.e. $m,n,N\geq L$. Otherwise, the effect of boundaries will induce to the whole range of the signal domain during recovery. Particularly, this is of importance in temporal domain when limited number of frames is provided. We considered $N=32$ frames for our analysis. The consequence of choosing $N$ close to $L$ is boundary artefacts affects frames $\{1:(L-1)/2\}$ and $\{N-(L-1)/2+1:N\}$. This is clearly observable in Section \ref{sec_frame_index_analysis} in Table \ref{Table:FrameIndex}, where the accuracy of reconstruction is deviated for frame indexes $\{1,2,\hdots,13\}$ and $\{20,21,\hdots,32\}$.} The rest of the parameters selected for BC-ADMM is as follows. {We acquire the compressed samples uniformly at random from every individual frame using two linear transforms discussed in Section \ref{Sec:CVSsampling} and \ref{SubSec:f}: (a) ${\bf\Phi}$: matrix formed by populating the entries with i.i.d. Gaussian variables followed by orthonormalization of the rows and (b) Walsh-Hadamard matrix.} The sparse representation is favoured by Wavelet dictionaries using symmlets of order $10$ with $4$ level of decomposition deployed by \textsc{UviWave 300}\footnote{available online at \url{http://cas.ensmp.fr/~chaplais/UviWave/About_UviWave.html}} wavelet toolbox. Peak-Signal-to-Noise-Ratio (PSNR) was used for metric evaluations between the reconstructed and the original videos. The regularization parameters were manually set to $c_{1}=80$, $c_{2}=10$, $c_{3}=1000$ to yield the best PSNR for all experimental setup. Moreover, the  parameters were adjusted to $\mu_{(1)}=4$, $\mu_{(2)}=4$, $\mu_{(3)}=40$ and $\rho=(\sqrt{5}+1)/2$ for optimal convergence with $\epsilon=10^{-4}$ tolerance error. The summary of parameter selection of BC-ADMM is presented in Table \ref{tab:BC-ADMM-Params}.
\vspace{-.2in}
\begin{table}[ht]
\caption{BC-ADMM parameter design}\vspace{-.1in}
\label{tab:BC-ADMM-Params}
\centering
\begin{tabular}{ll}    \toprule
Description & Parameter values \\\midrule
Kernel ${\bf d}$ used to construct ${\bf\nabla}$ &  $L=27$, $p=25$ \\ 
Type of  boundary& P-BC, AR-BC \\
Sampling transform operator & ${\bf\Phi}$:  orthonorm. i.i.d Gauss. matrix/  \\
& Walsh-Hadamard matrix \\
Sparse transform operator & ${\bf\Psi}$: Symmlets, ord. 10, dec. level 4 \\
Regularizing parameters & $c_1=80$, $c_2=10$, $c_3=1000$ \\
Optimality convergence parameters & $\mu_{(1)}=4$, $\mu_{(2)}=4$, $\mu_{(3)}=40$\\
& $\rho=(\sqrt{5}+1)/2$, $\epsilon=10^{-4}$ \\\bottomrule
\hline
\end{tabular}
\end{table}\vspace{-.25in}
\subsection{Competing Algorithms}\label{sec_competing_algorithms}
The performance of the BC-ADMM is evaluated and compared with the following current state-of-the-art algorithms:
\begin{itemize}
\item \text{TFOCS}\footnote{available online at \url{http://cvxr.com/tfocs/}}\cite{BeckerCandesGrant:2011}: is a flexible conic formulation solver that provides a solution to $\text{TV}\ell_1/\ell_1$ in (\ref{eq4_2}) using a proximal gradient method for the minimization in several steps: first it finds the conic formulation of the objective followed by smoothening the dual function using Nesterov's approach and finally solves the smooth function by first order method (gradient projection technique) \cite{BeckerCandesGrant:2011}. The gradient projection method adapts a backtracking line-search algorithm for proper step-size selection. The following are the tuned parameter that are selected by trial and error to achieve the best performance results. Scaling parameter $\mu=5$, data fidelity error $\delta=10$, and regularizing parameters for TV and $\ell_1$ norms are set to $\alpha=10^4$ and $\beta=10^2$, respectively. Maximum of $100$ iterations are allowed per each continuation with the total of $4$ steps. Tolerance decreasing rate at each continuation is fixed to $\textsc{betaTol} = 2$.

\item \text{NESTA}\footnote{available online at \url{http://www-stat.stanford.edu/~candes/nesta/}}\cite{BeckerBobinCandes:2011}: is a robust and efficient large-scale optimization method for recovering sparse signals using Nesterove's smoothing techniques \cite{BeckerBobinCandes:2011}. The smoothing parameter is dynamically updated using continuation methods in \cite{HaleYinZhang:2008} to accelerate the convergence of the algorithm. We chose TV model for basis pursuit denoising problem, i.e. $\|{\bf F}\|_{\text{TV}}+\gamma/2\|\overline{\bf\Phi}_{\bf\Omega}{\text{vec}\left({\bf F}\right)}-{\bf b}_{\bf\Omega}\|^2_2$. We set $10$ number of continuation steps and tune the smoothing parameter to $\mu=0.2$ to achieve the best results.

\item \text{TVAL3}\footnote{available online at \url{http://www.caam.rice.edu/~optimization/L1/TVAL3/}}\cite{TVAL3}: is an ADMM based approach for TV minimization models to reconstruct images from their measurements. The solver uses steepest descent method to update ADMM sub-problems favoured by Barzilai and Borwein (BB) method \cite{TVAL3} for step-length selection. We have selected TV/L2 model for basis pursuit denoising problem, i.e. $\|{\bf F}\|_{\text{TV}}+\gamma/2\|\overline{\bf\Phi}_{\bf\Omega}{\text{vec}\left({\bf F}\right)}-{\bf b}_{\bf\Omega}\|^2_2$. The primary and secondary penalty parameters are set to $\gamma=2^{10}$ and $\beta=2^7$, respectively, decided by trial and error to achieve the best results. Also, anisotropic model for TV regularization is selected, i.e. $\xi=1$.
\end{itemize}

We have modified the MATLAB codes in TFOCS, NESTA, and TVAL3 from 2D to 3D image problem. The gradients for TV regularizer are fashioned with $[-1,1]$ filter scheme and initialized with P-BC. {To avoid signal discontinuities at the end of video frames, we have duplicated the data in reflective form to preserve temporal continuity and hence to be more suitable for P-BC.} Similar tolerance error is also considered, i.e. $10^{-4}$ for stoping criteria. {TVAL3 and NESTA are the reduced version of $\text{TV}\ell_1/\ell_1$ problem in (\ref{eq4_2}) using only $\text{TV}\ell_1$ form, while TFOCS solves the same problem. The rationale of choosing the reduced versions, is to study the amount of accuracy compensated by sparse regularization, i.e. $\ell_1$-norm in TFOCS compared to TVAL3 and NESTA.}\vspace{-.1in}

\subsection{Sampling Rate Analysis}\label{sec_sampling_rate_analysis}
\begin{figure*}[ht]
\centerline{
\subfigure[]{\includegraphics[width=0.19\textwidth]{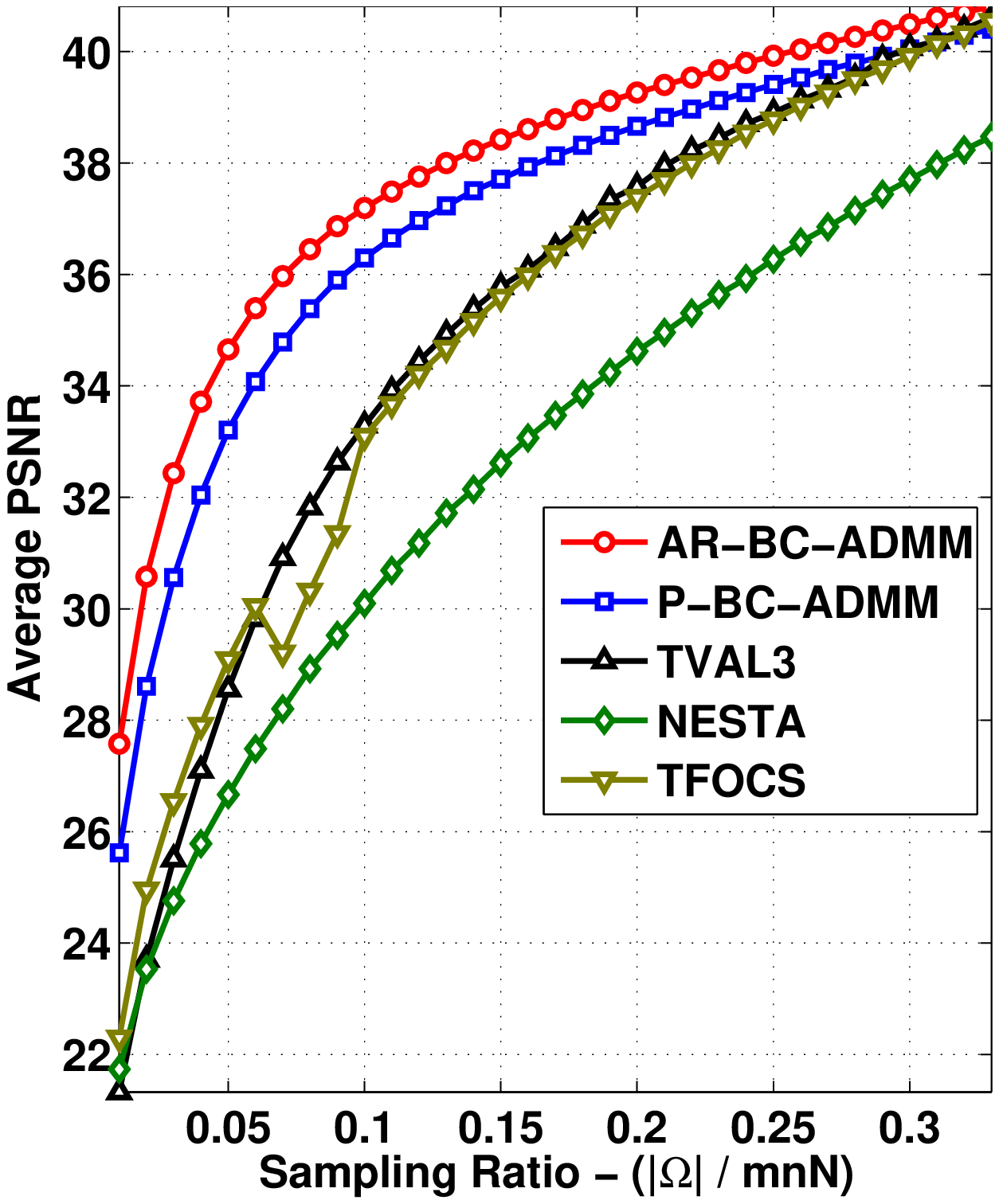}
\label{Fig:VaryingSamplingHall}}
\subfigure[]{\includegraphics[width=0.19\textwidth]{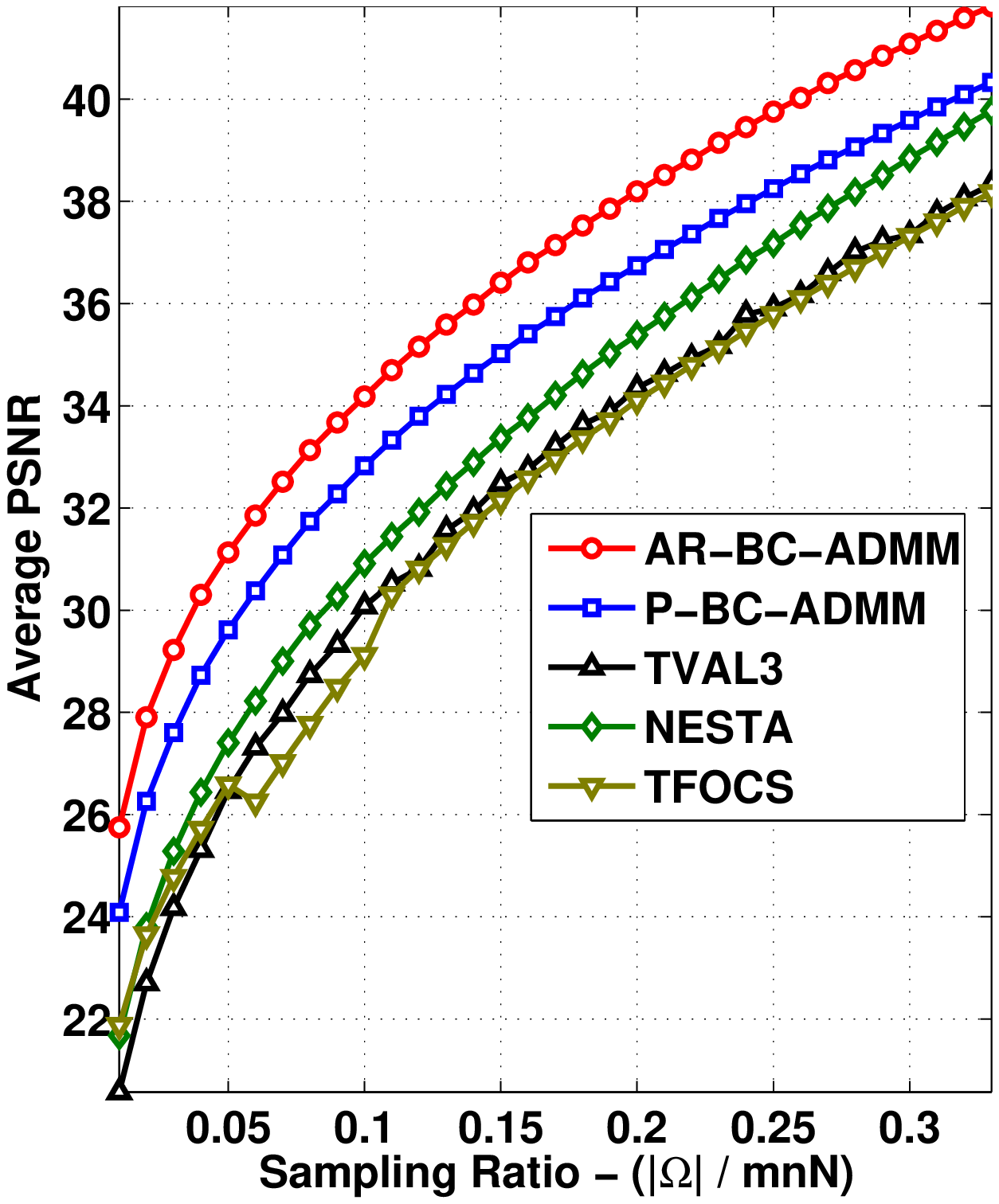}
\label{Fig:VaryingSamplingContainer}}
\subfigure[]{\includegraphics[width=0.19\textwidth]{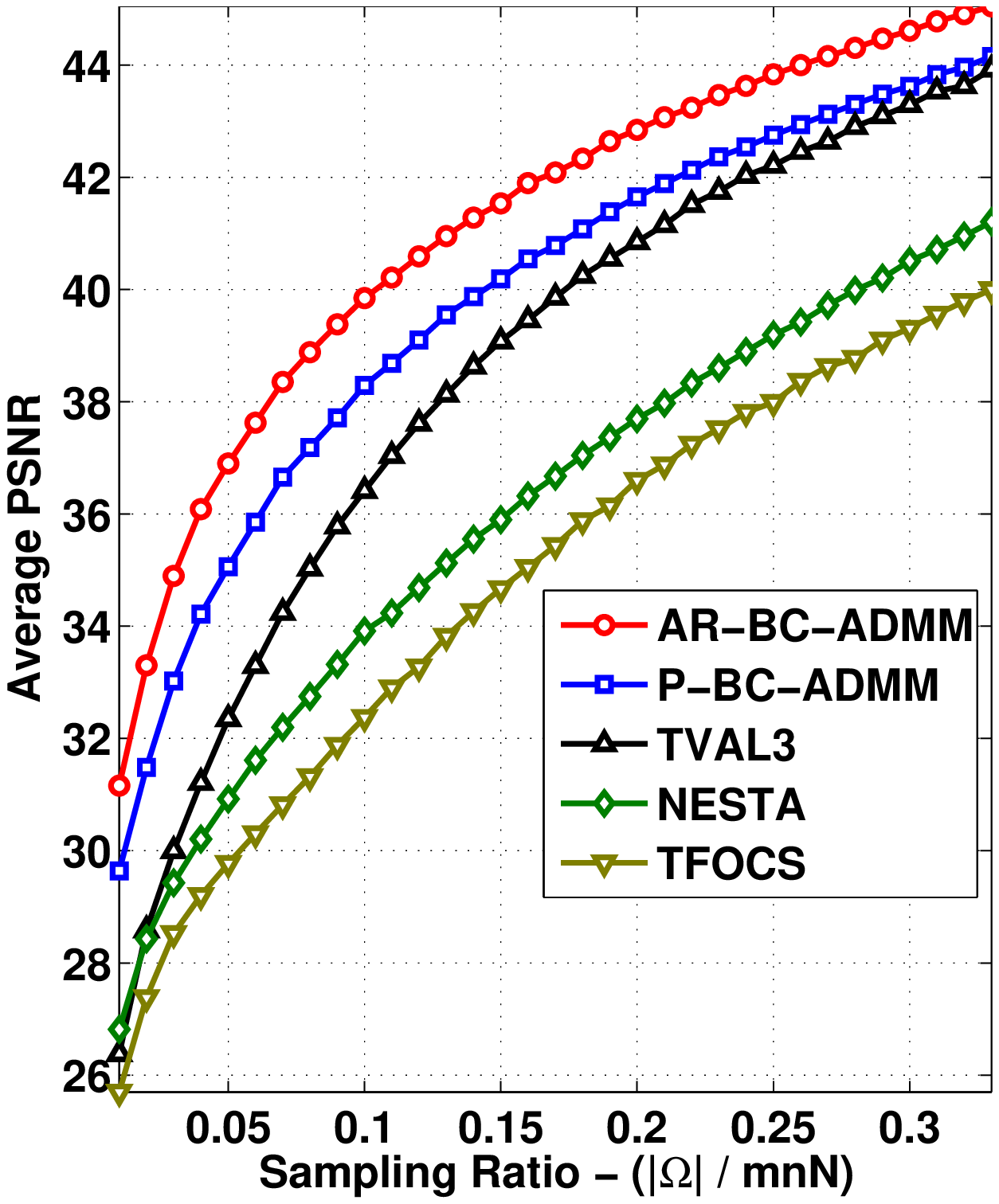}
\label{Fig:VaryingSamplingCurtain}}
\subfigure[]{\includegraphics[width=0.19\textwidth]{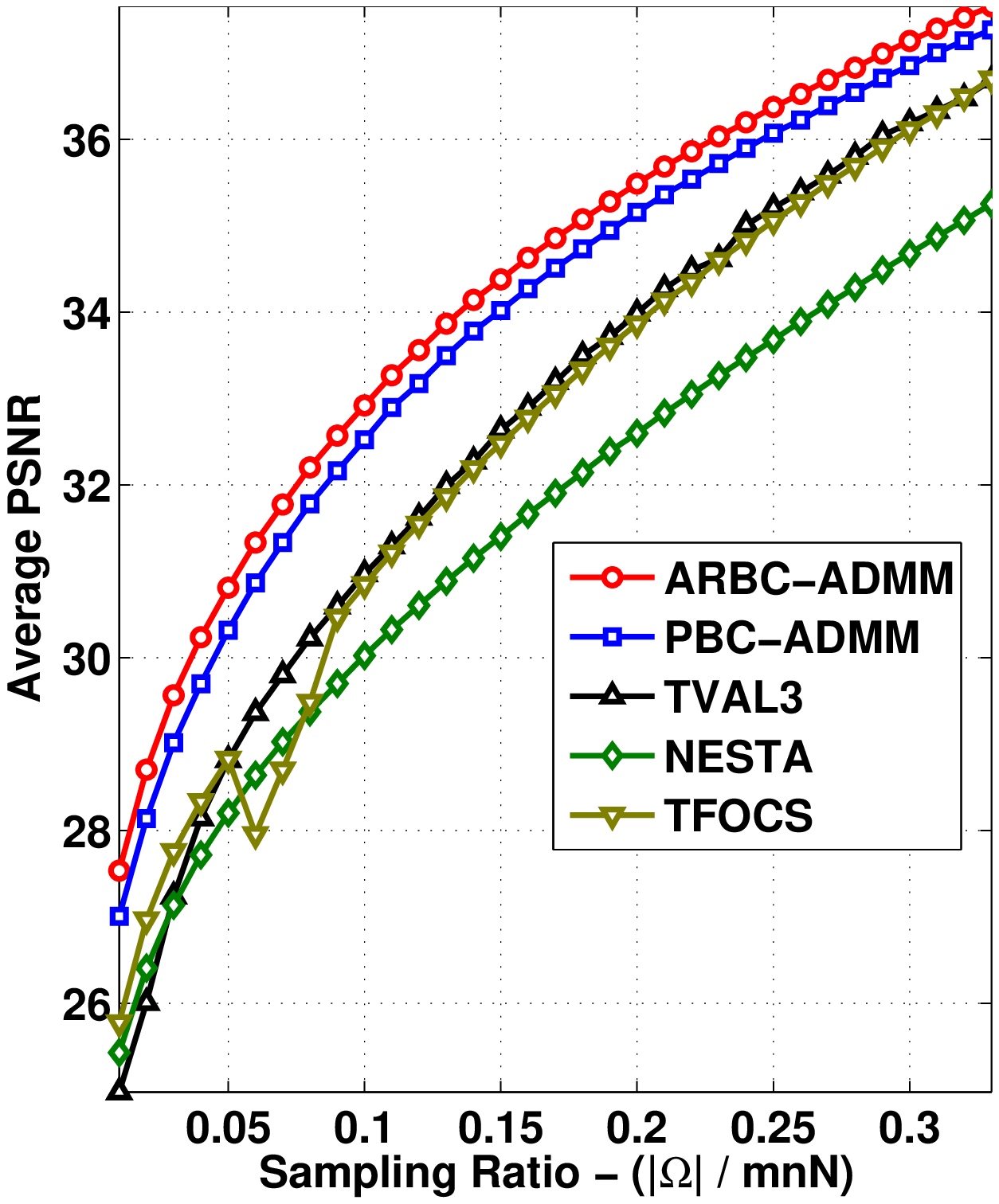}
\label{Fig:VaryingSamplingSquash}}
\subfigure[]{\includegraphics[width=0.19\textwidth]{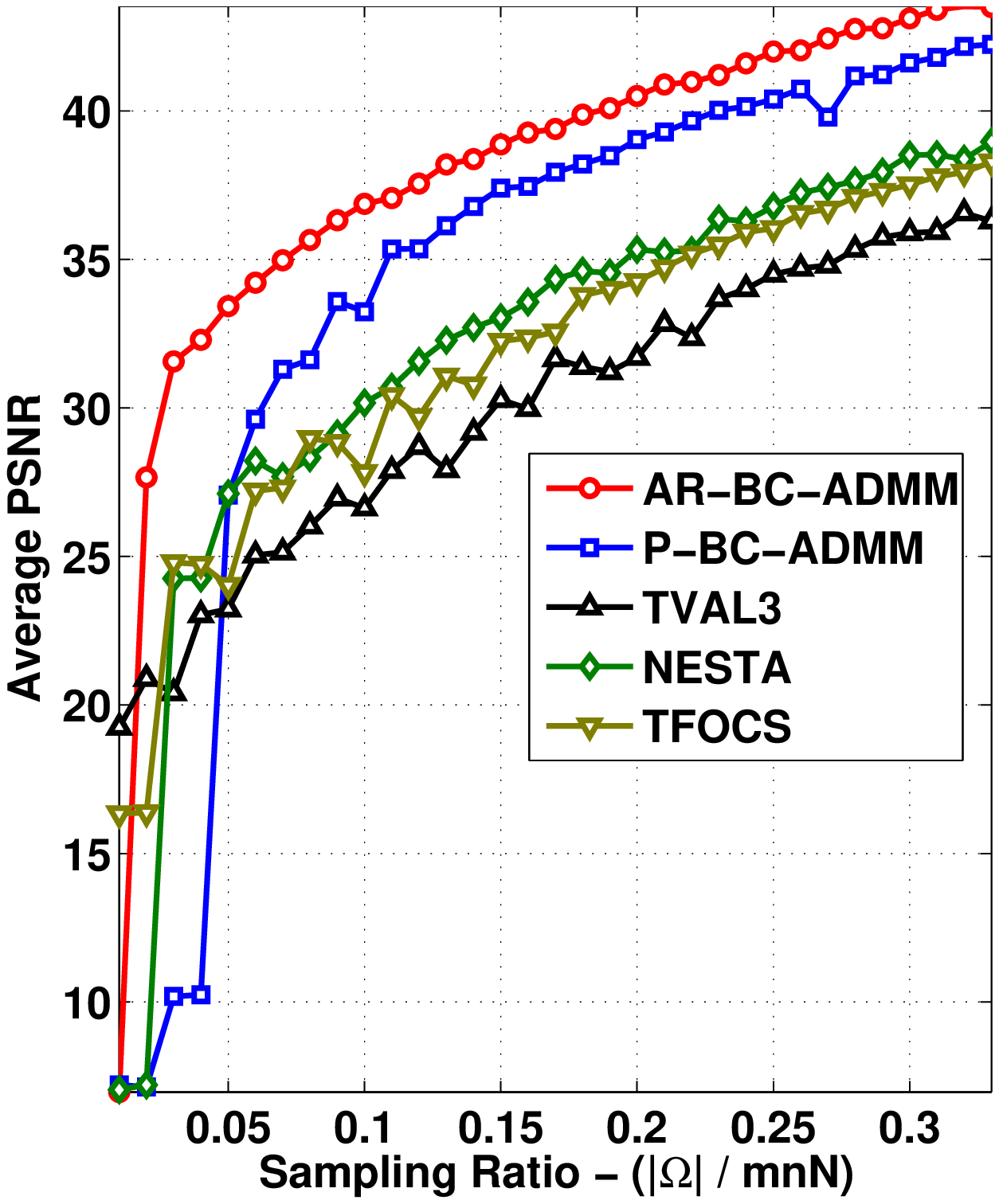}
\label{Fig:VaryingSamplingCurtain_WalshHadamard}}
}\vspace{-.125in}
\caption{Empirical performance of BC-ADMM, TFOCS, NESTA and TVAL3 over varying sampling rate and transforms tested on: (a) ``Hall-Monitor'', ${\bf\Phi}$-Gaussian; (b) ``Container'', ${\bf\Phi}$-Gaussian; (c) ``Office environment'', ${\bf\Phi}$-Gaussian; and (d) ``Squash2'', ${\bf\Phi}$-Gaussian; (e) ``Office environment'', ${\bf\Phi}$-Walsh-Hadamard}
\label{Fig:SamplingRateAnalysis}
\end{figure*}

The first experiment is conducted by varying different number of random samples $|{\bf\Omega}|$ acquired from all $32$ frames with total possible number of $mnN$. The sample ratio $|{\bf\Omega}|/mnN$ is varied from $1\%$ to $33\%$ with $1\%$ incremental steps. {The Monte Carlo simulation is carried out for $10$ experimental trials on different sampling rate generated at random combinations with uniform distribution and finally averaged the PSNR accuracy over all frames in the video.} Figure \ref{Fig:SamplingRateAnalysis} demonstrates the empirical performance algorithms on four different video clips using two sampling transforms ${\bf\Phi}$: random i.i.d. Gaussian, and Walsh-Hadamard transform. The reconstruction accuracy is shown as a function of sampling ratio. The rank observation of the proposed BC-ADMM outperforms over all sampling rates compared to TFOCS, NESTA and TVAL3. Despite stable recovery via random i.i.d. Gaussian sampling, Walsh-Hadamard transform provides unstable recovery in very low sampling rates, i.e. $|{\bf\Omega}|/mnN\sim 0.01$. This is expected since Walsh bases, compared to random Gaussian, does not have uniform RIP condition \cite{HowardCalderbankSearle:2008}. Notice the significant improvement achieved by BC-ADMM in low sampling rates. This is due to HO accuracy derivative kernel implemented to connect consecutive frames for TV regularization. The HO accuracy kernel establishes  high order polynomials of filter length ($L=27$ in this experiment) which provides high feasibility range in low sampling rates for derivative approximation compared to conventional $[-1,1]$ scheme. For instance, in Figure \ref{Fig:VaryingSamplingHall} using only $5\%$ CVS samples, the accuracy of ARBC-ADMM is $34.65$dB which is remarkably $5.57$dB higher than the best competitor TFOCS with $29.08$dB average PSNR.

Further evaluation is provided in Table \ref{Table:FrameIndex17RecoveryHallMonitor} by visualizing the reconstructed mid-frame (No. $17$) of \textit{Hall-Monitor} test clip by means different algorithms. Notice the robustness of the visual quality of BC-ADMM through degrading the sampling ratio. Furthermore, Table \ref{Table:FrameIndex17RecoveryHallMonitorResidue} illustrates the absolute residual error of the reconstructed frames in Table \ref{Table:FrameIndex17RecoveryHallMonitor} compared to the reference frame in Figure \ref{Fig:TestSequences}. {The variation of residual error is demonstrated in gray colour from white to back as it increases.} The Normalized-Mean-Square-Error (NMSE) is calculated to evaluate the performance. Improvement via BC-ADMM is highly evident in recovering edge information in the images especially in low sample rate. {These information are related to both motion and non-motion parts. To emphasize this argument, four different patches: P$1$, P$2$, P$3$, and P$4$ are considered in the recovery frames, where P$1$ and P$4$ contain motion residuals and P$2$ and P$3$ are from background image containing textural information.}

\begin{table*}
\renewcommand{\arraystretch}{1.3}
\caption{Reconstructed frame No. $17$ via five different CVS approaches and three various sampling rates}
\label{Table:FrameIndex17RecoveryHallMonitor}\vspace{-.1in}
\begin{center}
\begin{tabular}{|c|c|c|c|c|c|}
\cline{2-6}
\multicolumn{1}{c|}{} & AR-BC-ADMM & P-BC-ADMM & TFOCS & NESTA & TVAL3 \\
\hline
{\hspace{-.05in}\begin{sideways} $|{\bf\Omega}|/mnN=0.01$ \end{sideways}\hspace{-.05in}} &
{\hspace{-.05in}\includegraphics[width=0.185\textwidth]{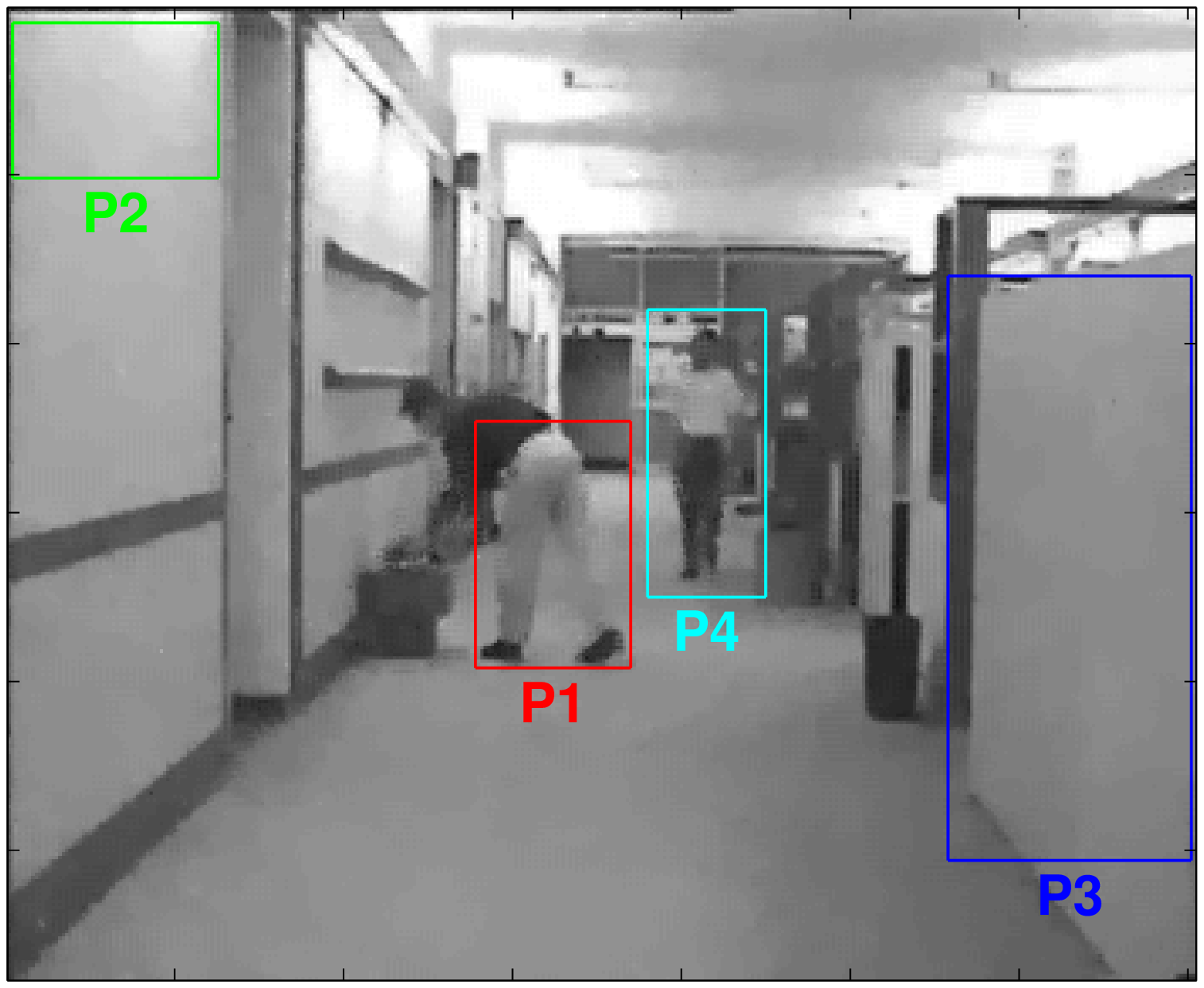}\hspace{-.05in}} &
{\hspace{-.05in}\includegraphics[width=0.185\textwidth]{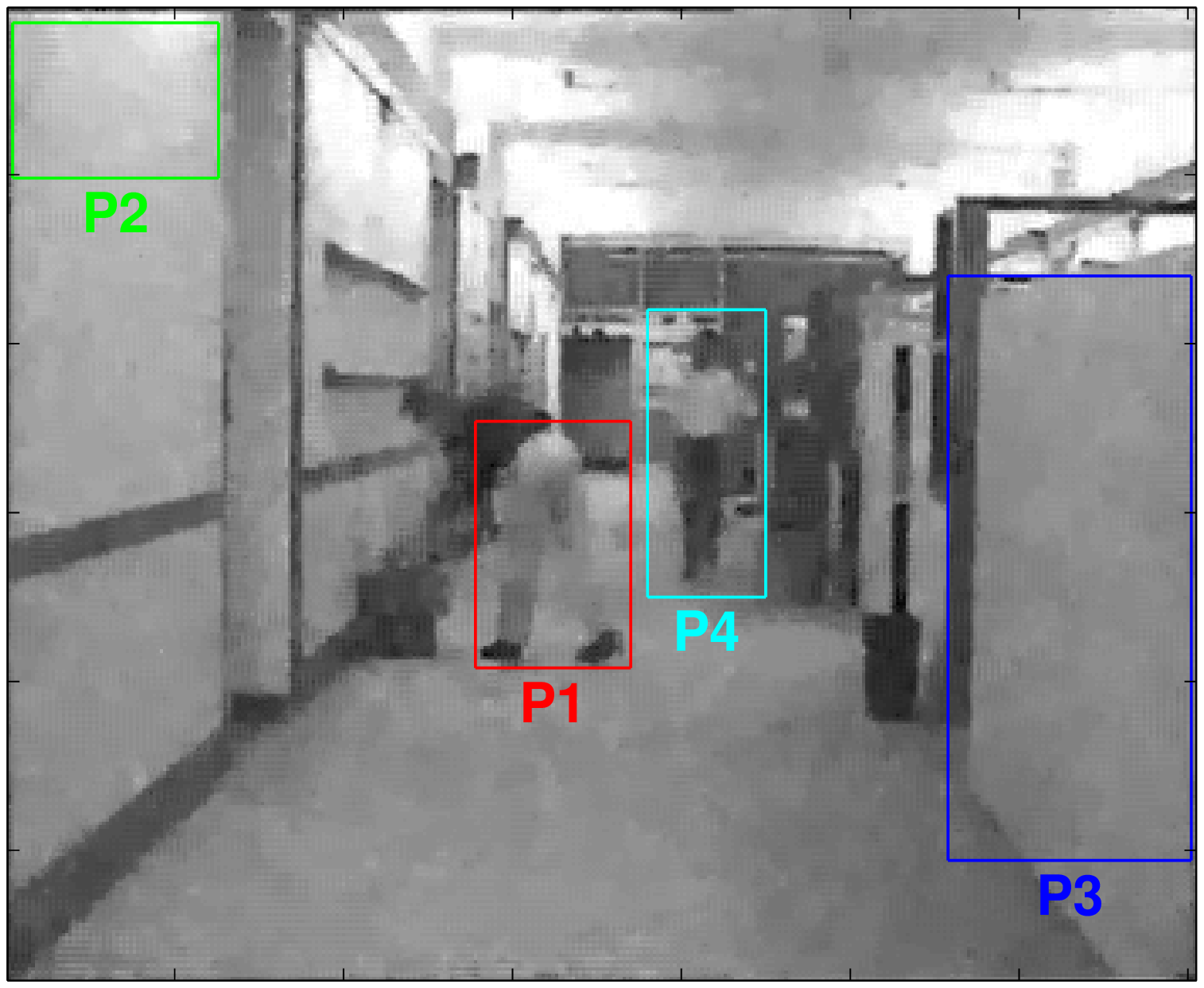}\hspace{-.05in}} &
{\hspace{-.05in}\includegraphics[width=0.185\textwidth]{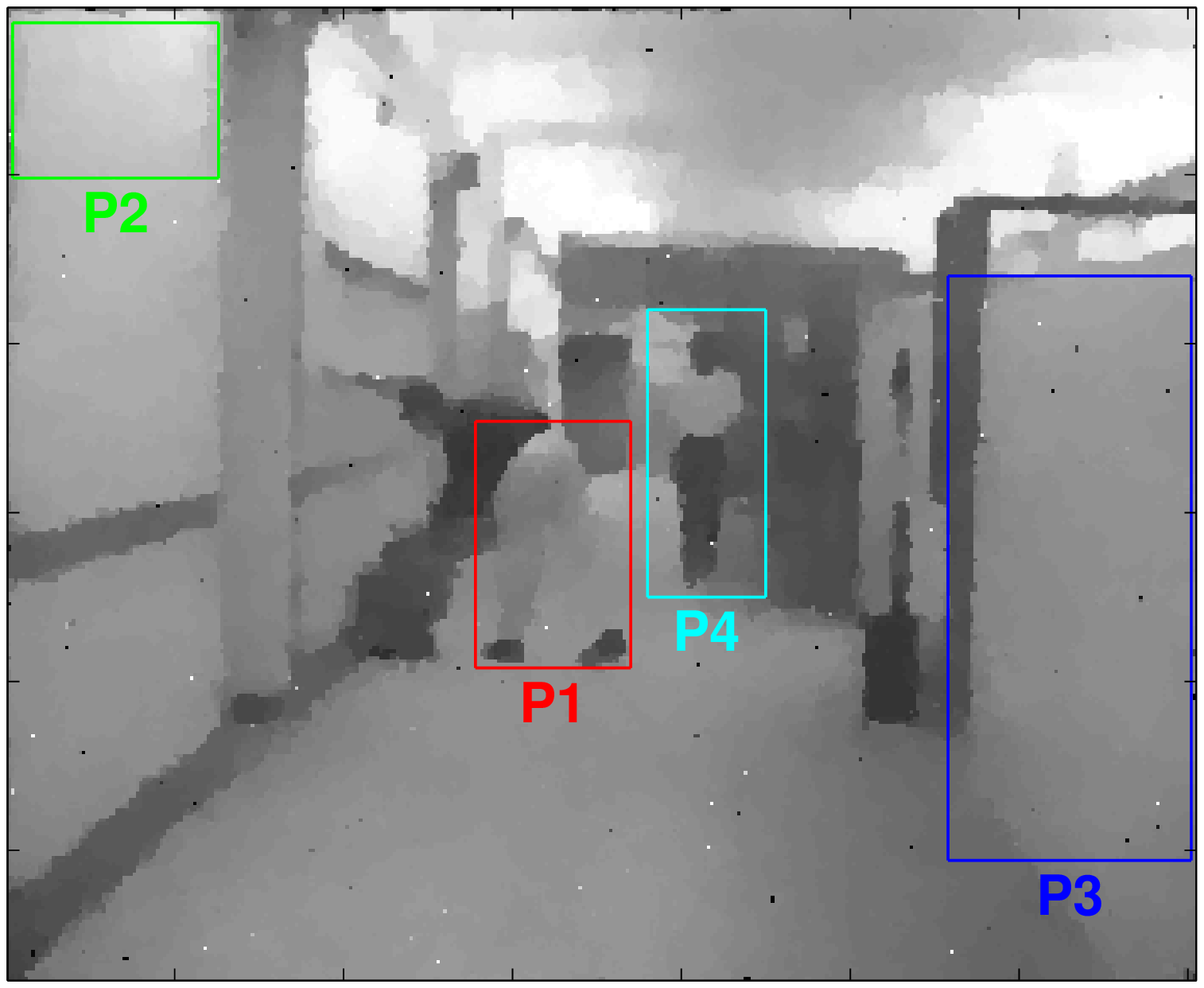}\hspace{-.05in}} &
{\hspace{-.05in}\includegraphics[width=0.185\textwidth]{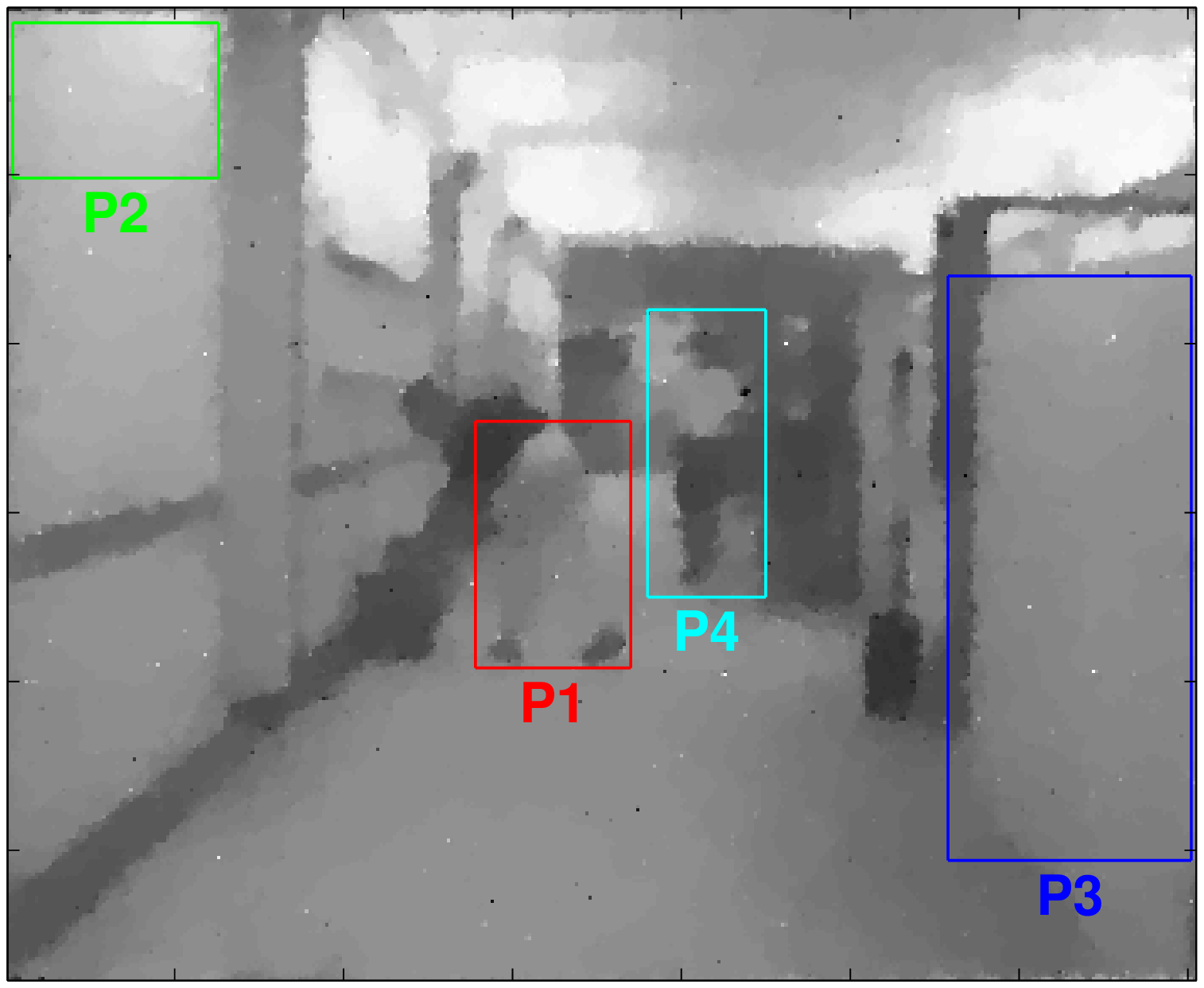}\hspace{-.05in}} &
{\hspace{-.05in}\includegraphics[width=0.185\textwidth]{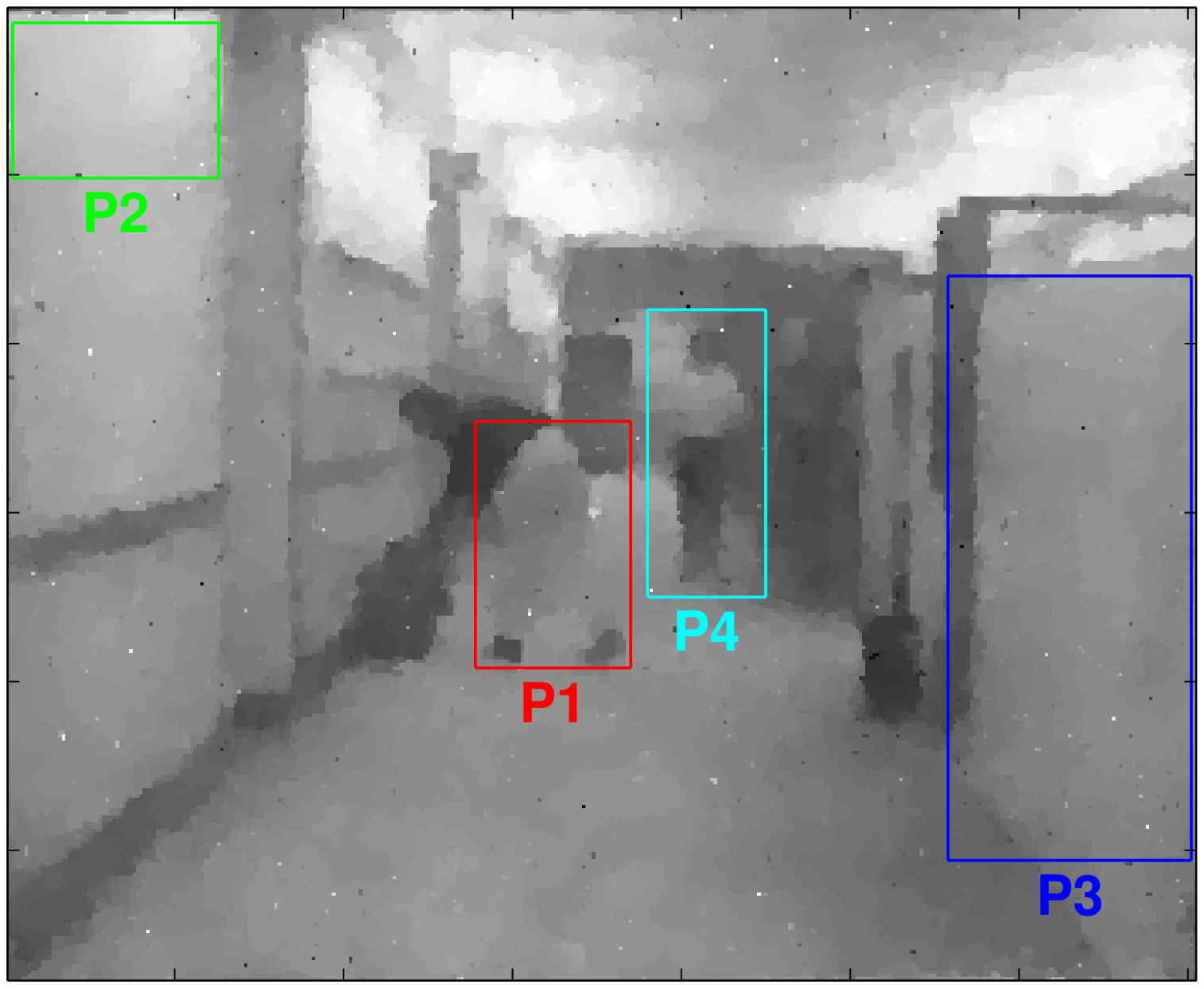}\hspace{-.05in}} \\
& PSNR = 28.02  & PSNR = 26.03 & PSNR = 21.46 & PSNR = 21.82 & PSNR = 22.36 \\
\hline
{\hspace{-.05in}\begin{sideways} $|{\bf\Omega}|/mnN=0.03$ \end{sideways}\hspace{-.05in}} &
{\hspace{-.05in}\includegraphics[width=0.185\textwidth]{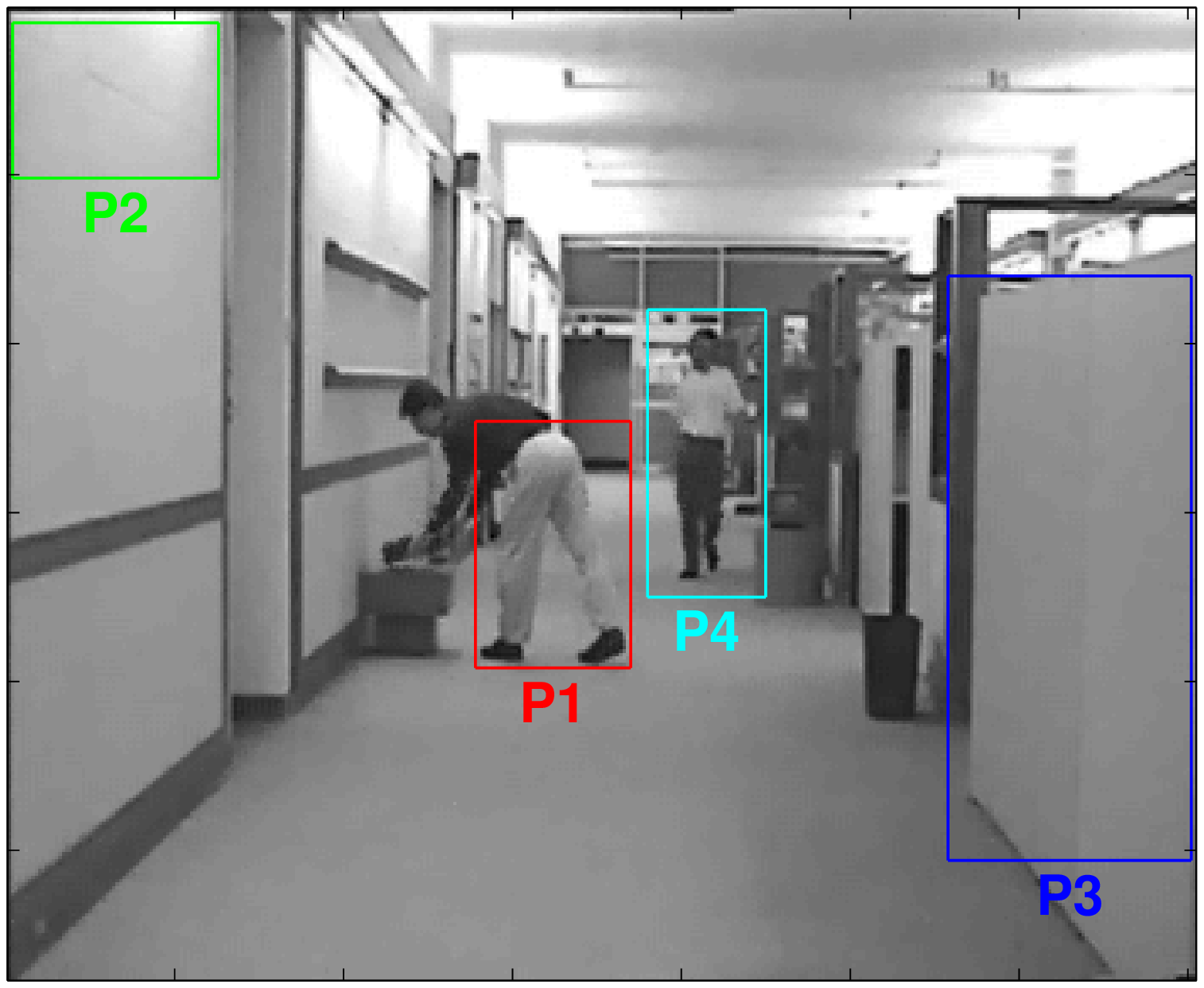}\hspace{-.05in}} &
{\hspace{-.05in}\includegraphics[width=0.185\textwidth]{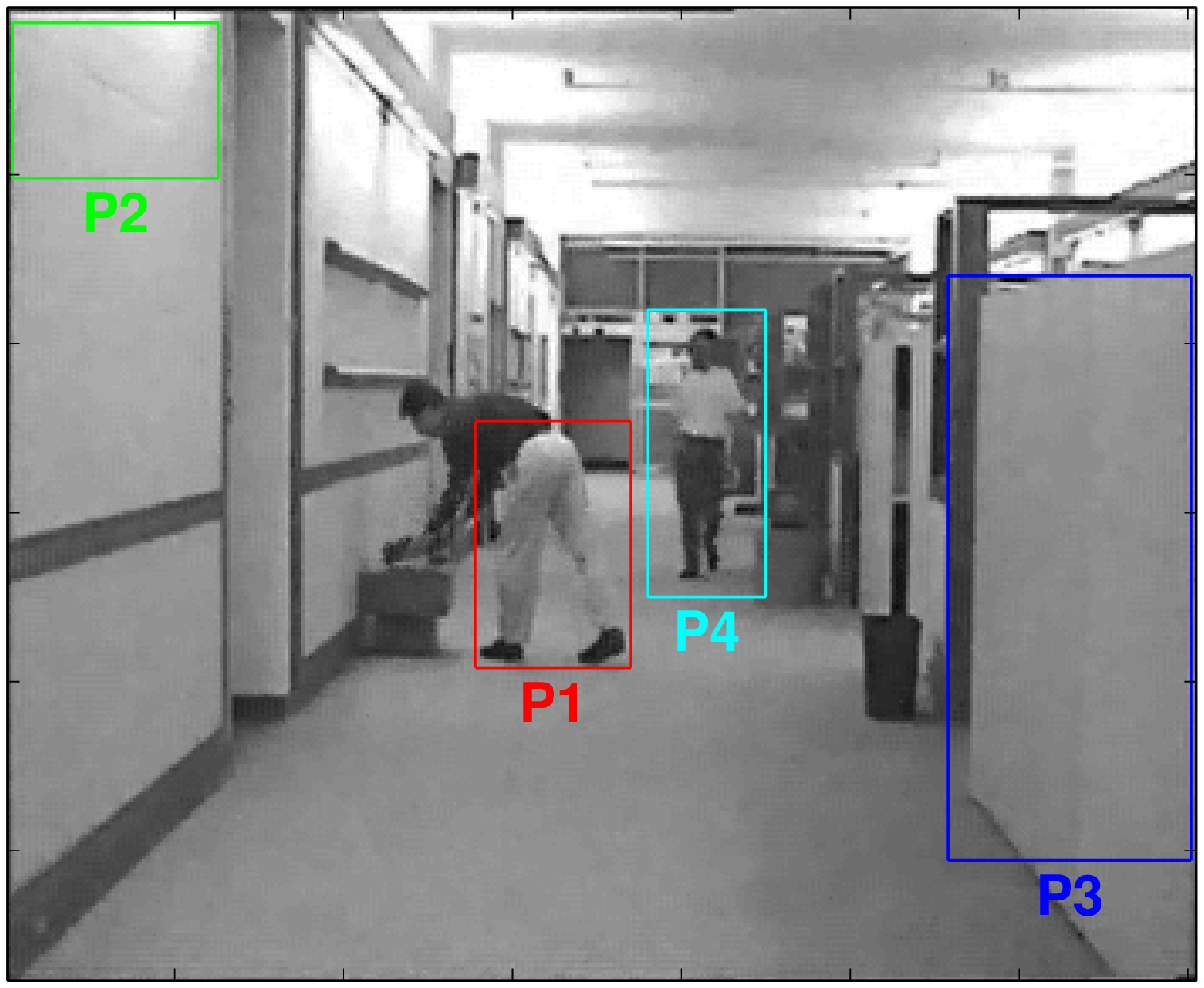}\hspace{-.05in}} &
{\hspace{-.05in}\includegraphics[width=0.185\textwidth]{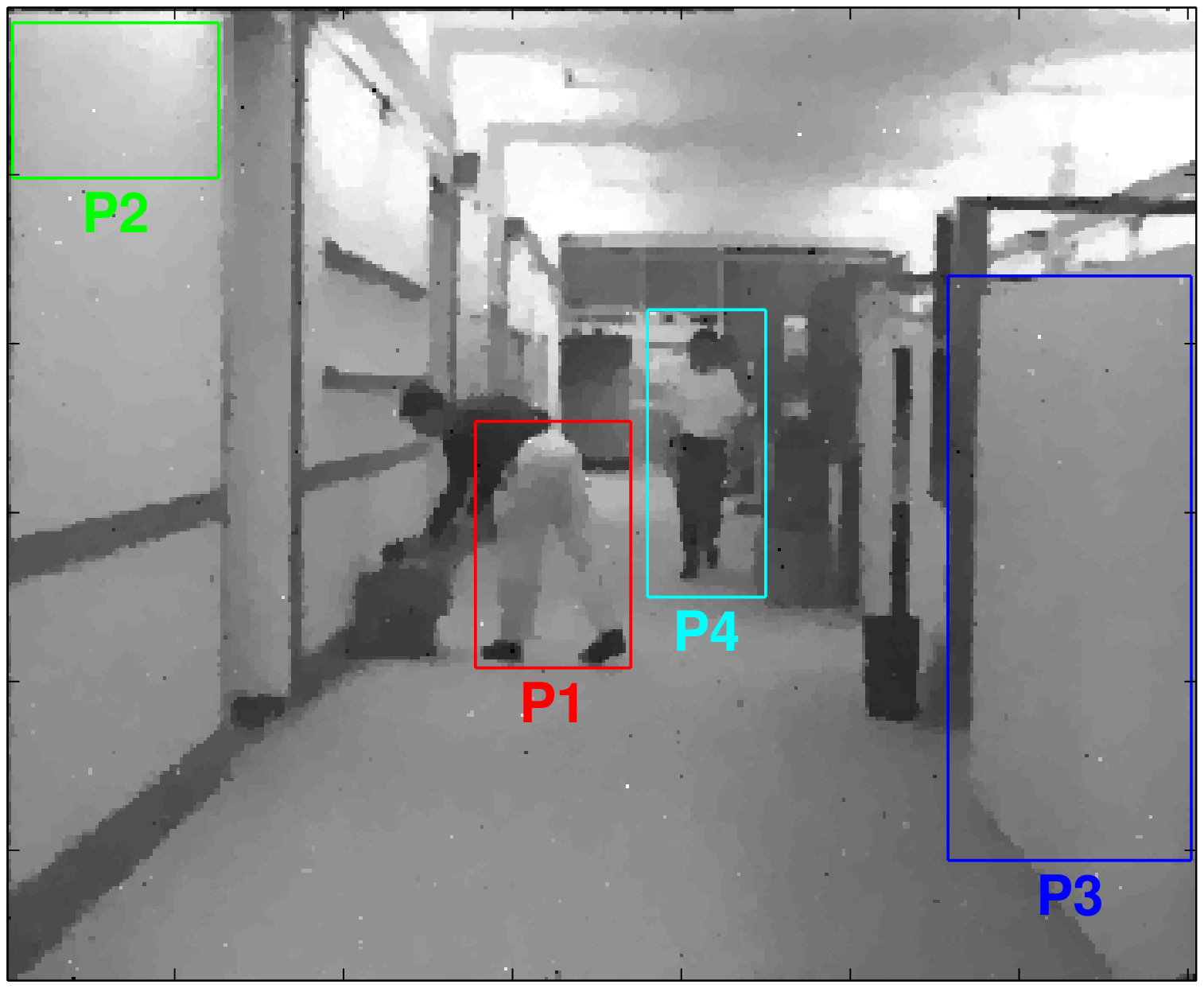}\hspace{-.05in}} &
{\hspace{-.05in}\includegraphics[width=0.185\textwidth]{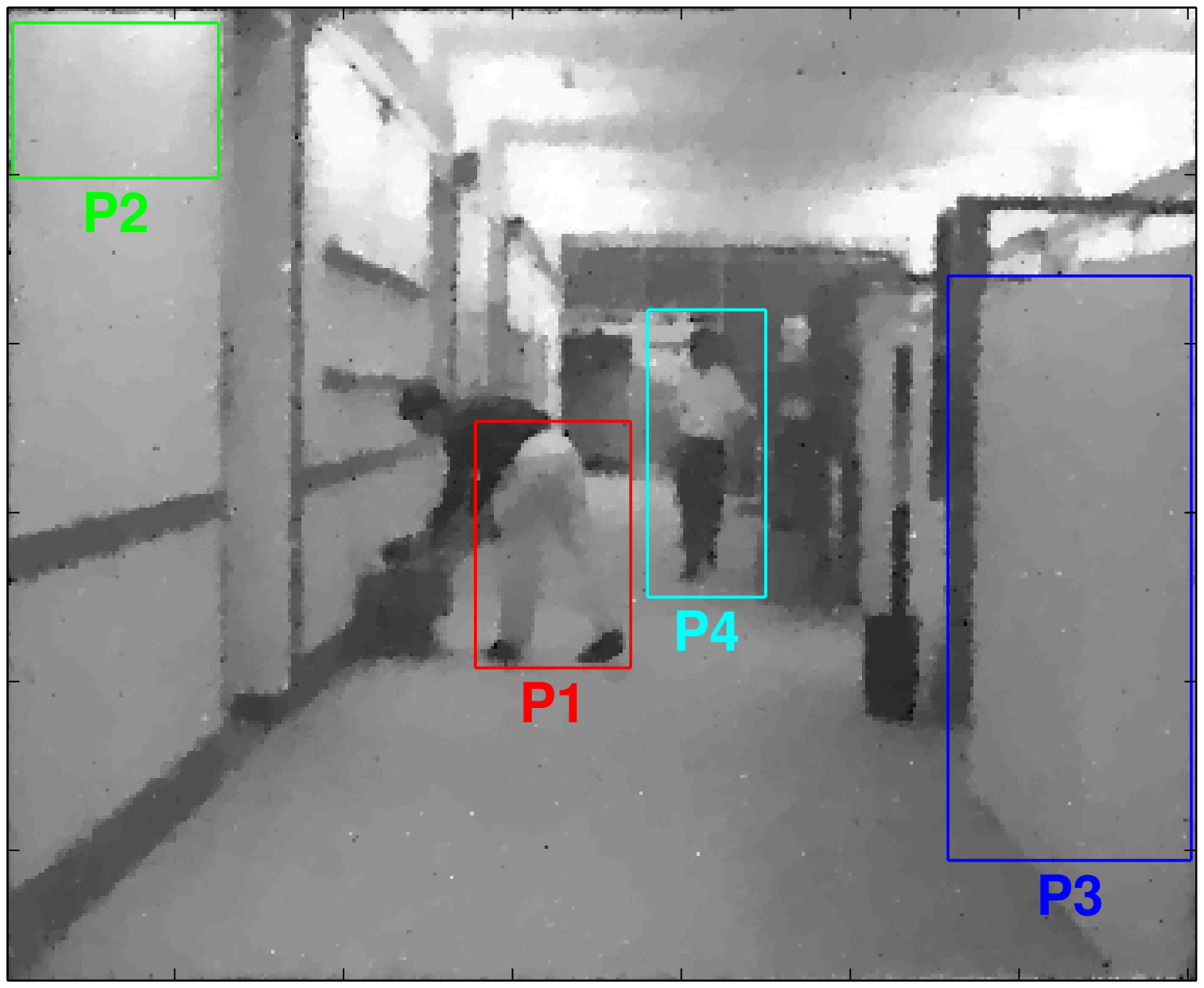}\hspace{-.05in}} &
{\hspace{-.05in}\includegraphics[width=0.185\textwidth]{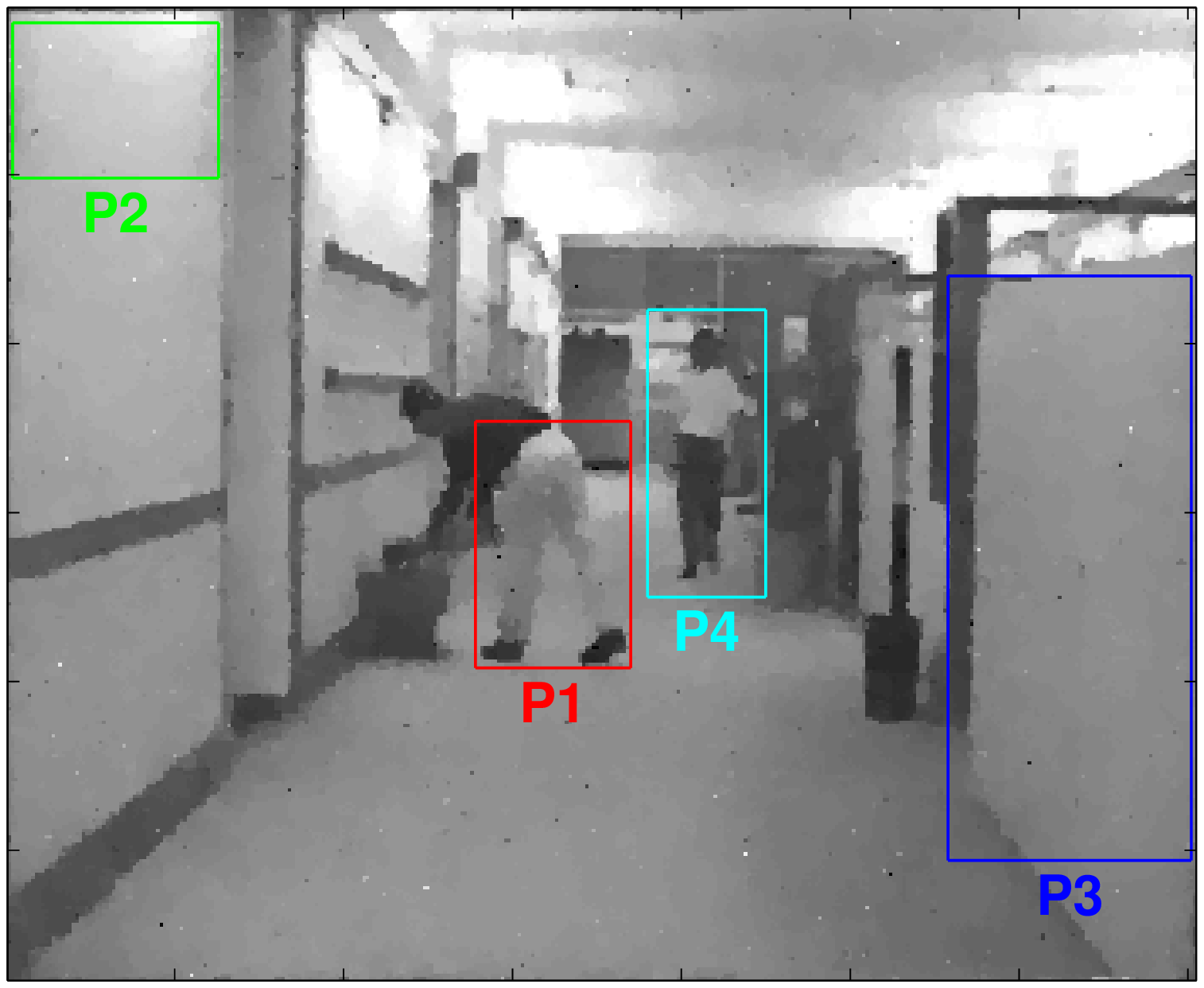}\hspace{-.05in}} \\
& PSNR = 33.30  & PSNR = 31.28  & PSNR = 25.44  & PSNR = 24.79  & PSNR = 26.67\\
\hline
{\hspace{-.05in}\begin{sideways} $|{\bf\Omega}|/mnN=0.01$ \end{sideways}\hspace{-.05in}} &
{\hspace{-.05in}\includegraphics[width=0.185\textwidth]{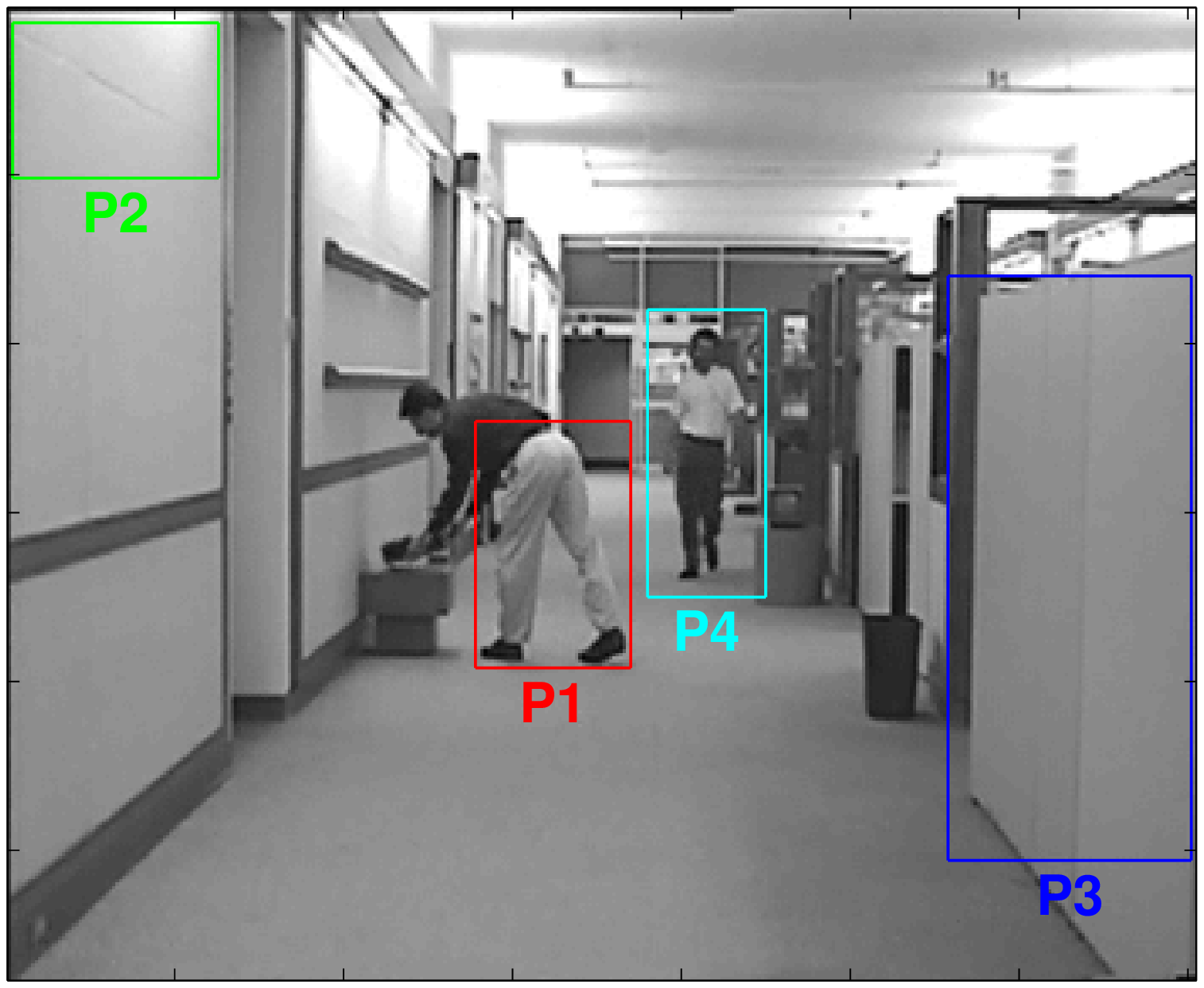}\hspace{-.05in}} &
{\hspace{-.05in}\includegraphics[width=0.185\textwidth]{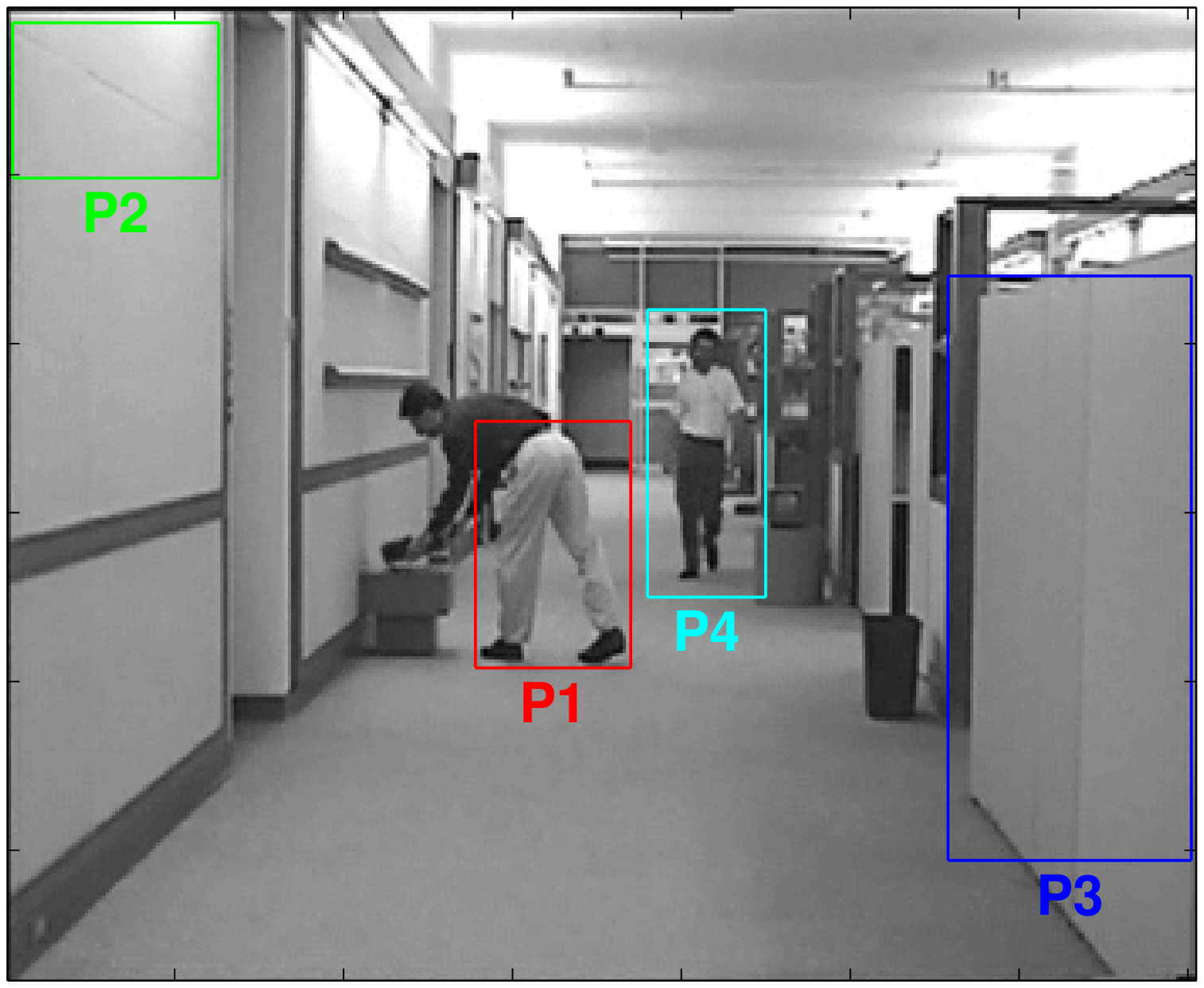}\hspace{-.05in}} &
{\hspace{-.05in}\includegraphics[width=0.185\textwidth]{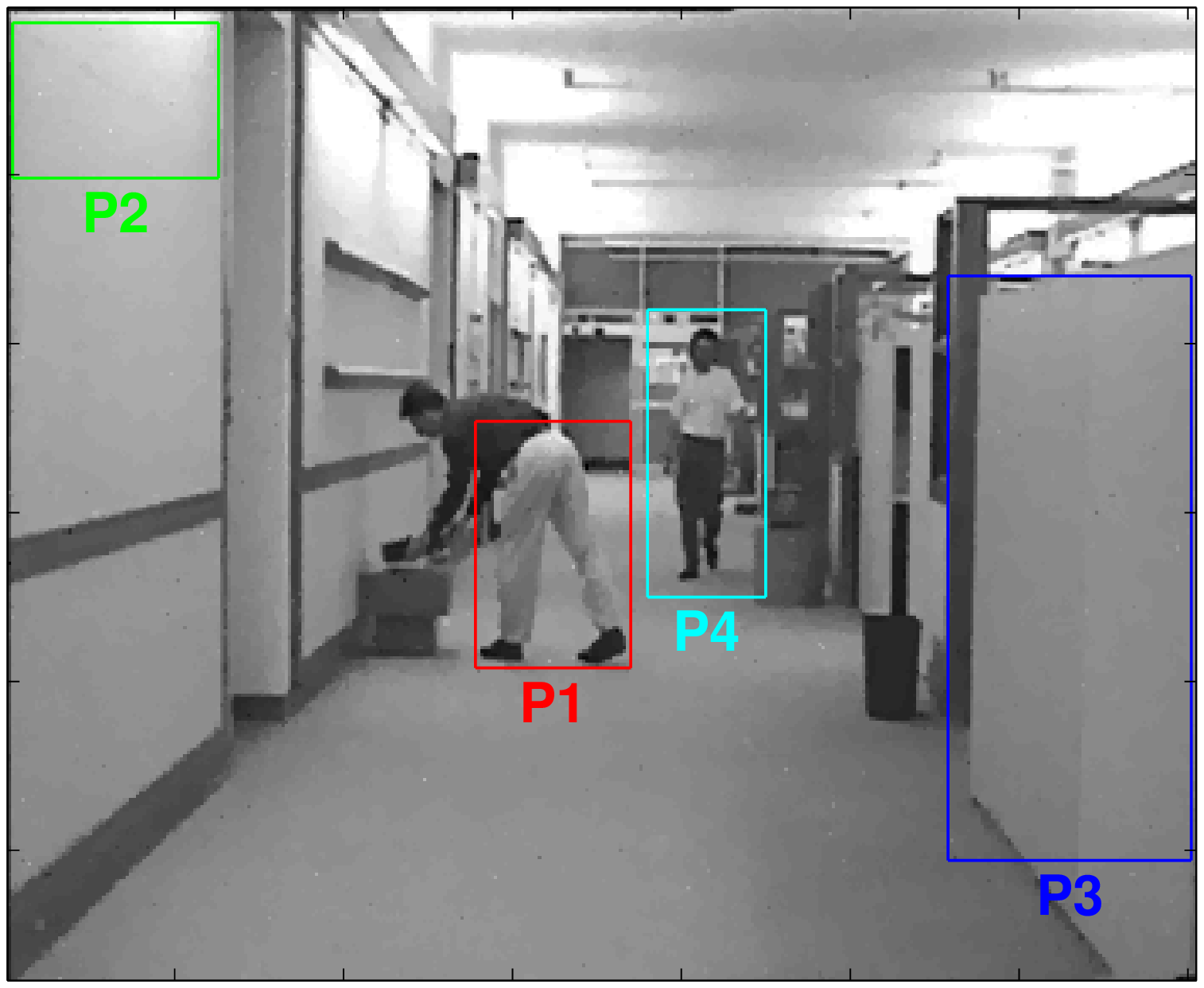}\hspace{-.05in}} &
{\hspace{-.05in}\includegraphics[width=0.185\textwidth]{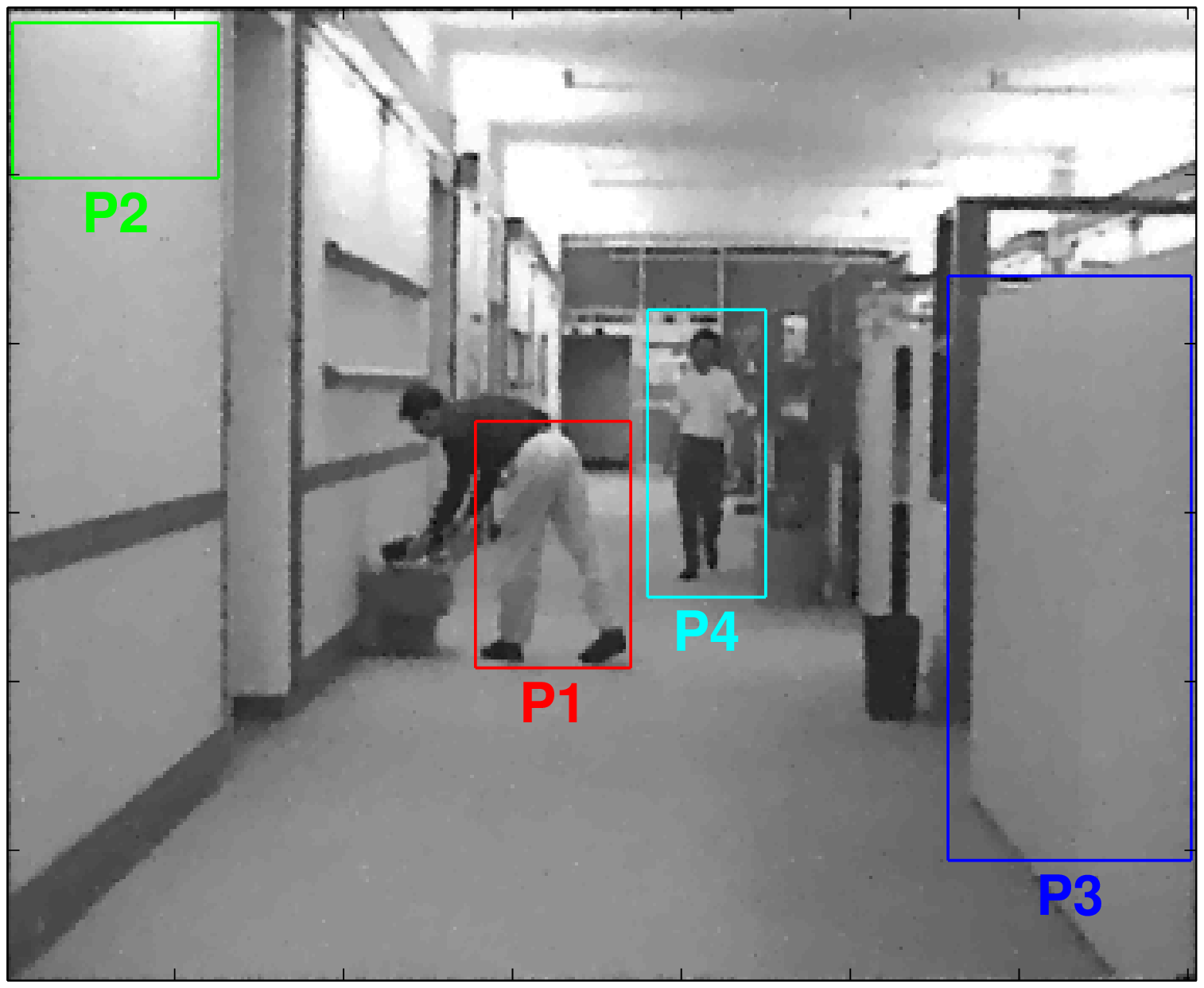}\hspace{-.05in}} &
{\hspace{-.05in}\includegraphics[width=0.185\textwidth]{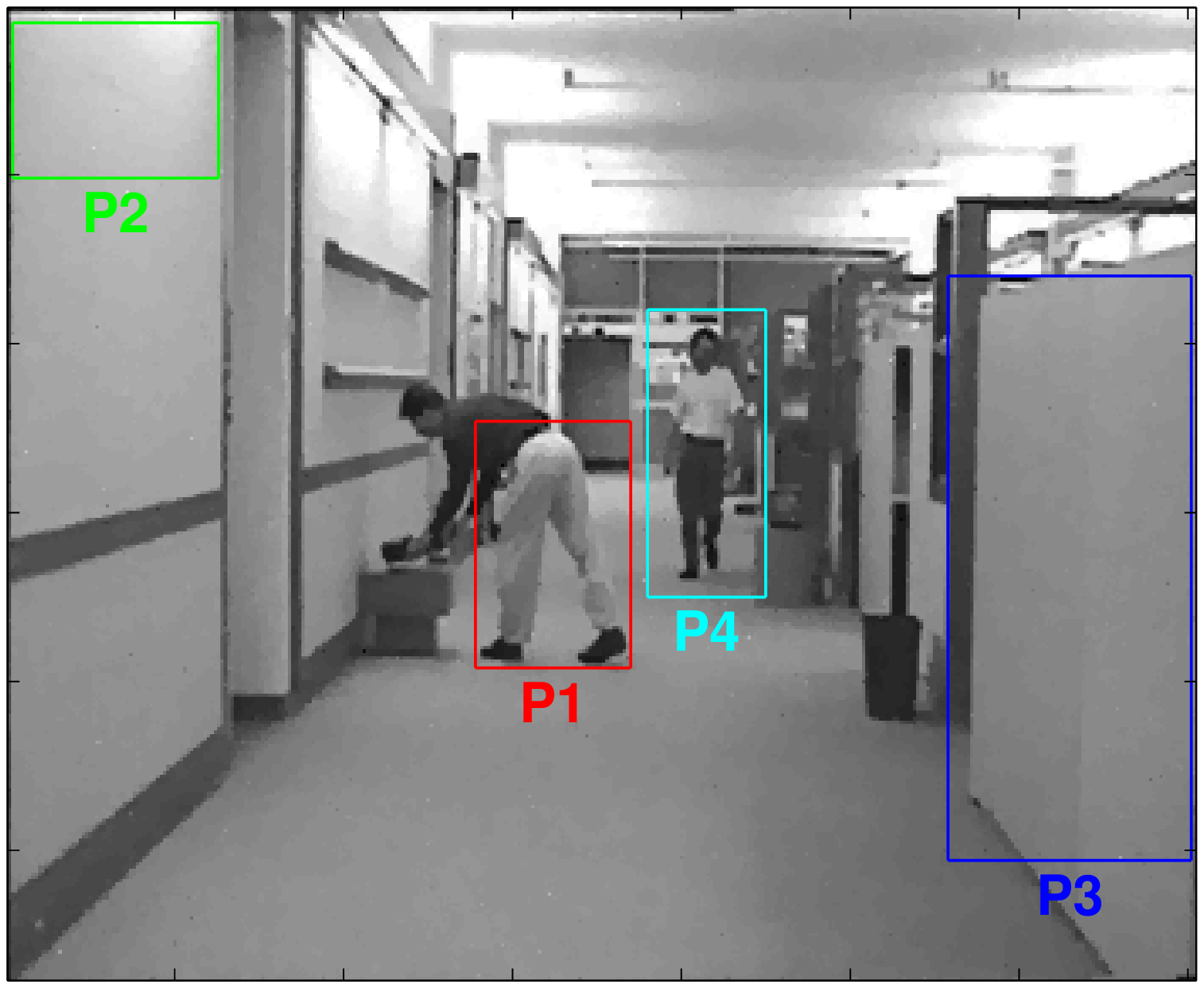}\hspace{-.05in}} \\
& PSNR = 38.27 & PSNR = 37.40 & PSNR = 33.32 & PSNR = 30.09  & PSNR = 33.35  \\
\hline
\end{tabular}
\end{center}
\end{table*}

\begin{table*}
\renewcommand{\arraystretch}{1.3}
\caption{Absolute residual error of reconstructed frame No. $17$ with respect to reference frame.}
\label{Table:FrameIndex17RecoveryHallMonitorResidue}\vspace{-.1in}
\begin{center}
\begin{tabular}{|c|c|c|c|c|c|}
\cline{2-6}
\multicolumn{1}{c|}{} & AR-BC-ADMM & P-BC-ADMM & TFOCS & NESTA & TVAL3 \\
\hline
{\hspace{-.05in}\begin{sideways} $|{\bf\Omega}|/mnN=0.01$ \end{sideways}\hspace{-.05in}} &
{\hspace{-.05in}\includegraphics[width=0.185\textwidth]{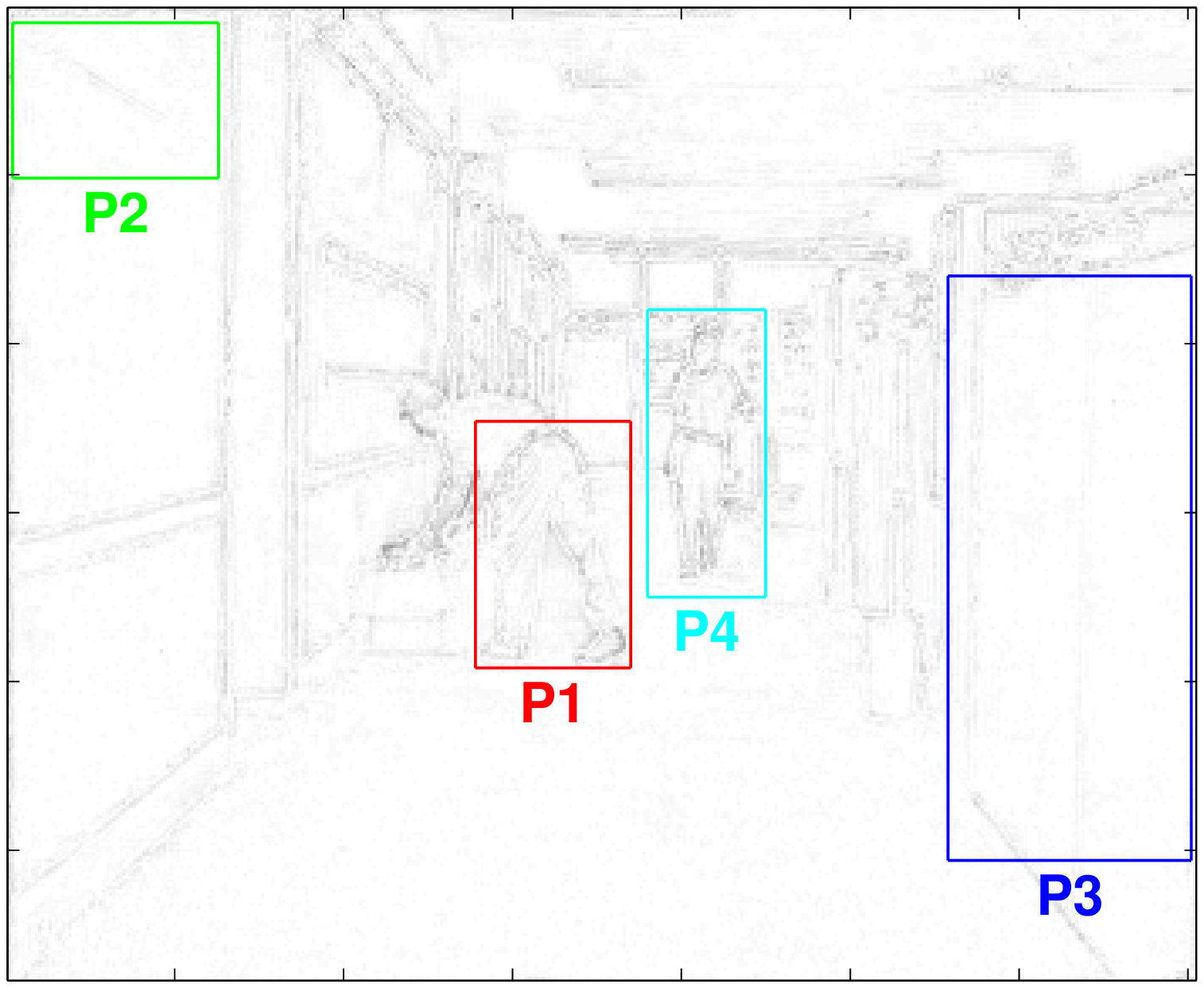}\hspace{-.05in}} &
{\hspace{-.05in}\includegraphics[width=0.185\textwidth]{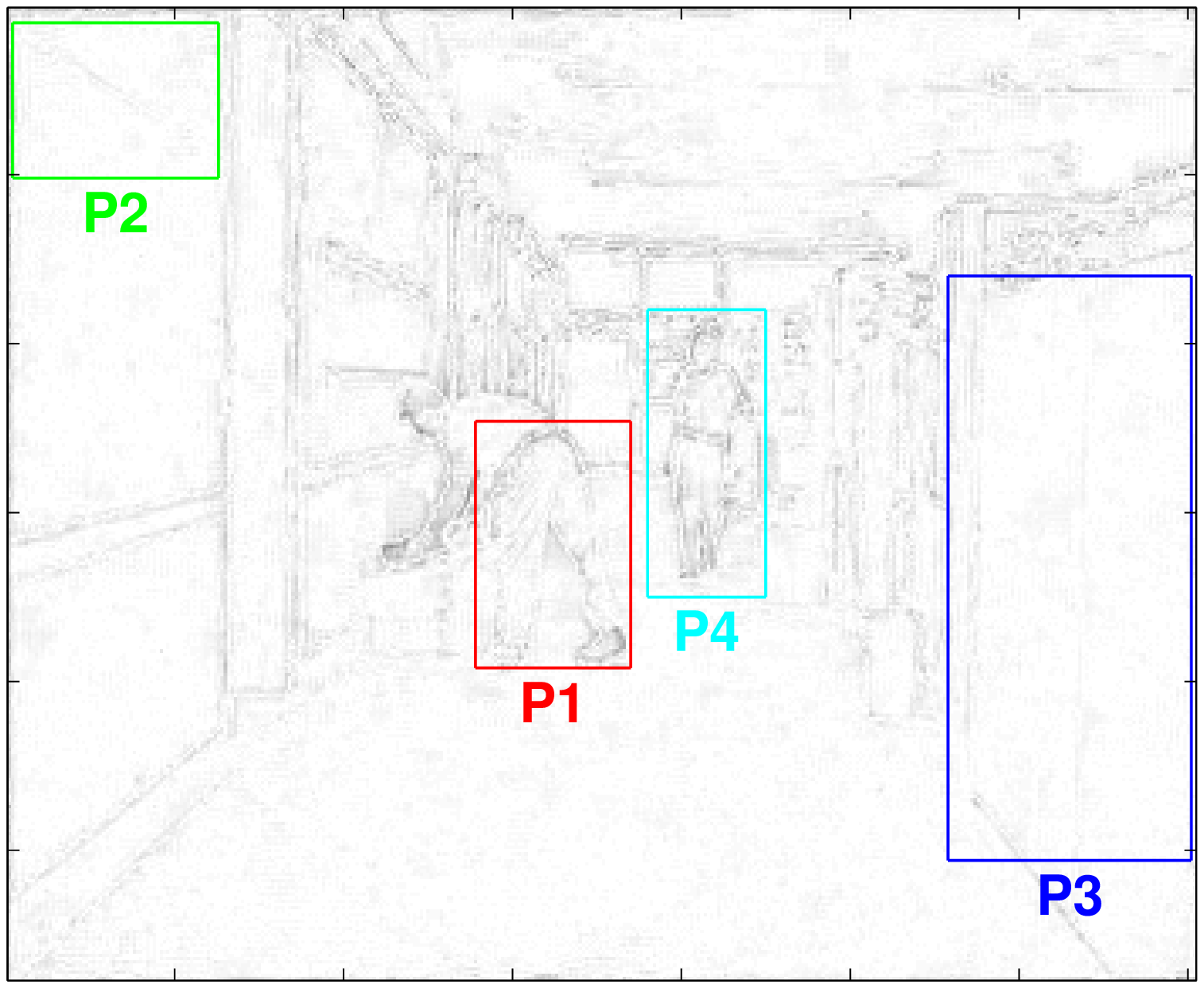}\hspace{-.05in}} &
{\hspace{-.05in}\includegraphics[width=0.185\textwidth]{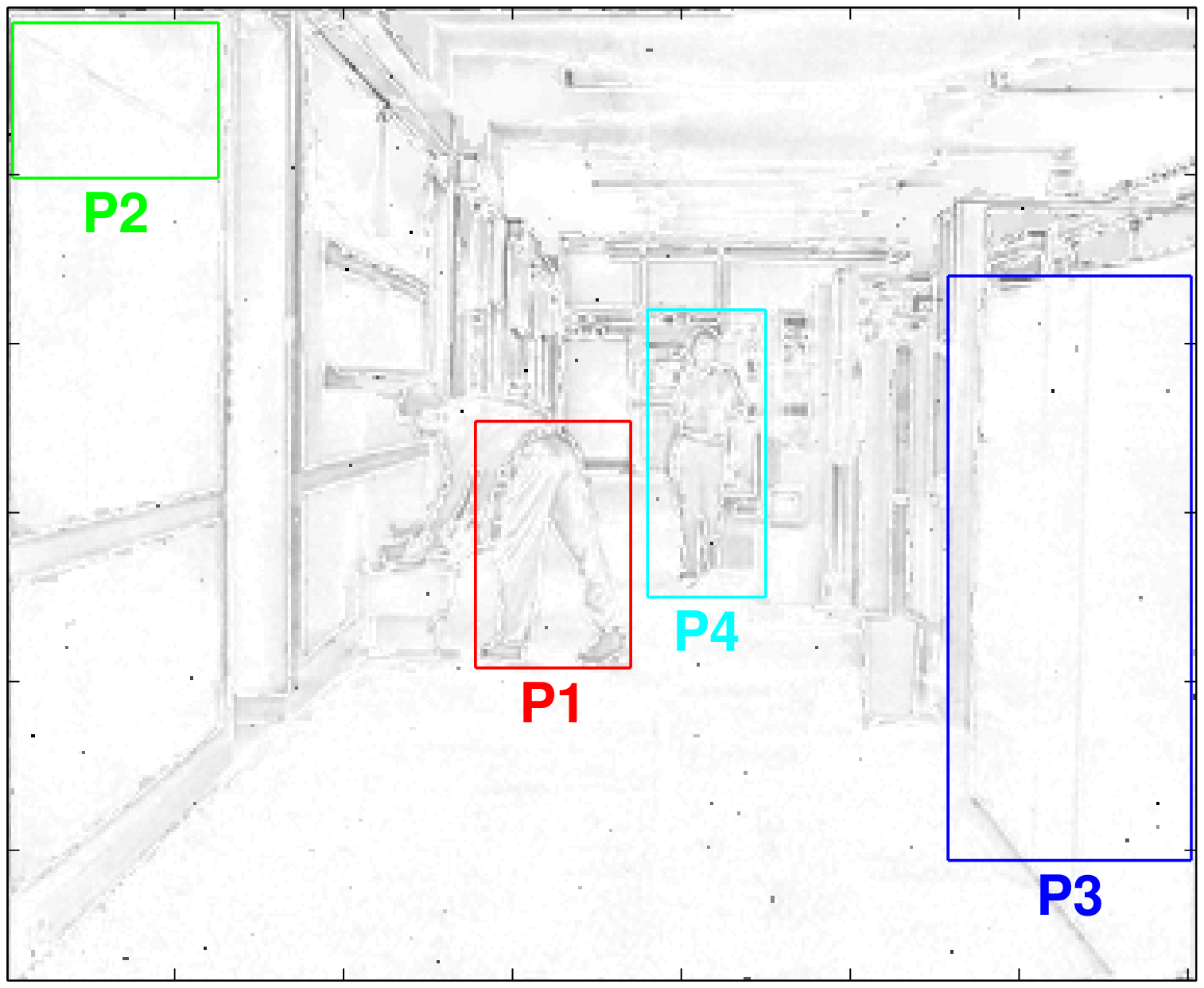}\hspace{-.05in}} &
{\hspace{-.05in}\includegraphics[width=0.185\textwidth]{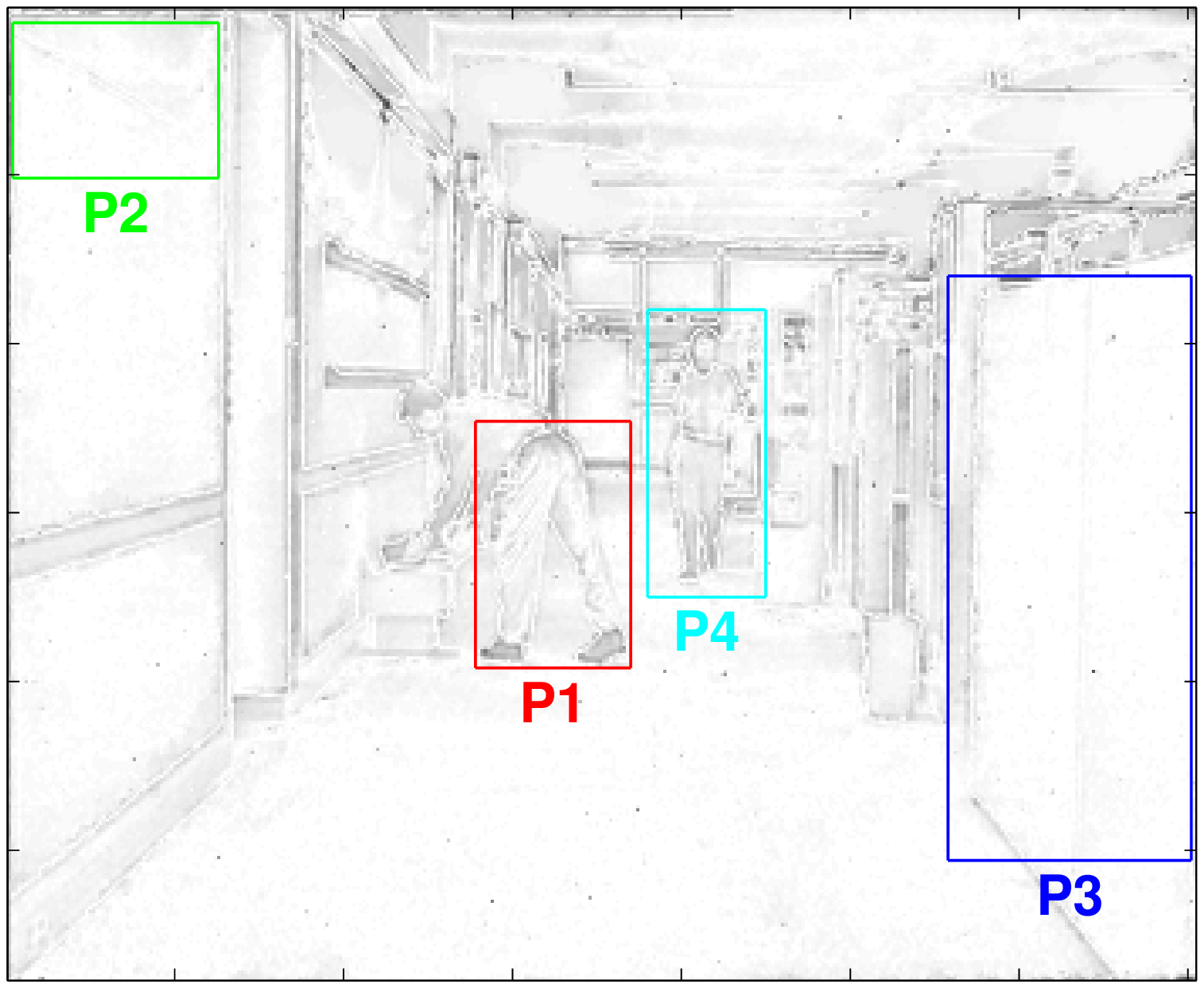}\hspace{-.05in}} &
{\hspace{-.05in}\includegraphics[width=0.185\textwidth]{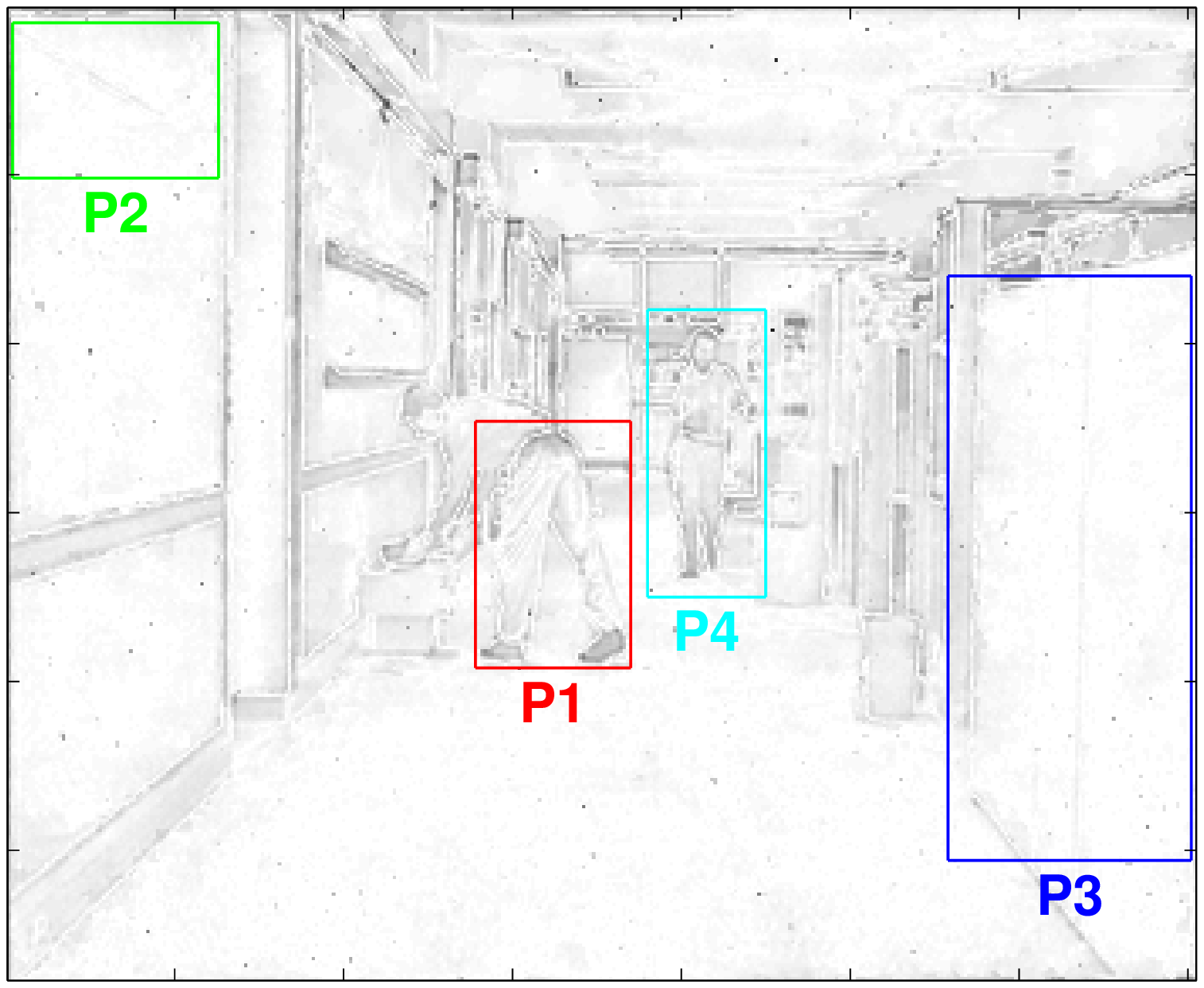}\hspace{-.05in}} \\
& NMSE = 0.0669  & NMSE = 0.0842 & NMSE = 0.1285 & NMSE = 0.1367 & NMSE = 0.1424 \\
\hline
{\hspace{-.05in}\begin{sideways} $|{\bf\Omega}|/mnN=0.03$ \end{sideways}\hspace{-.05in}} &
{\hspace{-.05in}\includegraphics[width=0.185\textwidth]{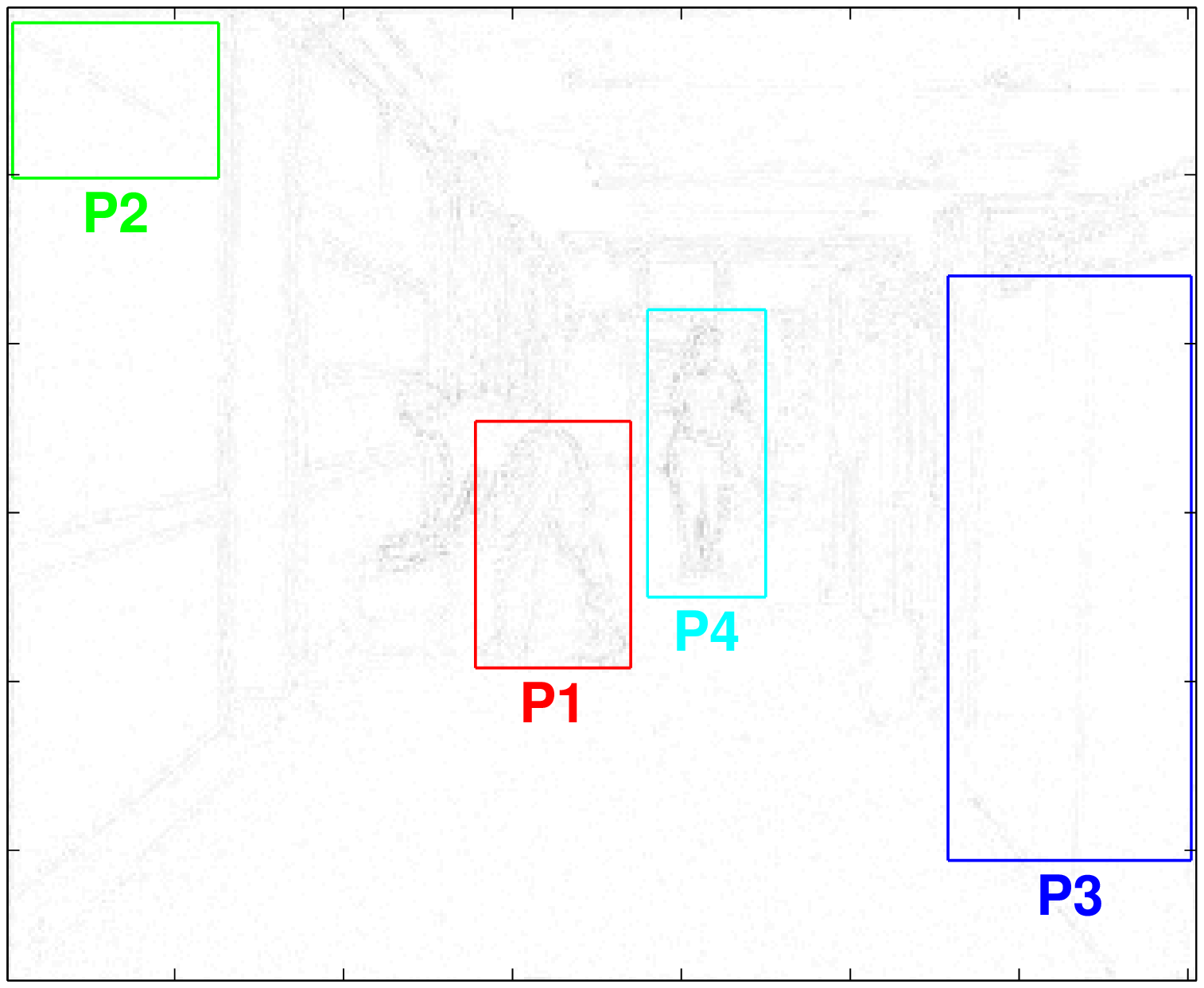}\hspace{-.05in}} &
{\hspace{-.05in}\includegraphics[width=0.185\textwidth]{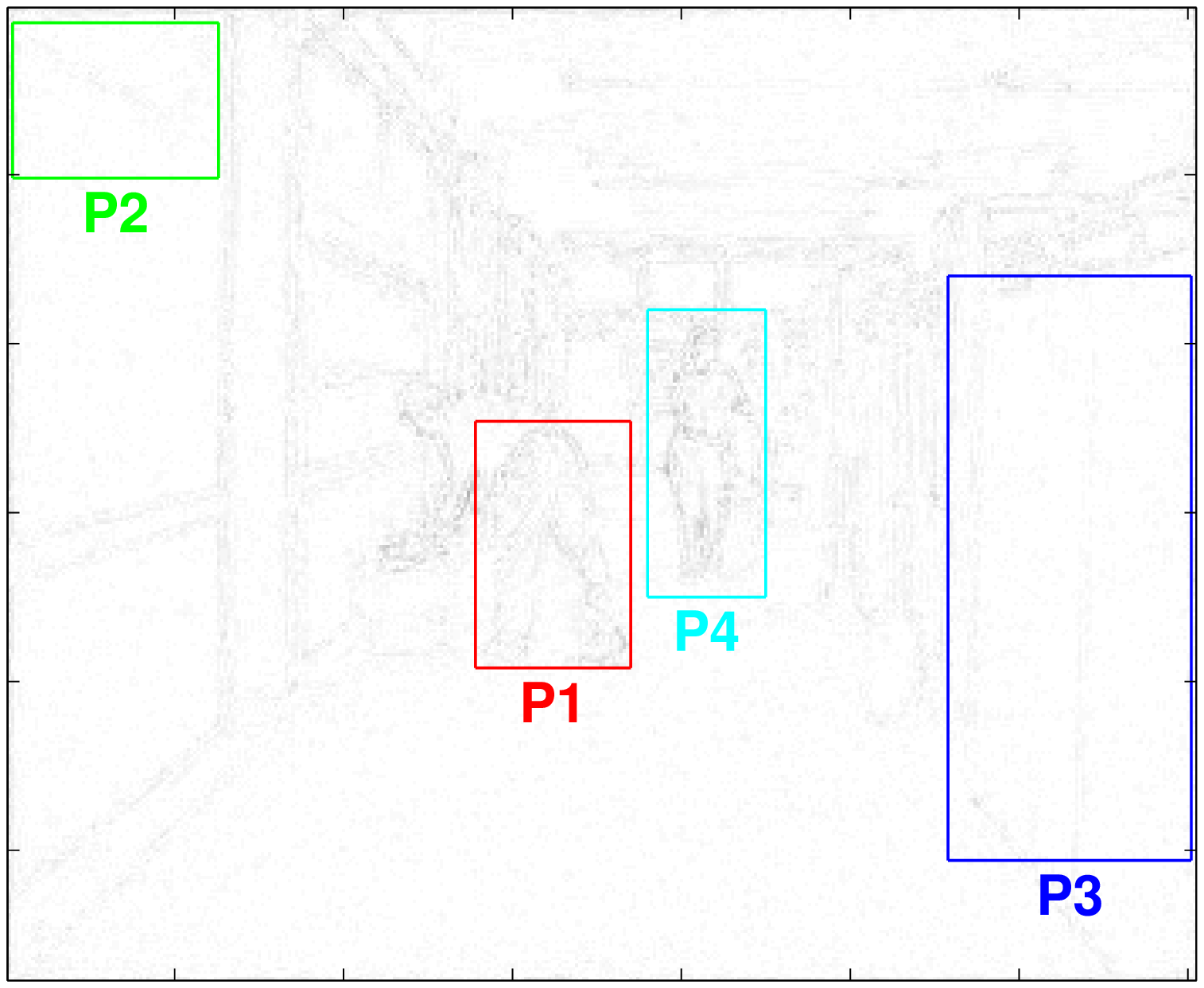}\hspace{-.05in}} &
{\hspace{-.05in}\includegraphics[width=0.185\textwidth]{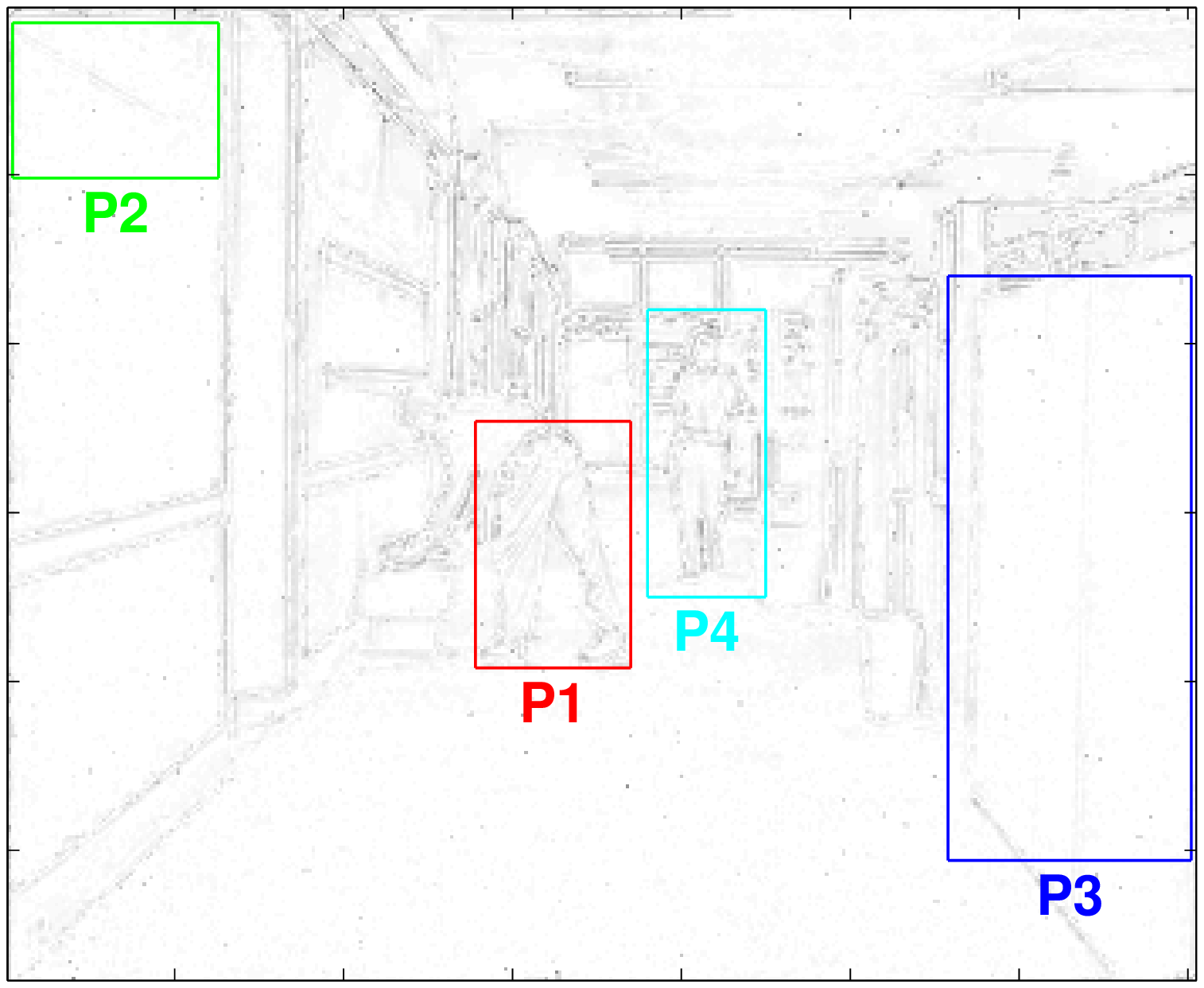}\hspace{-.05in}} &
{\hspace{-.05in}\includegraphics[width=0.185\textwidth]{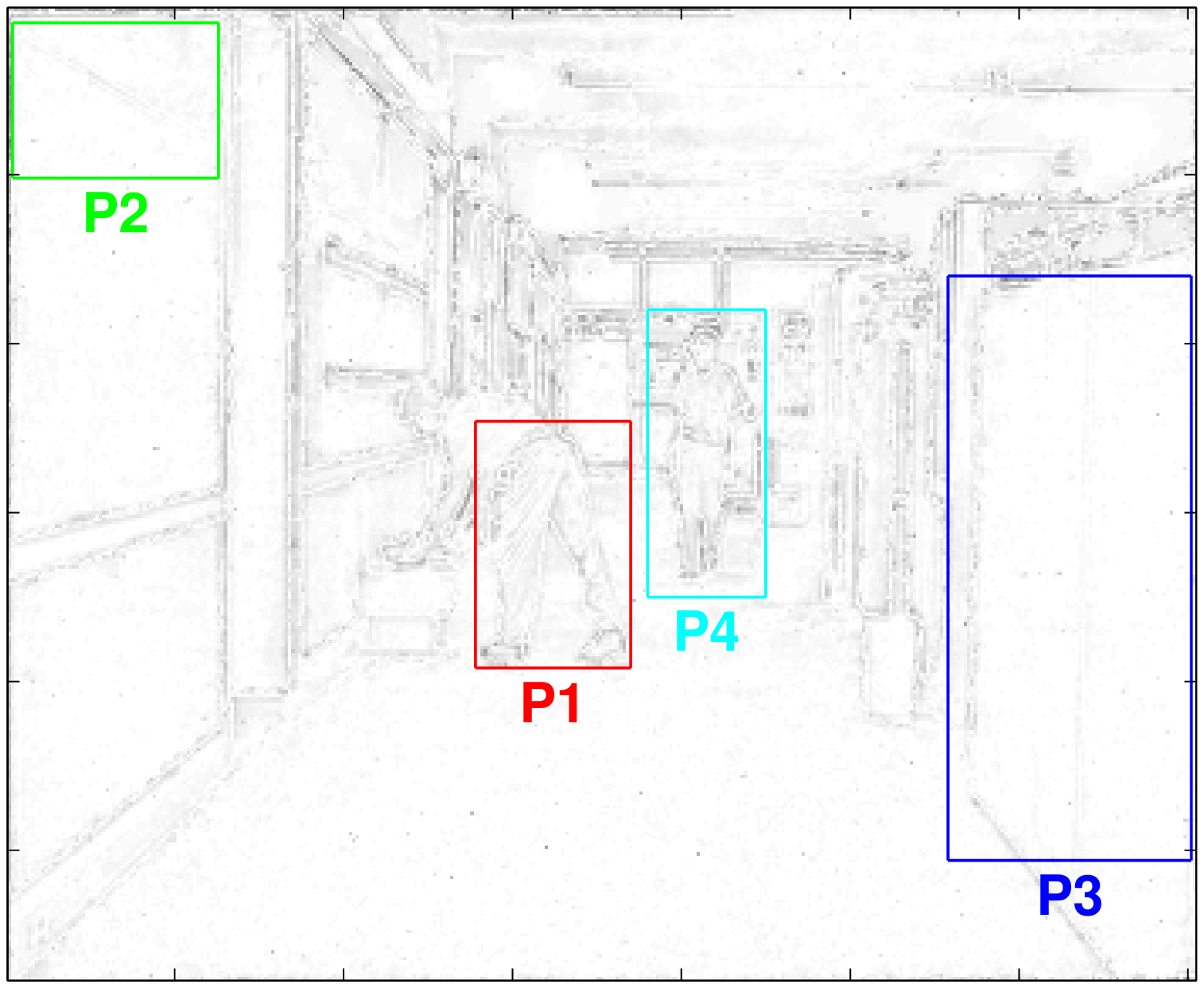}\hspace{-.05in}} &
{\hspace{-.05in}\includegraphics[width=0.185\textwidth]{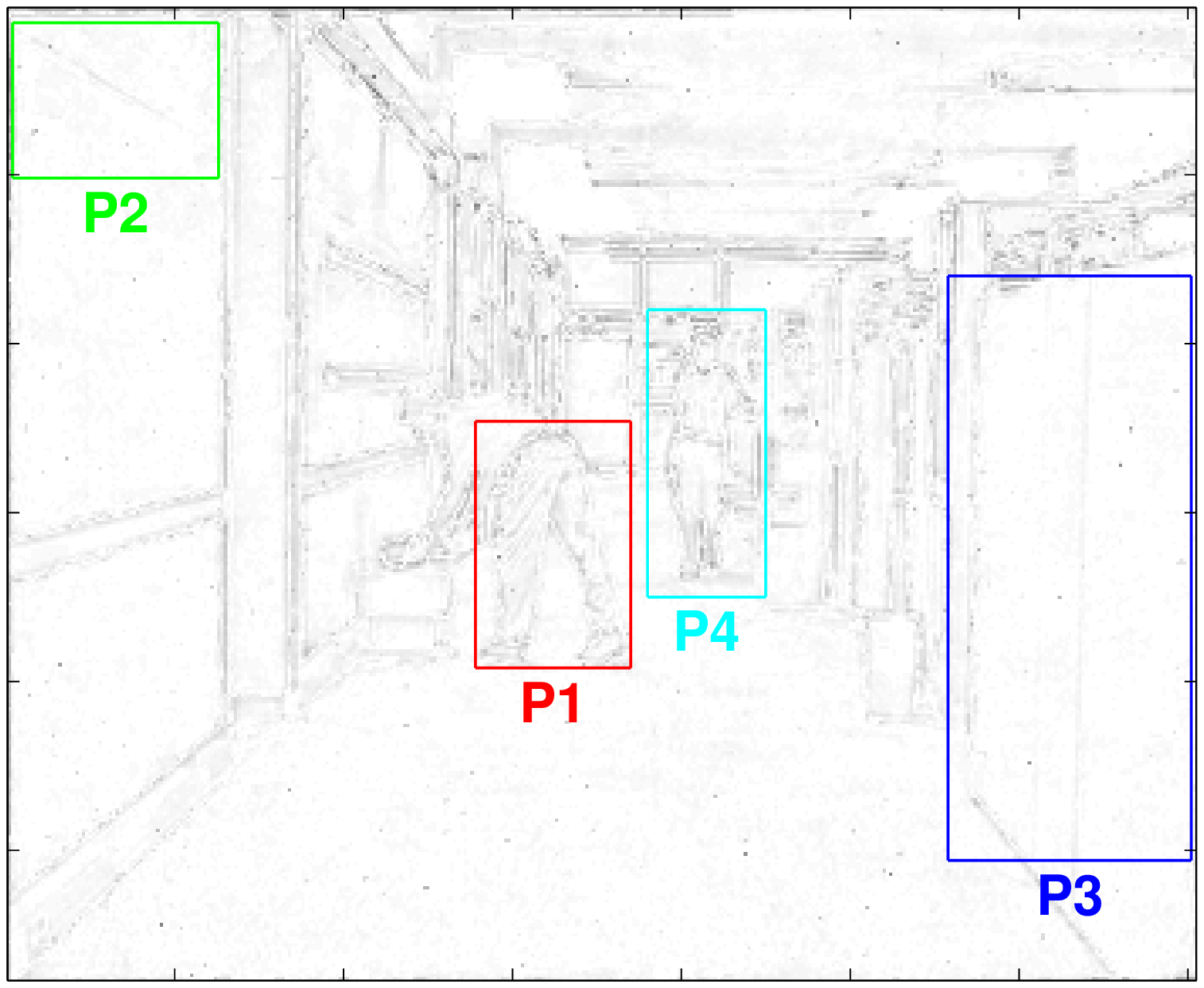}\hspace{-.05in}} \\
& NMSE = 0.0365  & NMSE = 0.0460  & NMSE = 0.0782  & NMSE = 0.0971  & NMSE = 0.0527\\
\hline
{\hspace{-.05in}\begin{sideways} $|{\bf\Omega}|/mnN=0.01$ \end{sideways}\hspace{-.05in}} &
{\hspace{-.05in}\includegraphics[width=0.185\textwidth]{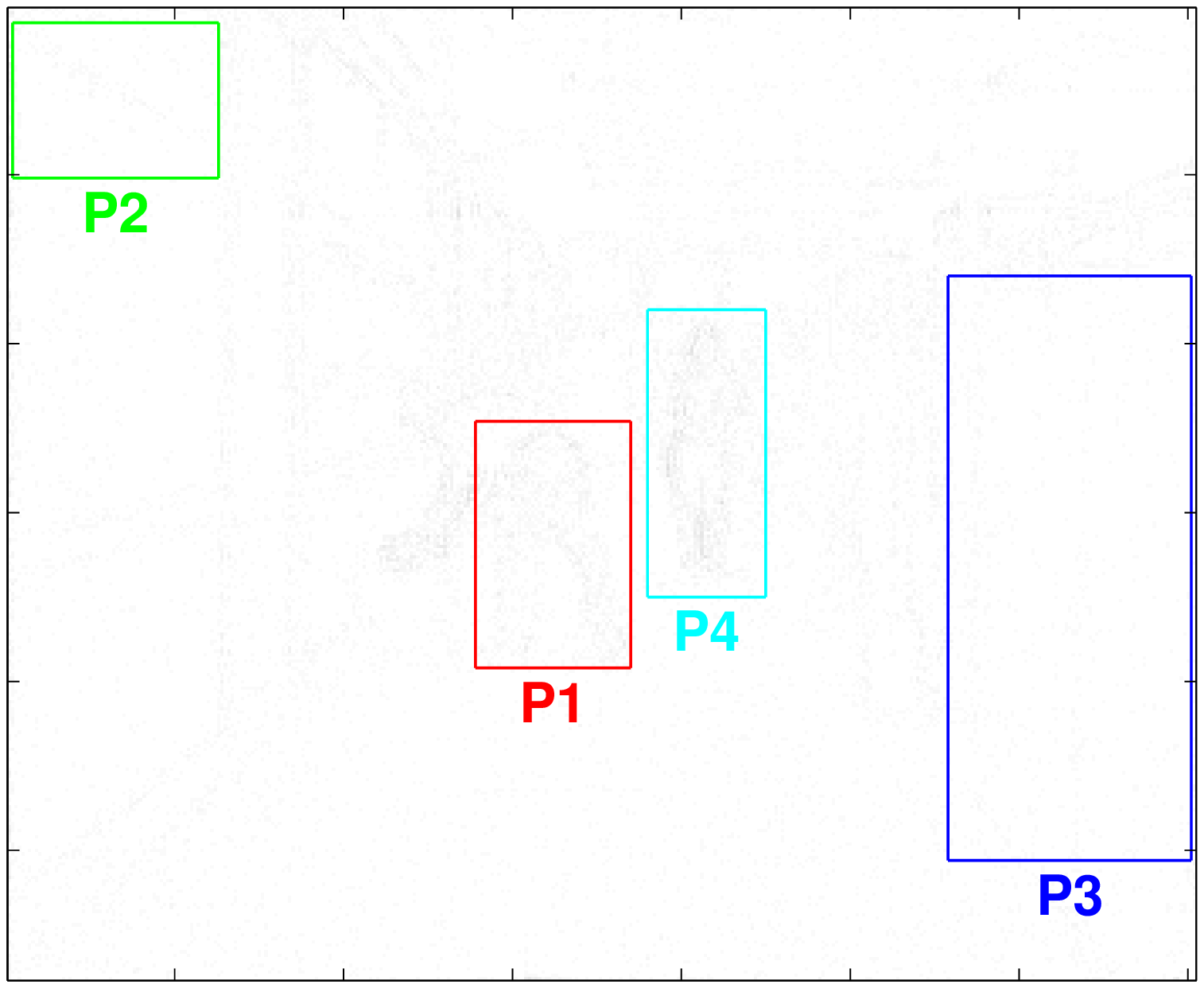}\hspace{-.05in}} &
{\hspace{-.05in}\includegraphics[width=0.185\textwidth]{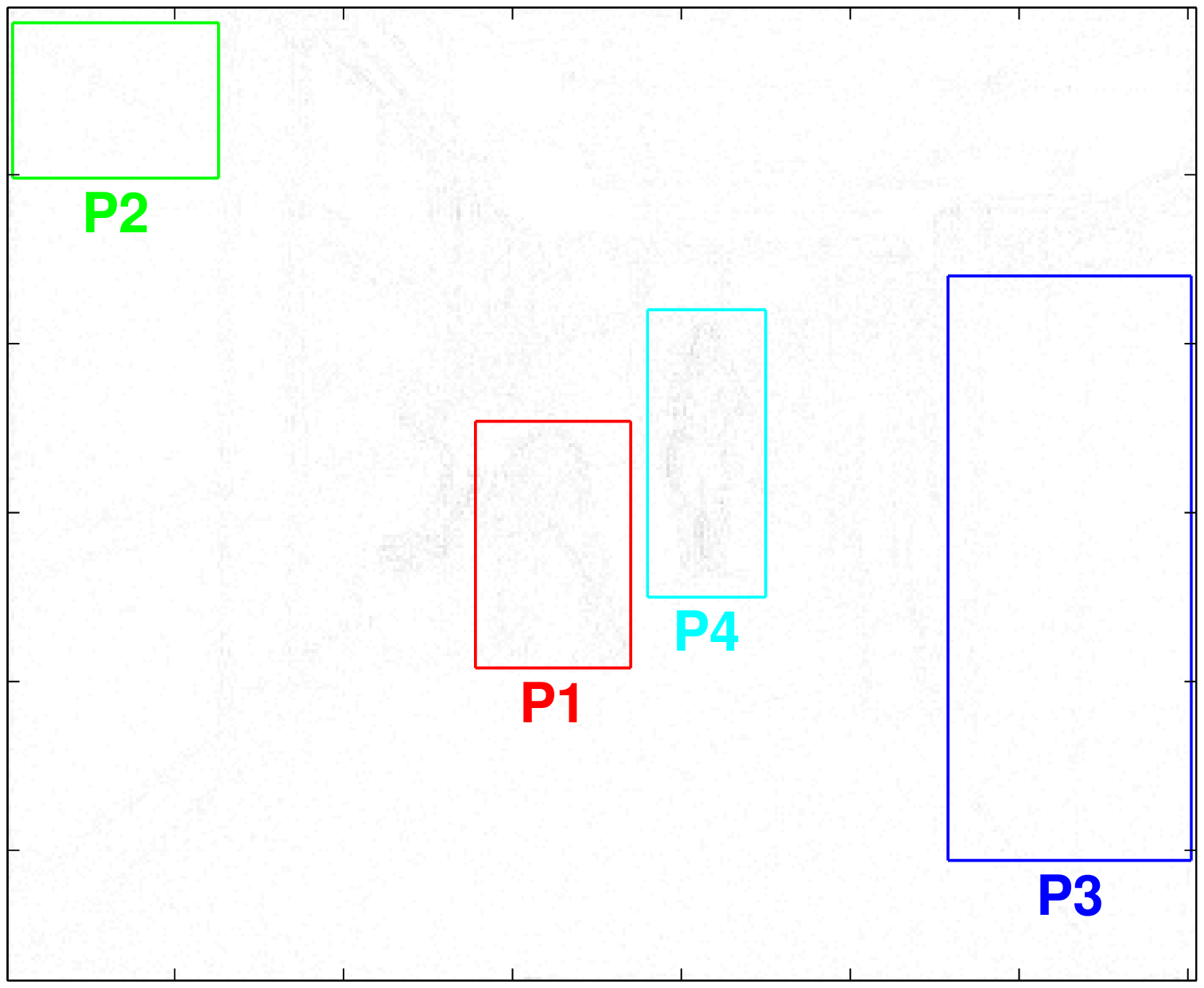}\hspace{-.05in}} &
{\hspace{-.05in}\includegraphics[width=0.185\textwidth]{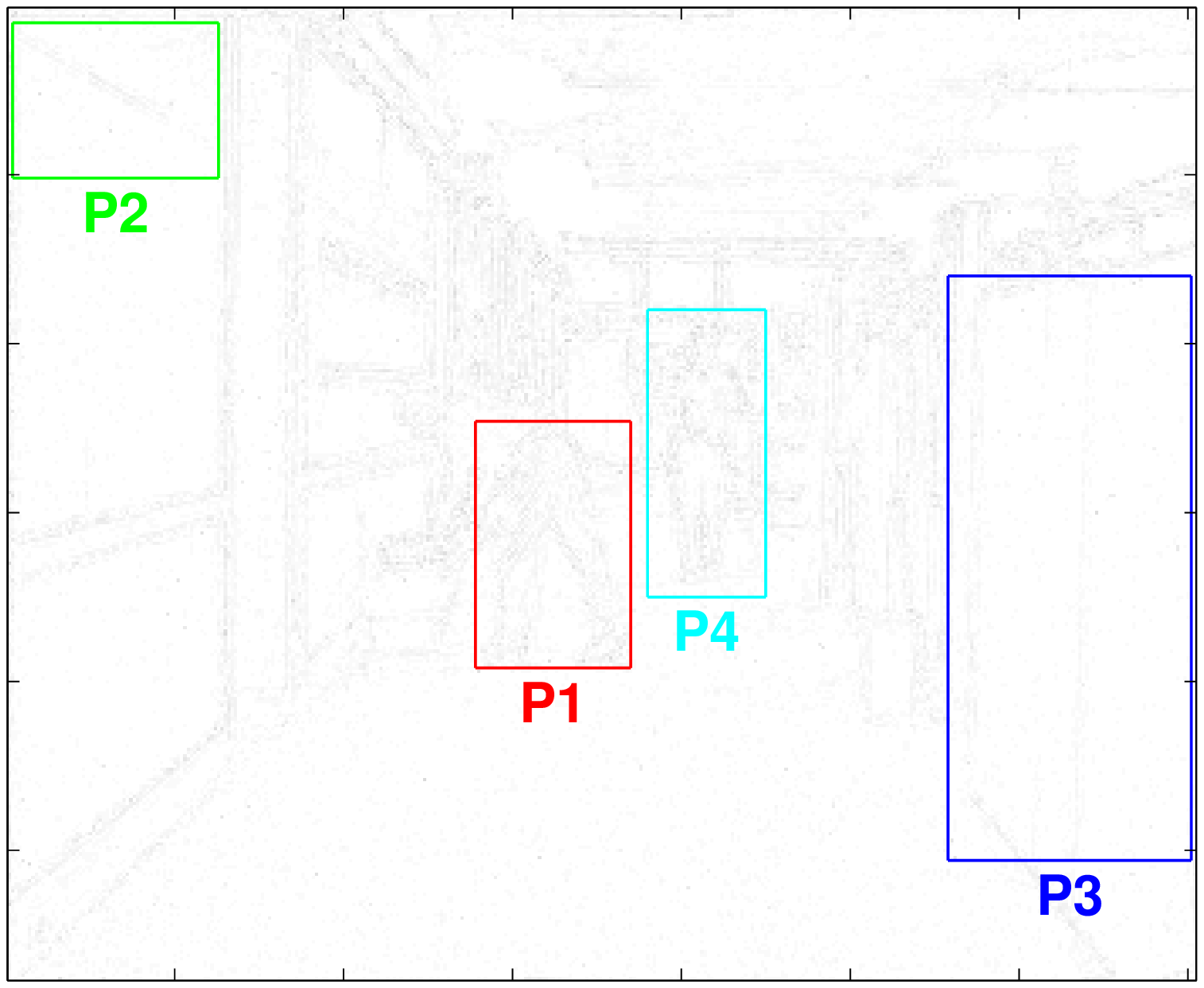}\hspace{-.05in}} &
{\hspace{-.05in}\includegraphics[width=0.185\textwidth]{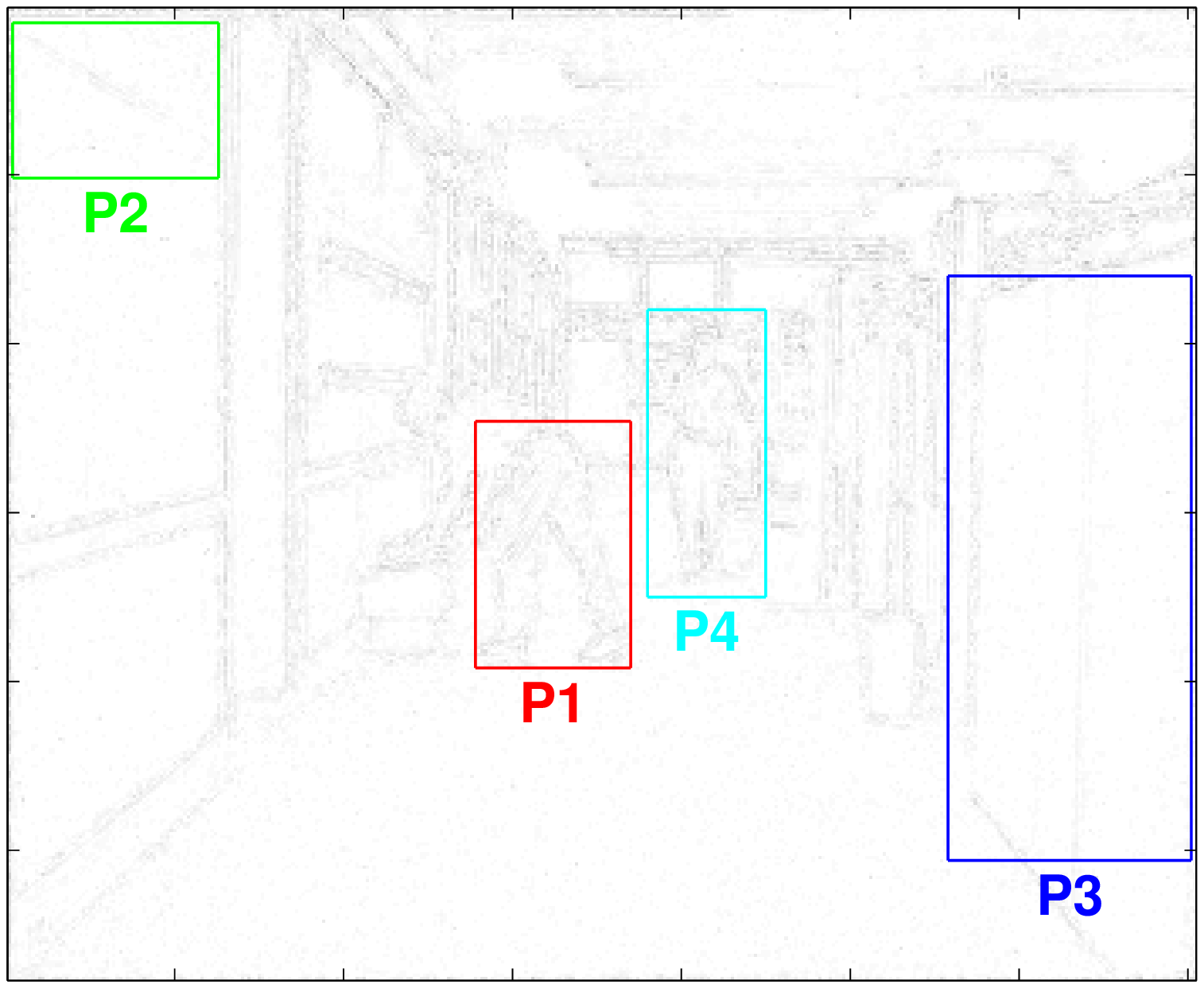}\hspace{-.05in}} &
{\hspace{-.05in}\includegraphics[width=0.185\textwidth]{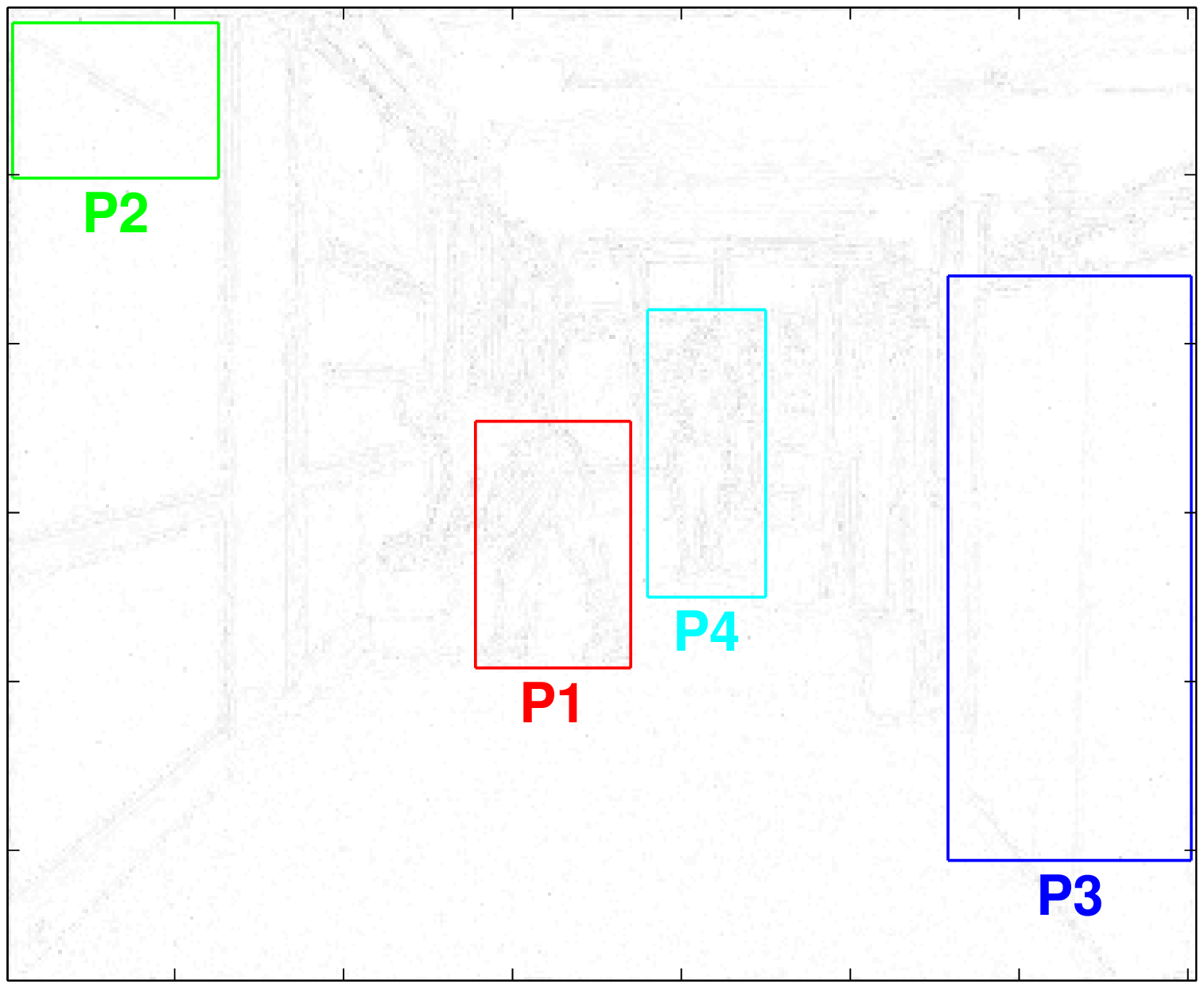}\hspace{-.05in}} \\
& NMSE = 0.0206  & NMSE = 0.0227 & NMSE = 0.0362 & NMSE = 0.0901  & NMSE = 0.0364  \\
\hline
\end{tabular}
\end{center}
\end{table*}

\vspace{-.1in}
\subsection{Frame Index (Boundary) Analysis}\label{sec_frame_index_analysis}
The second experiment is performed by freezing the sampling rate and evaluate the reconstruction accuracy of each individual frame ($32$ frames) in the test videos. The main goal is to study the effect of different BCs on our proposed method and compare the results with the existing state-of-the-art algorithms. Table \ref{Table:FrameIndex} elaborates on the results of using two different sampling rates fixed at $1\%$ and $10\%$ on four video test clips. The results state the efficiency of AR-BC-ADMM which remains superior with respect P-BC-ADMM and rest of the competing algorithms. Also, all three competing algorithms retain much lower accuracy compared to AR-BC-ADMM and P-BC-ADMM as a result of employing low accuracy differential scheme. As expected, P-BC-ADMM at high sampling rates, weakly performs on boundary frames. Such a weak performance is mostly noticeable at the portions with motion trajectories. This is because of using larger size of filter length for temporal differentiation. Such a kernel size has a greater impact on introducing artificial temporal discontinuity during periodization and severely affects the accuracy of reconstruction at the boundary frames. This issue is greatly enhanced by employing AR-BC-ADMM, where the continuity of the temporal derivatives are preserved accordingly.

\begin{table*}
\renewcommand{\arraystretch}{1.3}
\caption{Frame index recovery of BC-ADMM, TFOCS, NESTA and TVAL3 via different sampling ratios performed on four databases}
\label{Table:FrameIndex}\vspace{-.1in}
\begin{center}
\begin{tabular}{|c|c|c|c|c|}
\cline{2-5}
\multicolumn{1}{c|}{} & \textit{Hall-Monitor} & \textit{Container} & \textit{Office Environment} & \textit{Squash2} \\
\hline
{\hspace{-.05in}\begin{sideways} \hspace{.3in}$|{\bf\Omega}|/mnN=0.01$ \end{sideways}\hspace{-.05in}} &
{\hspace{-.05in}\includegraphics[width=0.2\textwidth]{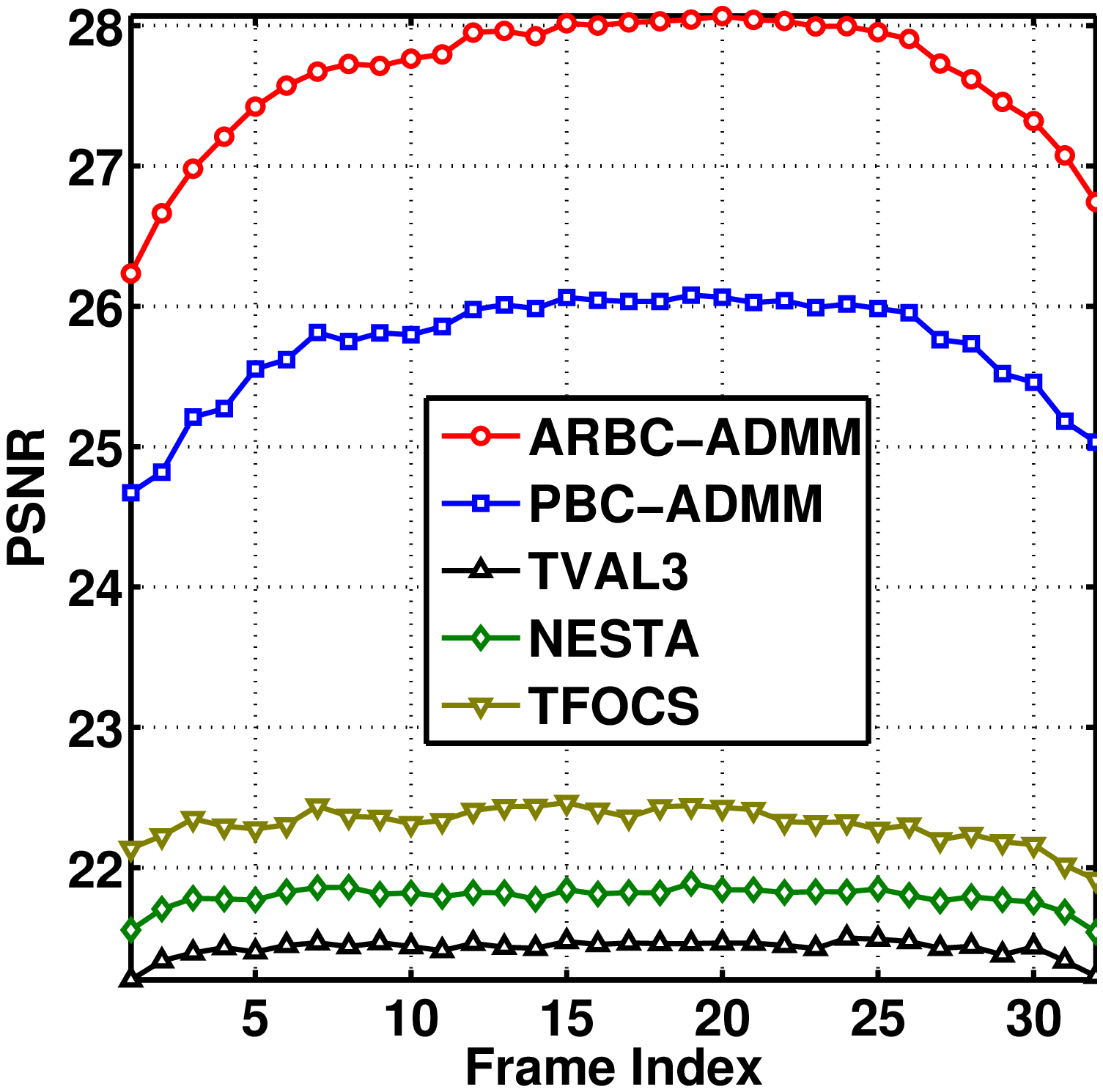}\hspace{-.05in}} &
{\hspace{-.05in}\includegraphics[width=0.2\textwidth]{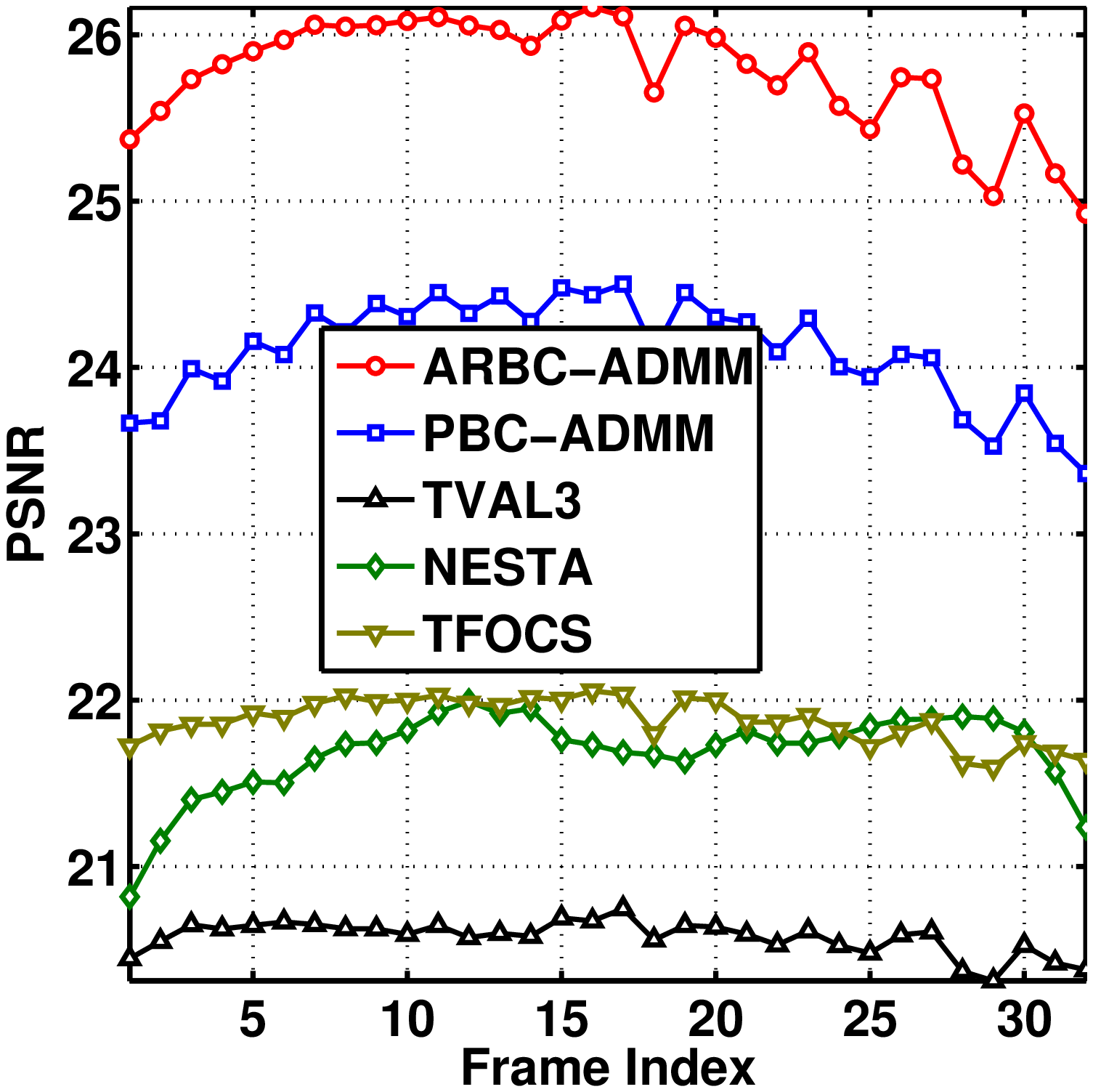}\hspace{-.05in}} &
{\hspace{-.05in}\includegraphics[width=0.2\textwidth]{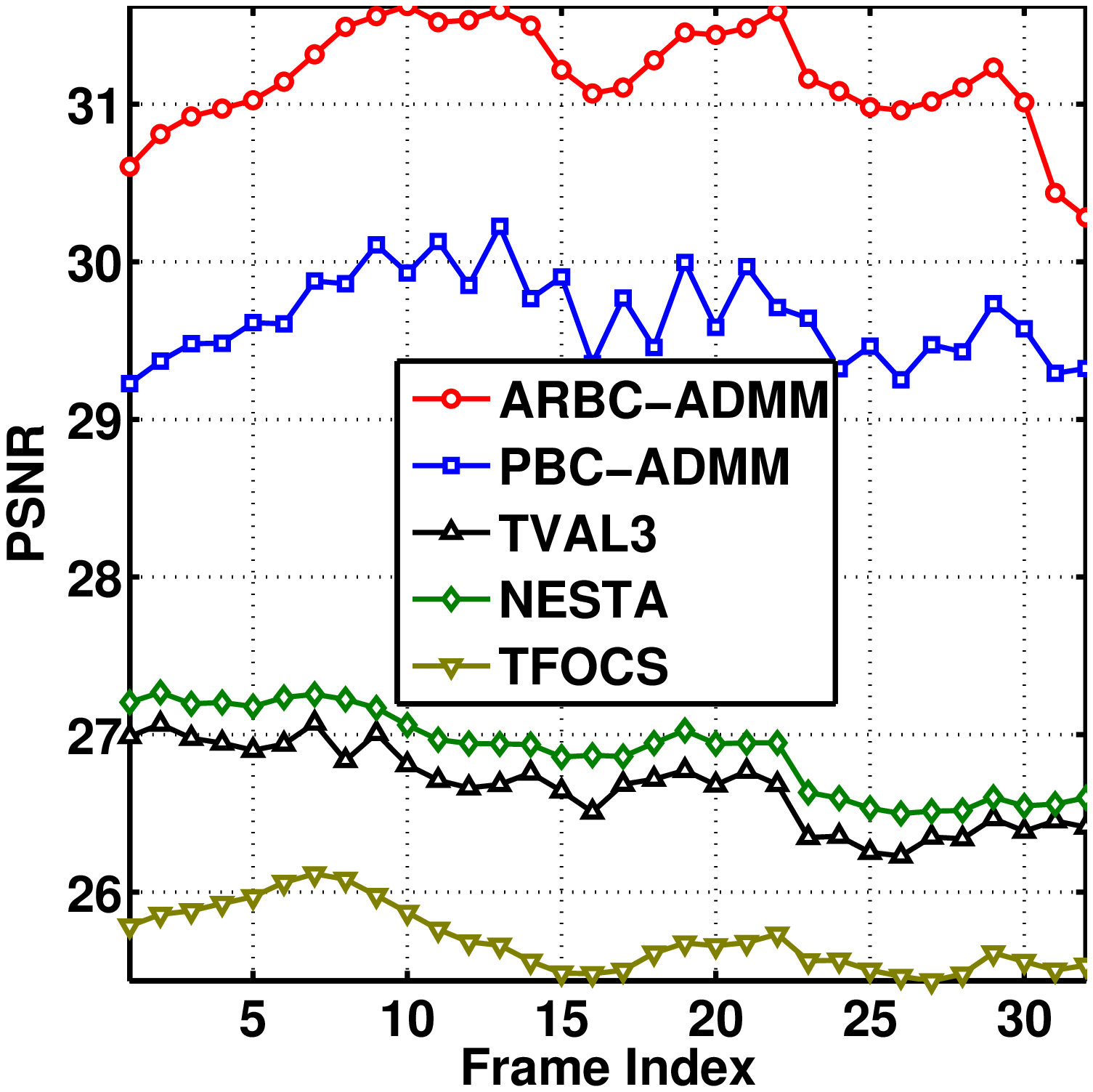}\hspace{-.05in}} &
{\hspace{-.05in}\includegraphics[width=0.2\textwidth]{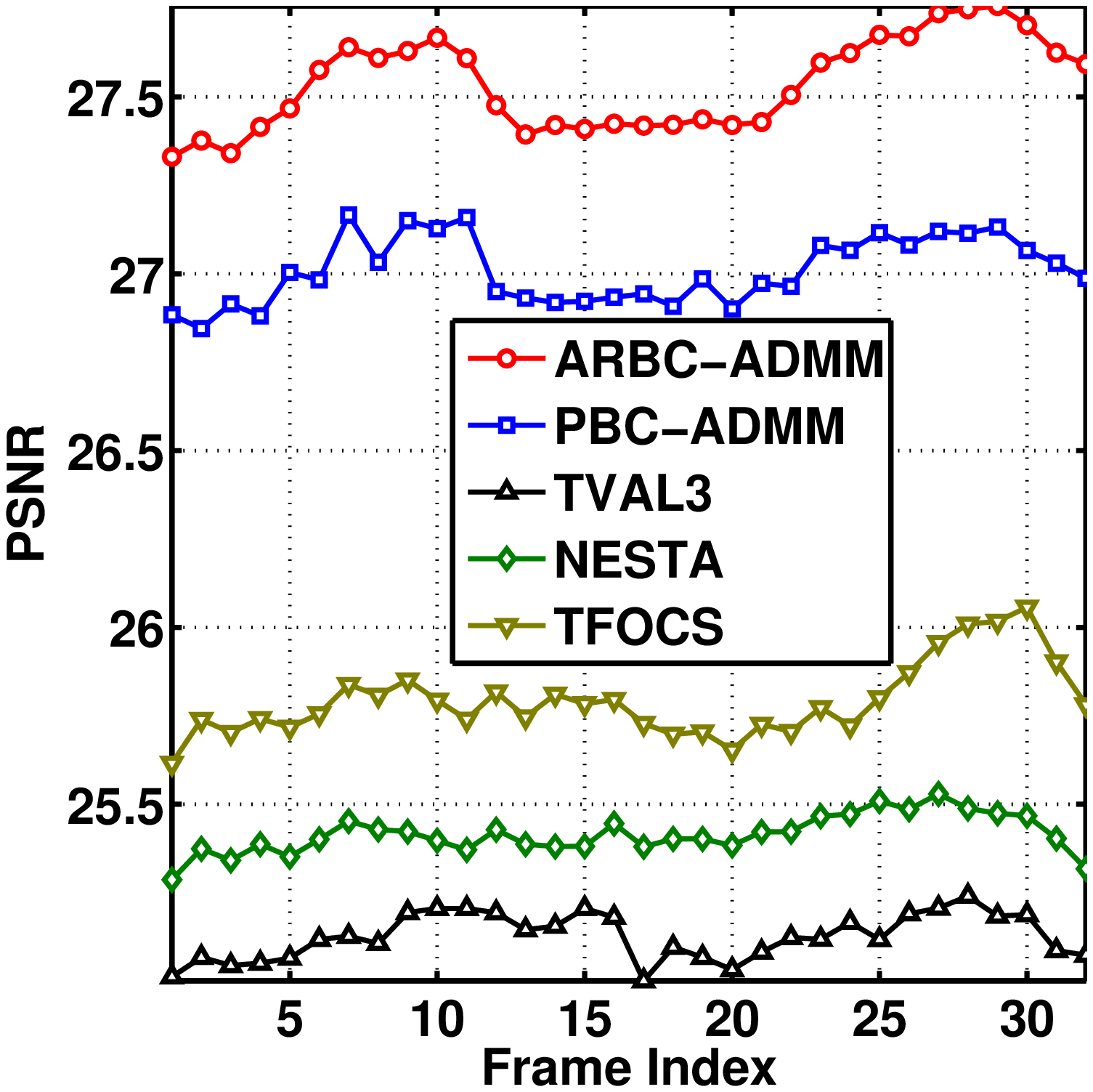}\hspace{-.05in}} \\
\hline 
{\hspace{-.05in}\begin{sideways} \hspace{.3in}$|{\bf\Omega}|/mnN=0.01$ \end{sideways}\hspace{-.05in}} &
{\hspace{-.05in}\includegraphics[width=0.2\textwidth]{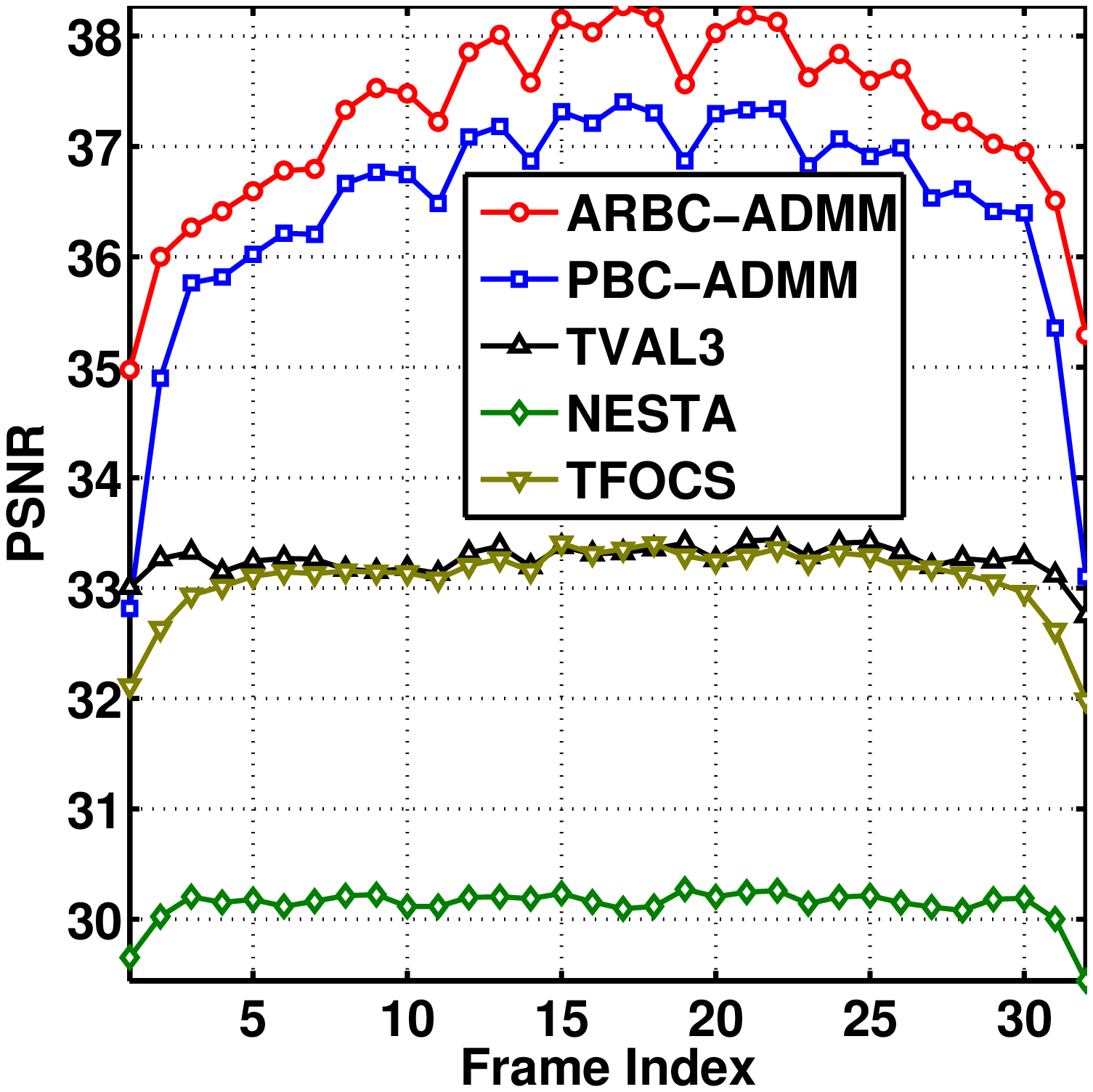}\hspace{-.05in}} &
{\hspace{-.05in}\includegraphics[width=0.2\textwidth]{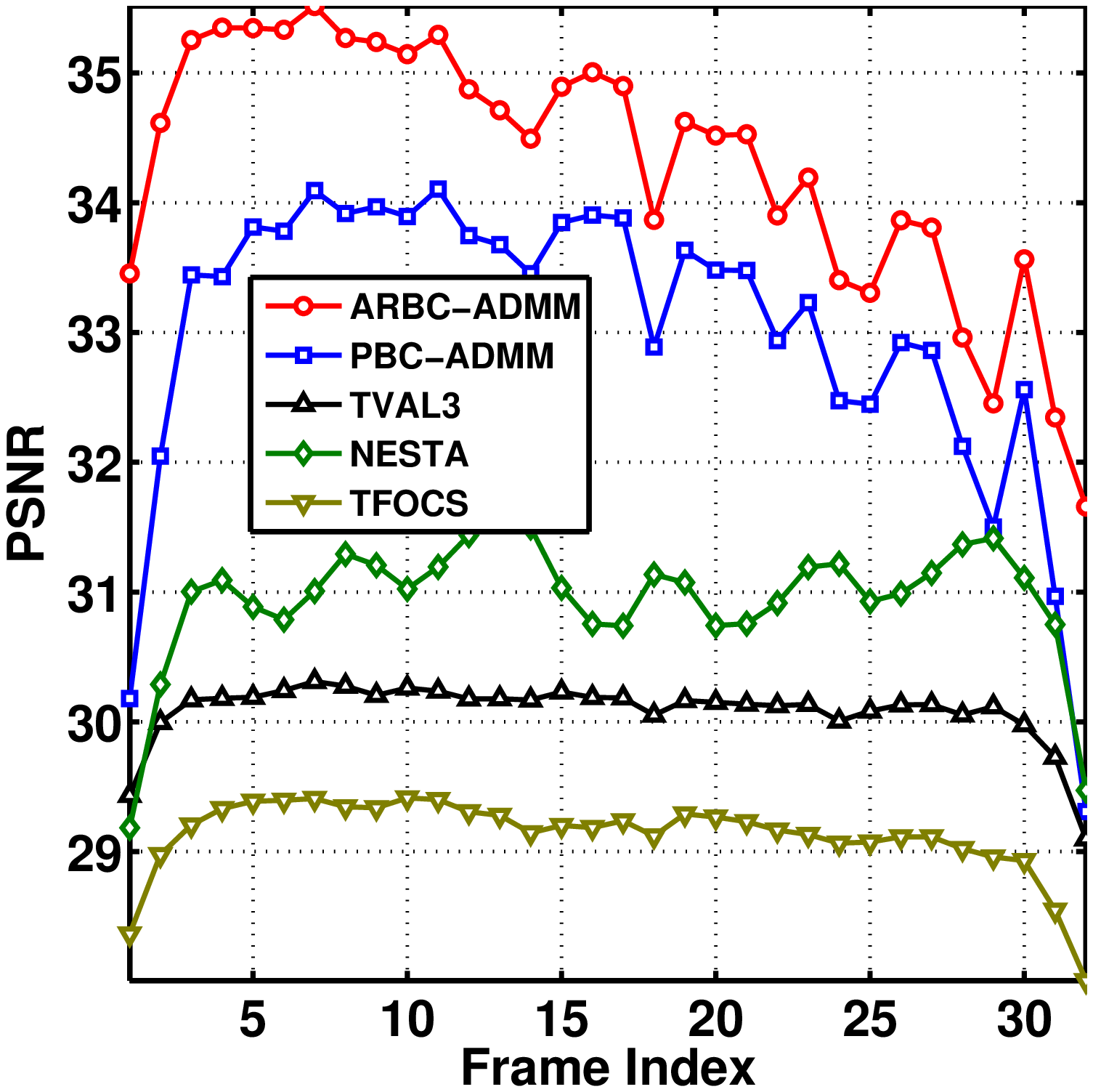}\hspace{-.05in}} &
{\hspace{-.05in}\includegraphics[width=0.2\textwidth]{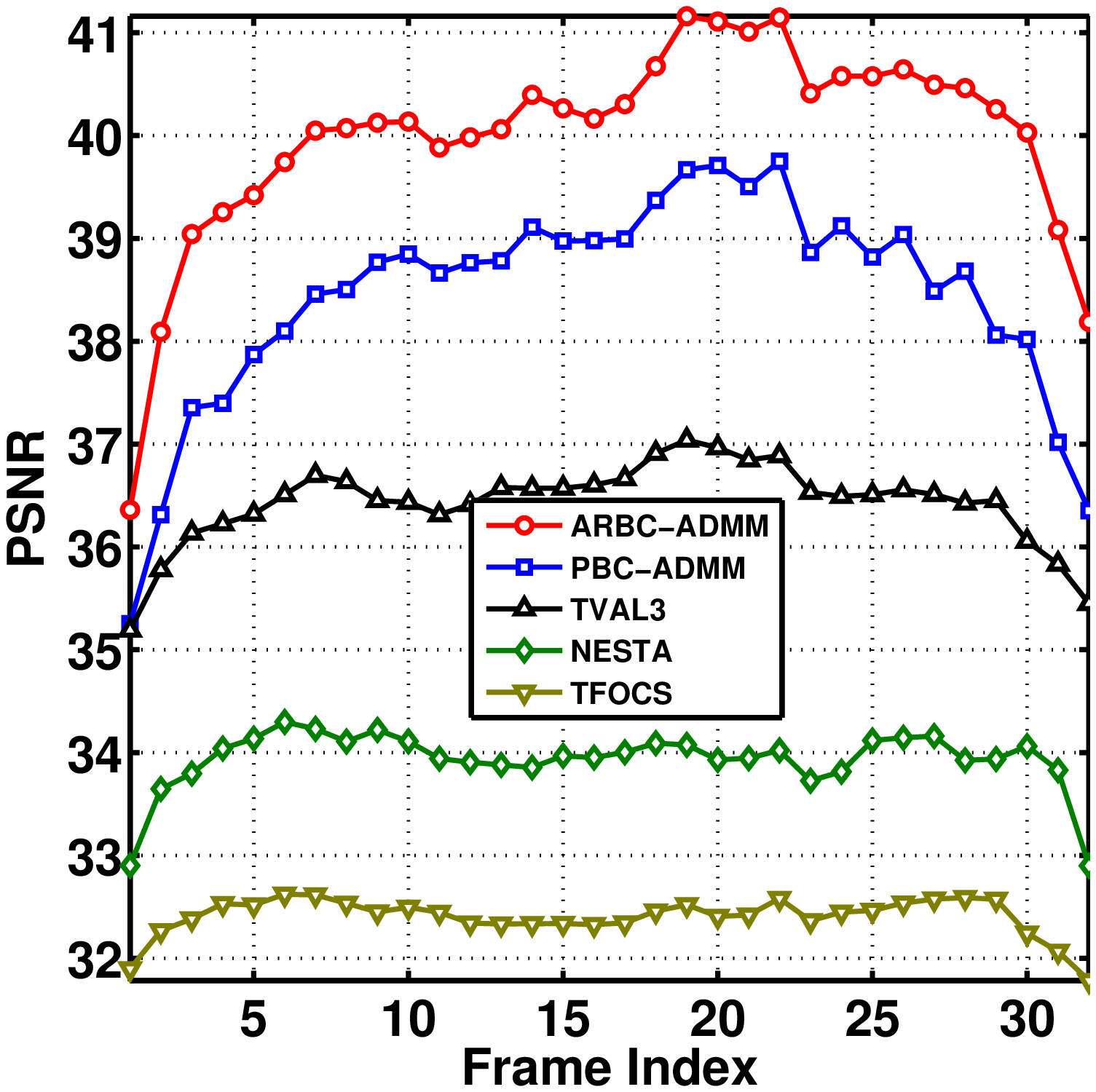}\hspace{-.05in}} &
{\hspace{-.05in}\includegraphics[width=0.2\textwidth]{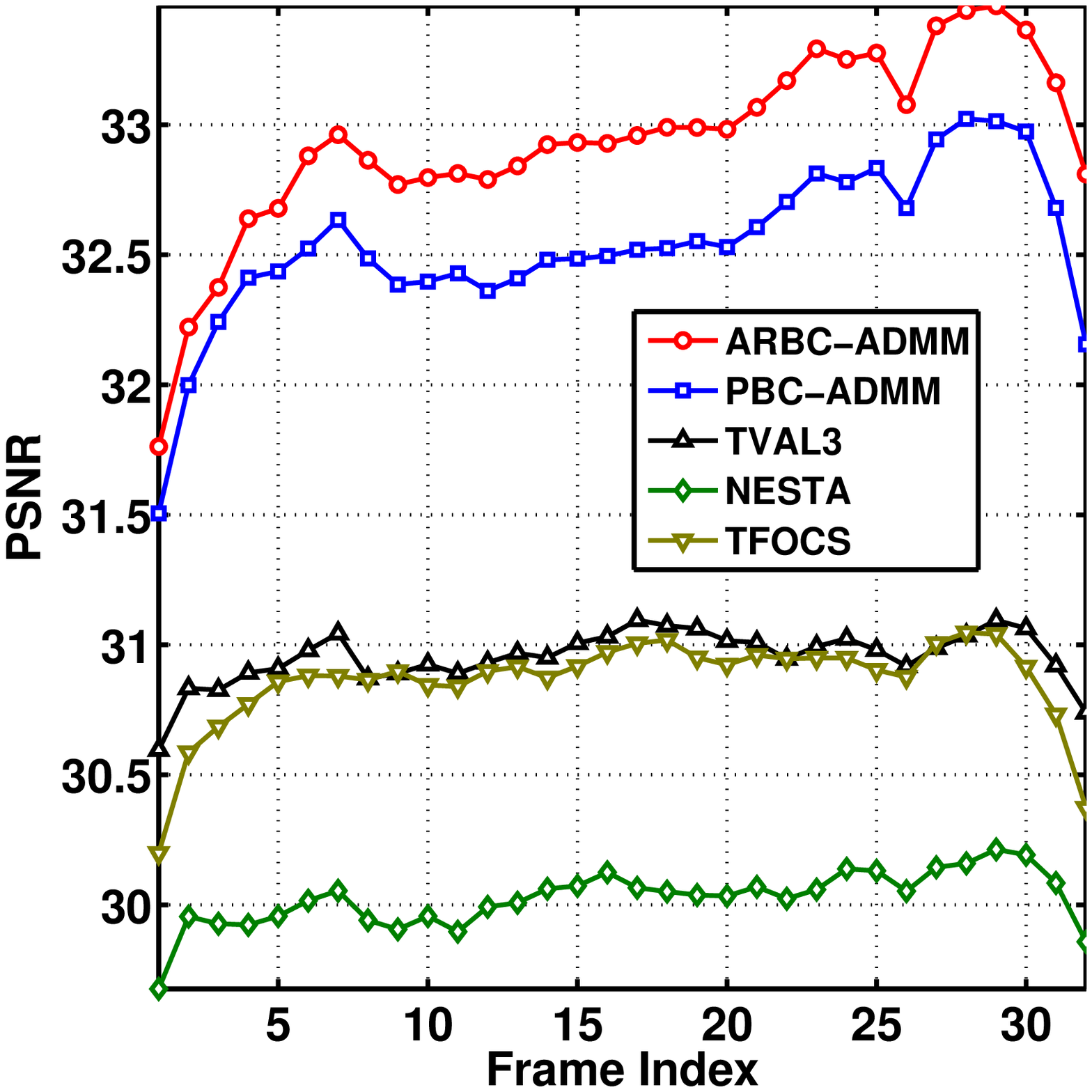}\hspace{-.05in}} \\
\hline
\end{tabular}
\end{center}
\end{table*}

\vspace{-.1in}
{
\subsection{Texture and Geometric Loss Analysis}
One of the biggest claim of this paper is addressing the deficiency issues in TV regularization: the texture and geometric loss. A texture in an image contain fine edge information which is pertinent to high frequency component in the signal. As we claimed in this paper, low accuracy differentiation such as $[-1,1]$ is unable to accurately transform such components and hence degrades the pertinent information. For the discussion we have selected two patches: P$1$ and P$2$ from the reconstructed frames in Table \ref{Table:FrameIndex17RecoveryHallMonitor} and magnify the resolution in Table \ref{Table:FrameIndex17SurfacePatch}. BC-ADMM is capable of reconstructing textural information on men' pant in image patch P$2$ from both under-sampling rates ${|{\bf\Omega}|}/{mnN}=\{0.03,0.1\}$, while TFOCS, NESTA, and TVAL3 almost lost every textural information in the image and degraded in quality with spiky and blurry recoveries. Similar loss is happened in patch P$1$ where the slanted line in the wall is disappeared in recovery from the competing methods.}

\begin{table*}
\renewcommand{\arraystretch}{1.3}
\caption{Reconstructed image patches via five different CVS approaches from two sampling rates.}
\label{Table:FrameIndex17SurfacePatch}\vspace{-.1in}
\begin{center}
\begin{tabular}{|c|c|c|c|c|c|c|c|}
\cline{3-8}
\multicolumn{2}{c|}{} & \text{Reference Patch} & \text{AR-BC-ADMM} & \text{P-BC-ADMM} & \text{TFOCS} & \text{NESTA} & \text{TVAL3} \\
\hline
{\hspace{-.05in}\begin{sideways} \hspace{.3in} $\frac{|{\bf\Omega}|}{mnN}=0.03$ \end{sideways}\hspace{-.05in}} &
{\hspace{-.05in}\begin{sideways} \hspace{.4in} Patch No.$1$ \end{sideways}\hspace{-.05in}} &
{\hspace{-.05in}\includegraphics[width=0.125\textwidth]{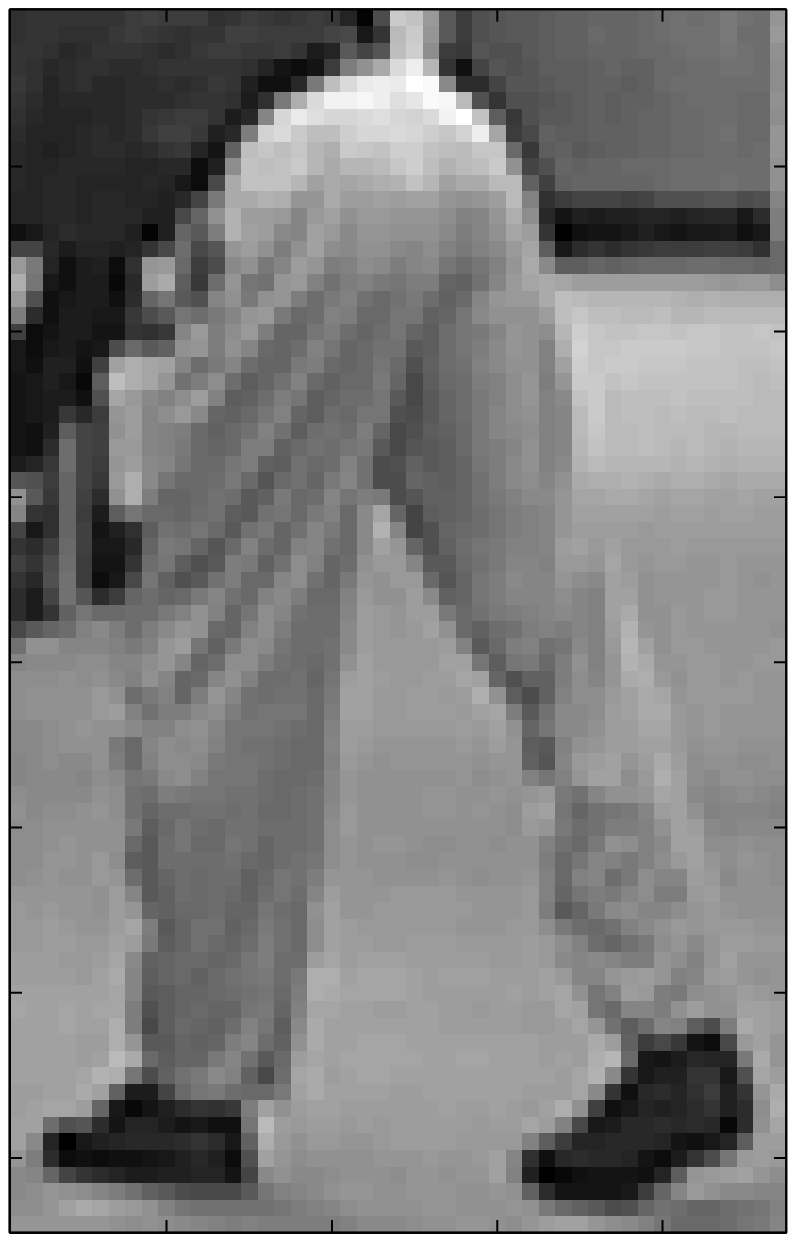}\hspace{-.05in}} &
{\hspace{-.05in}\includegraphics[width=0.125\textwidth]{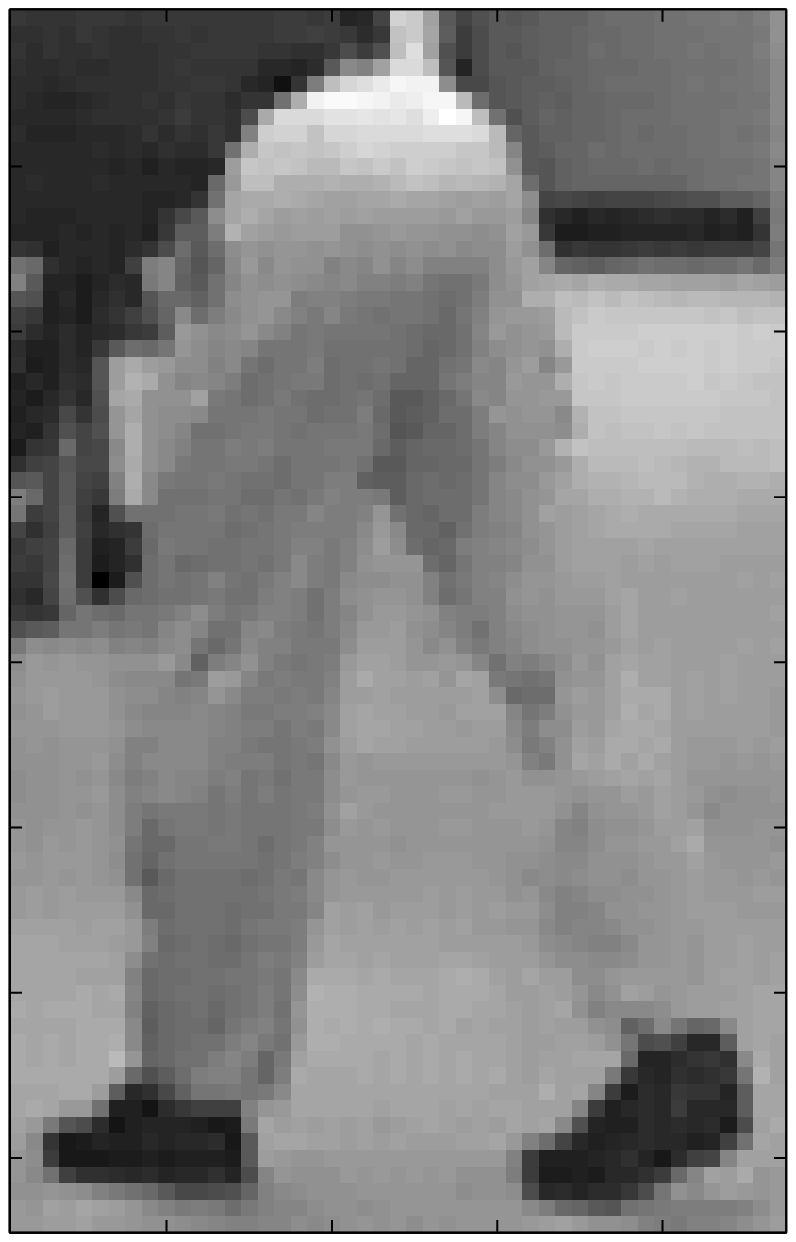}\hspace{-.05in}} &
{\hspace{-.05in}\includegraphics[width=0.125\textwidth]{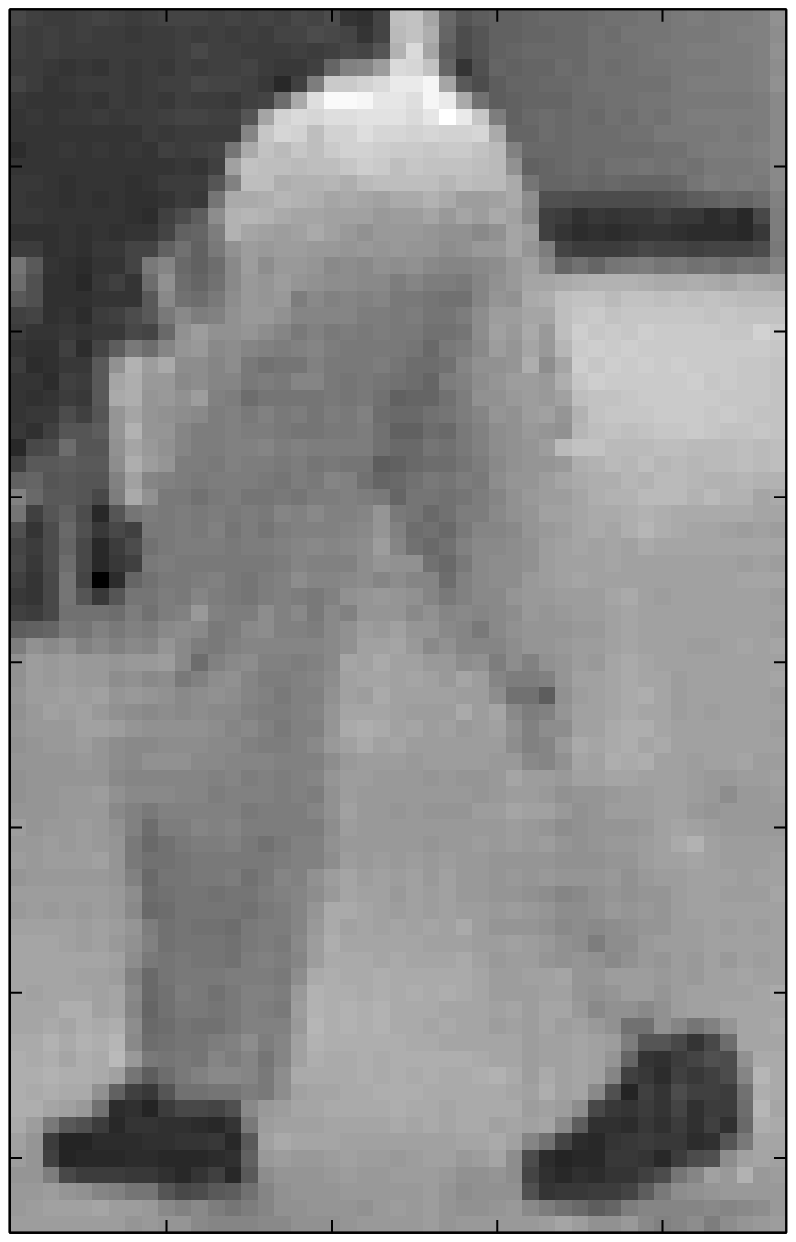}\hspace{-.05in}} &
{\hspace{-.05in}\includegraphics[width=0.125\textwidth]{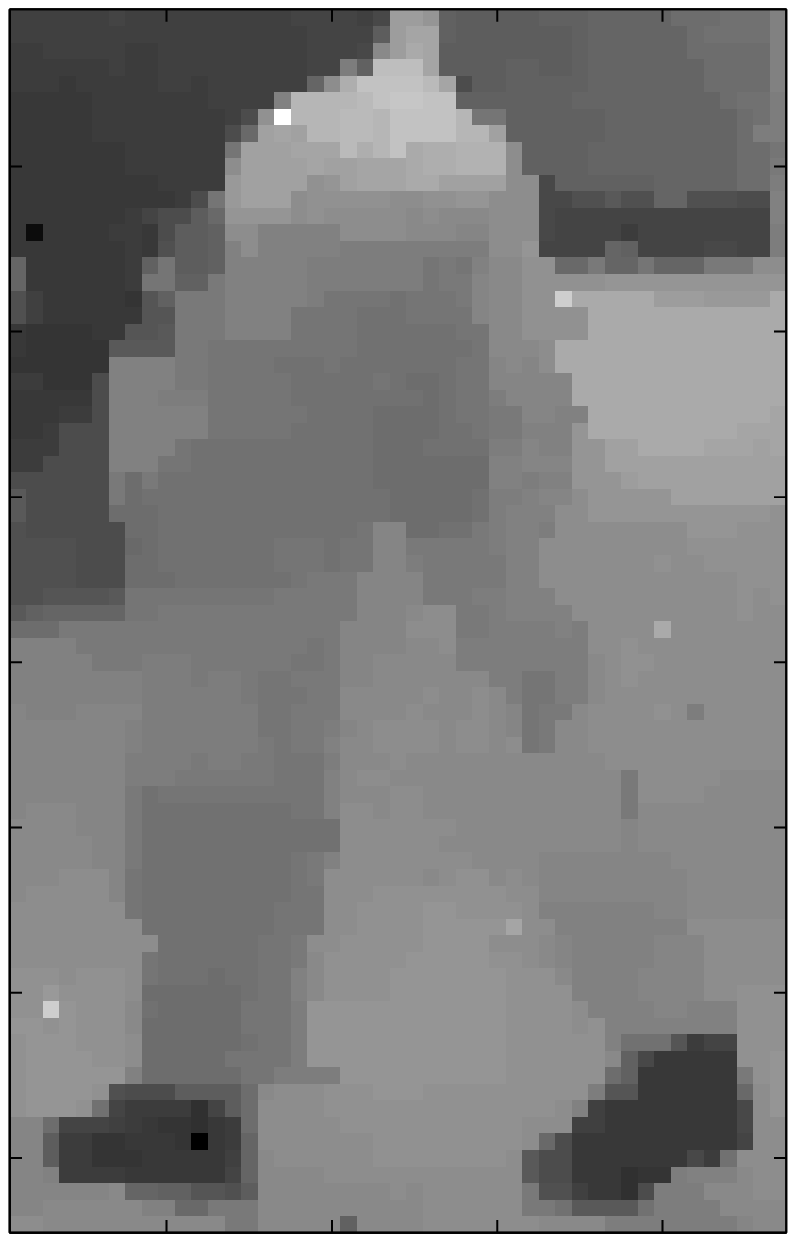}\hspace{-.05in}} &
{\hspace{-.05in}\includegraphics[width=0.125\textwidth]{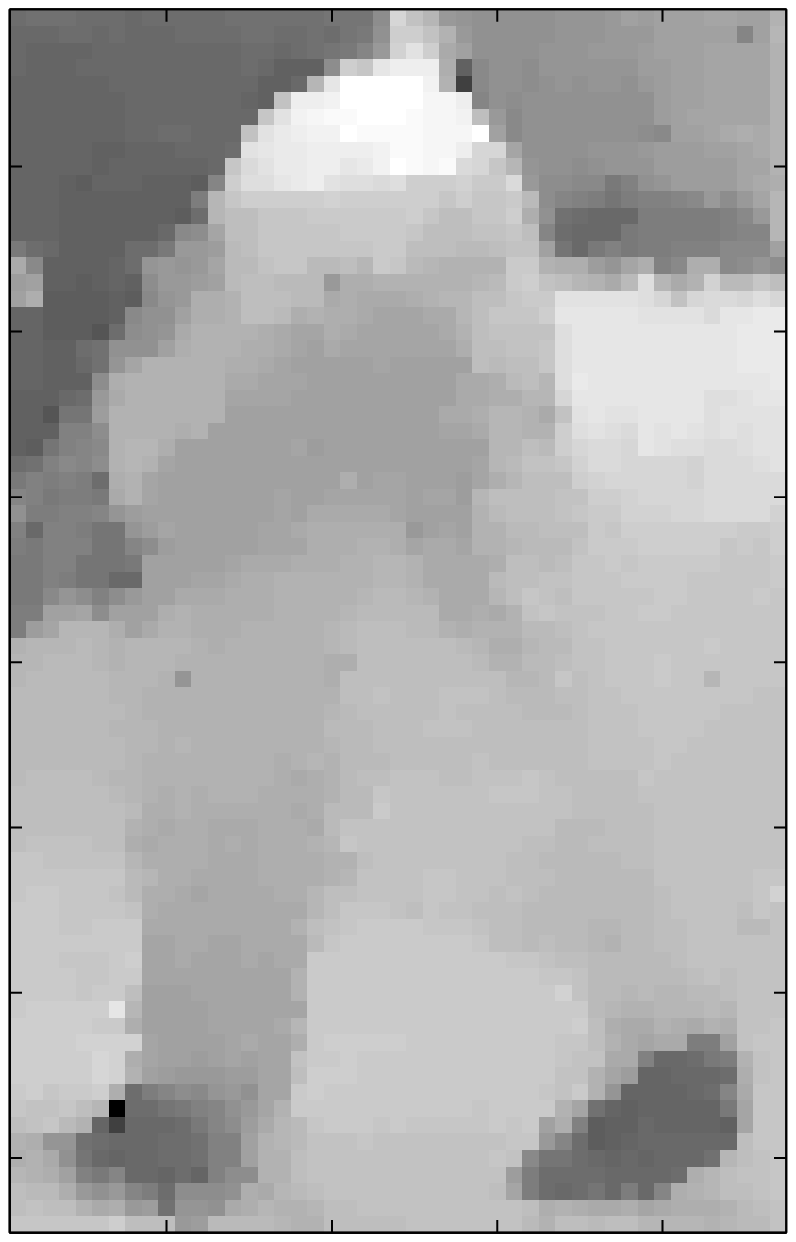}\hspace{-.05in}} &
{\hspace{-.05in}\includegraphics[width=0.125\textwidth]{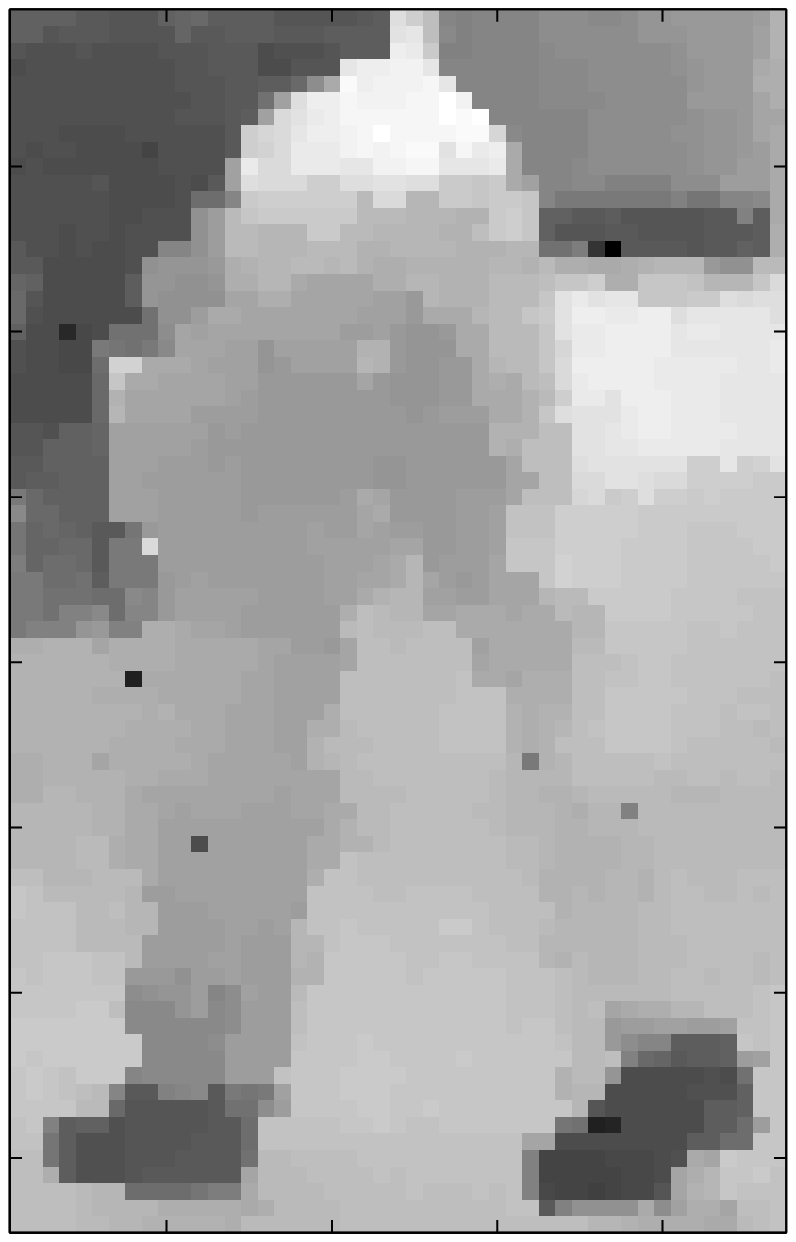}\hspace{-.05in}} \\
\hline 
{\hspace{-.05in}\begin{sideways} \hspace{-.05in} $\frac{|{\bf\Omega}|}{mnN}=0.03$ \end{sideways}\hspace{-.05in}} &
{\hspace{-.05in}\begin{sideways} \hspace{.075in}Patch No.$2$ \end{sideways}\hspace{-.05in}} &
{\hspace{-.05in}\includegraphics[width=0.125\textwidth]{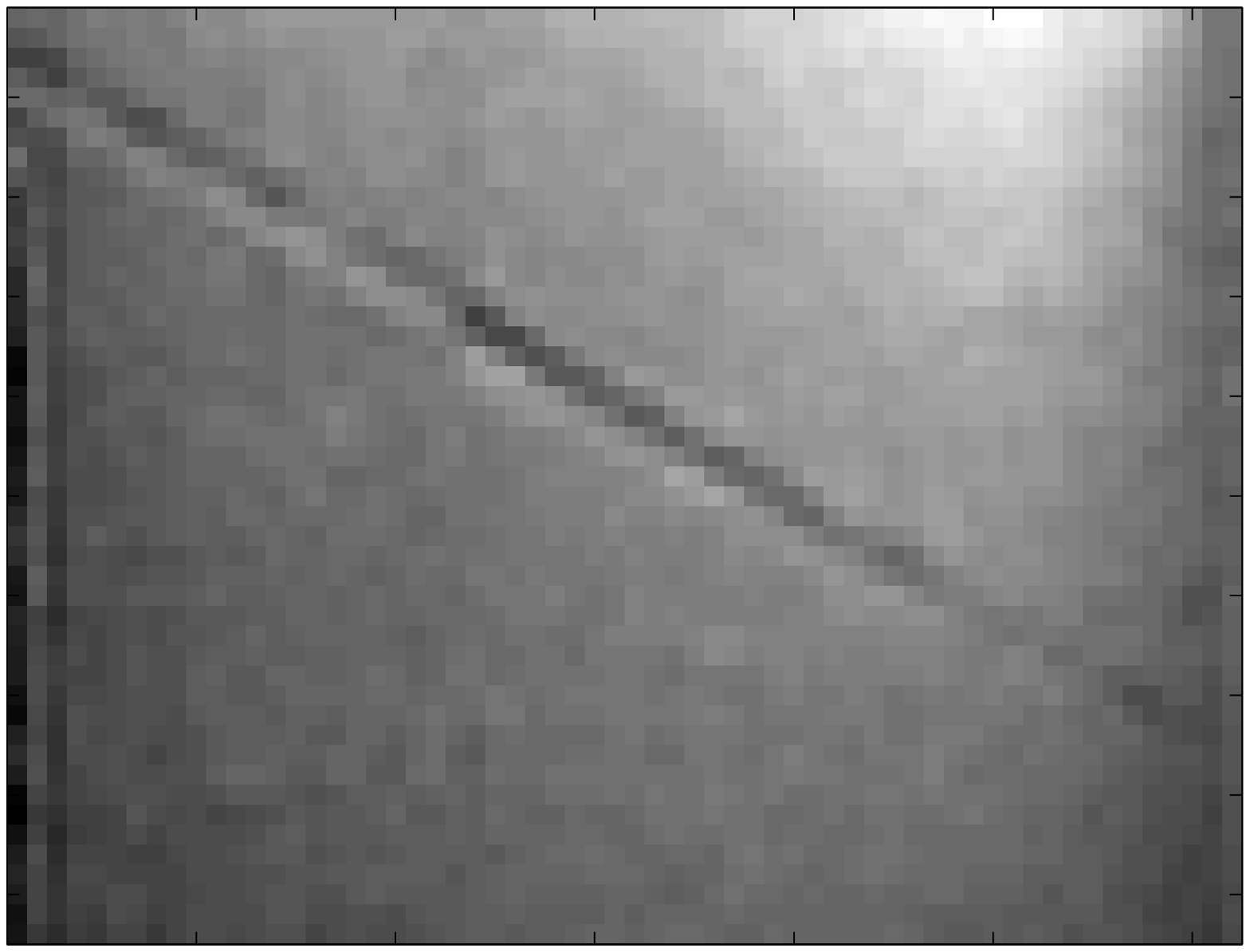}\hspace{-.05in}} &
{\hspace{-.05in}\includegraphics[width=0.125\textwidth]{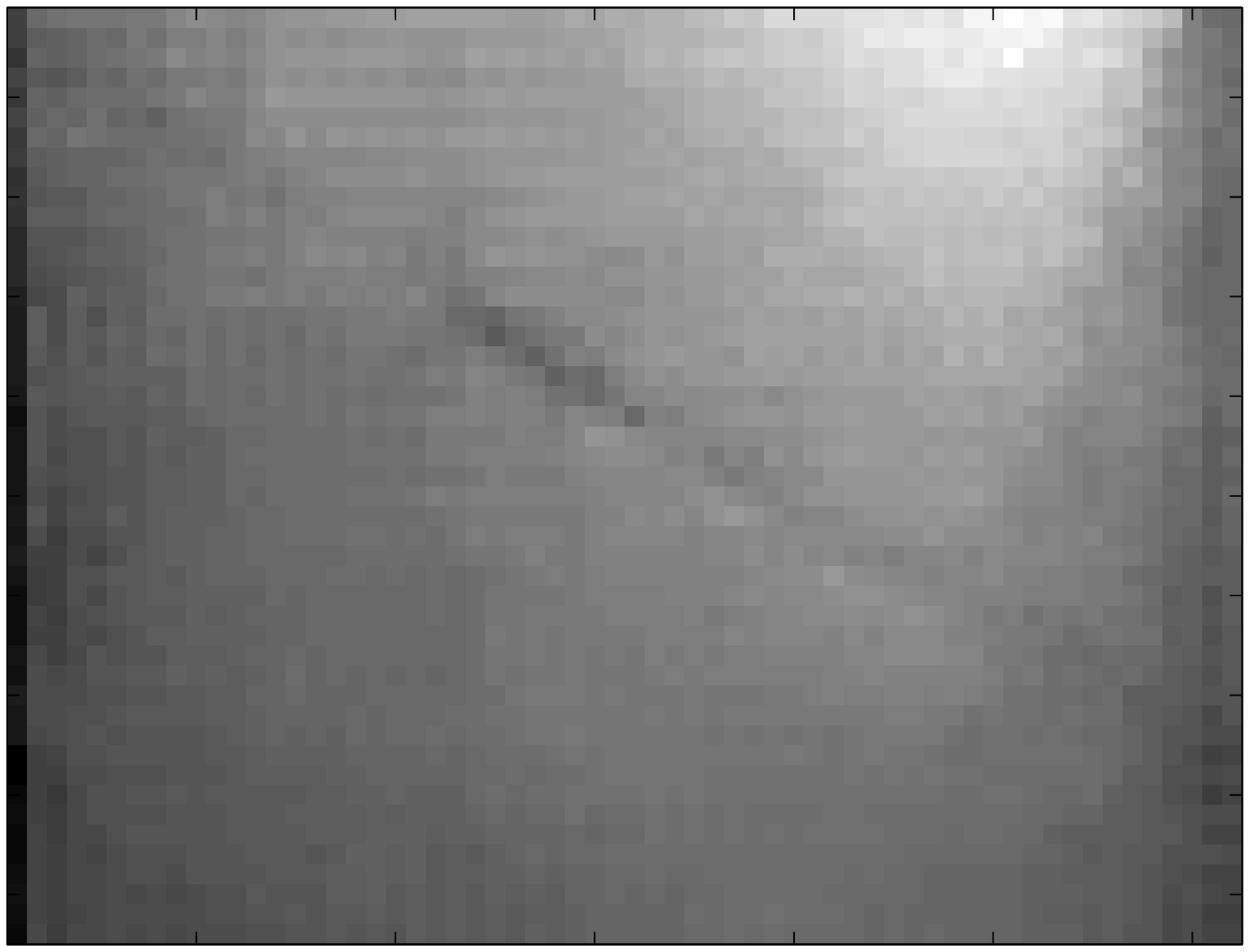}\hspace{-.05in}} &
{\hspace{-.05in}\includegraphics[width=0.125\textwidth]{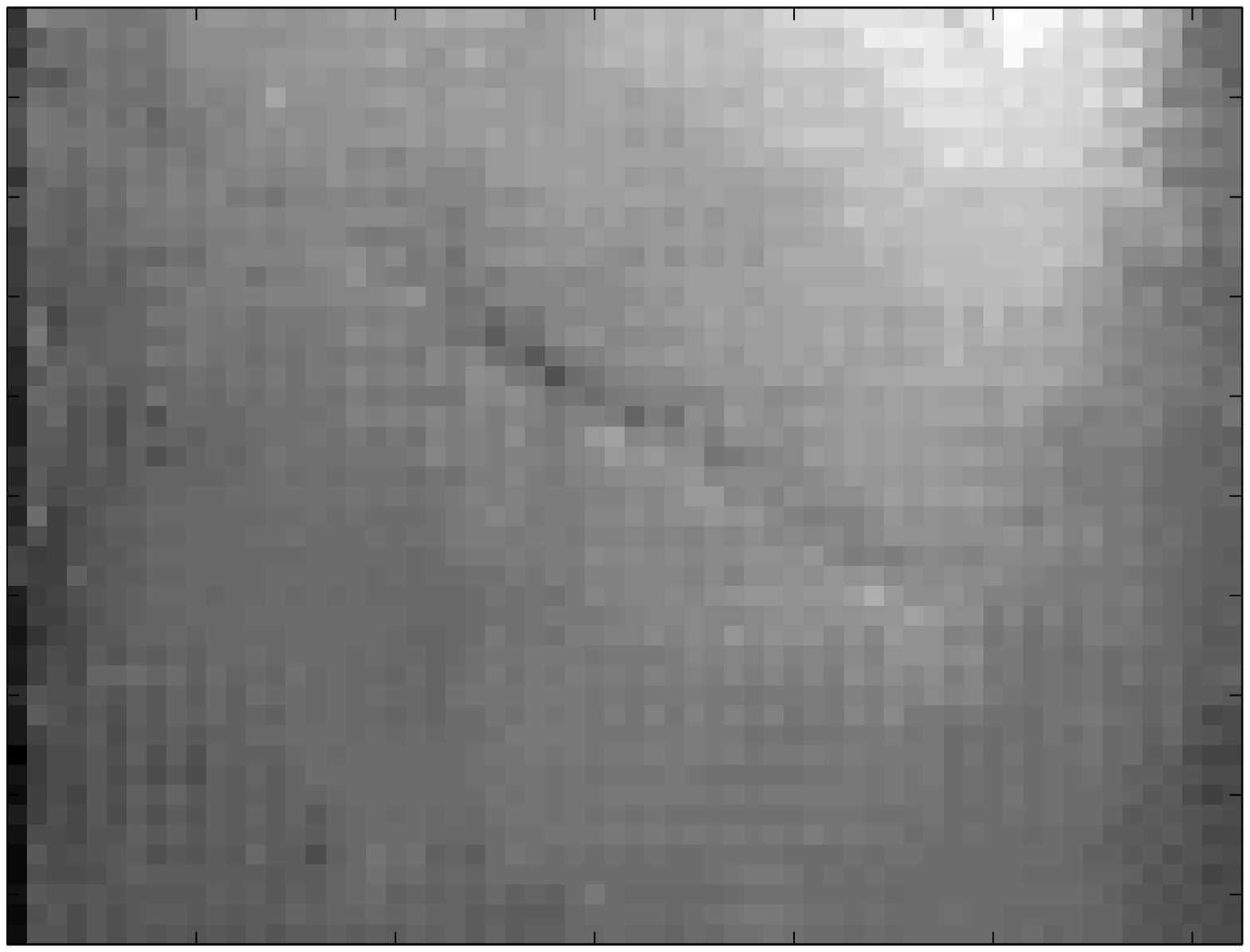}\hspace{-.05in}} &
{\hspace{-.05in}\includegraphics[width=0.125\textwidth]{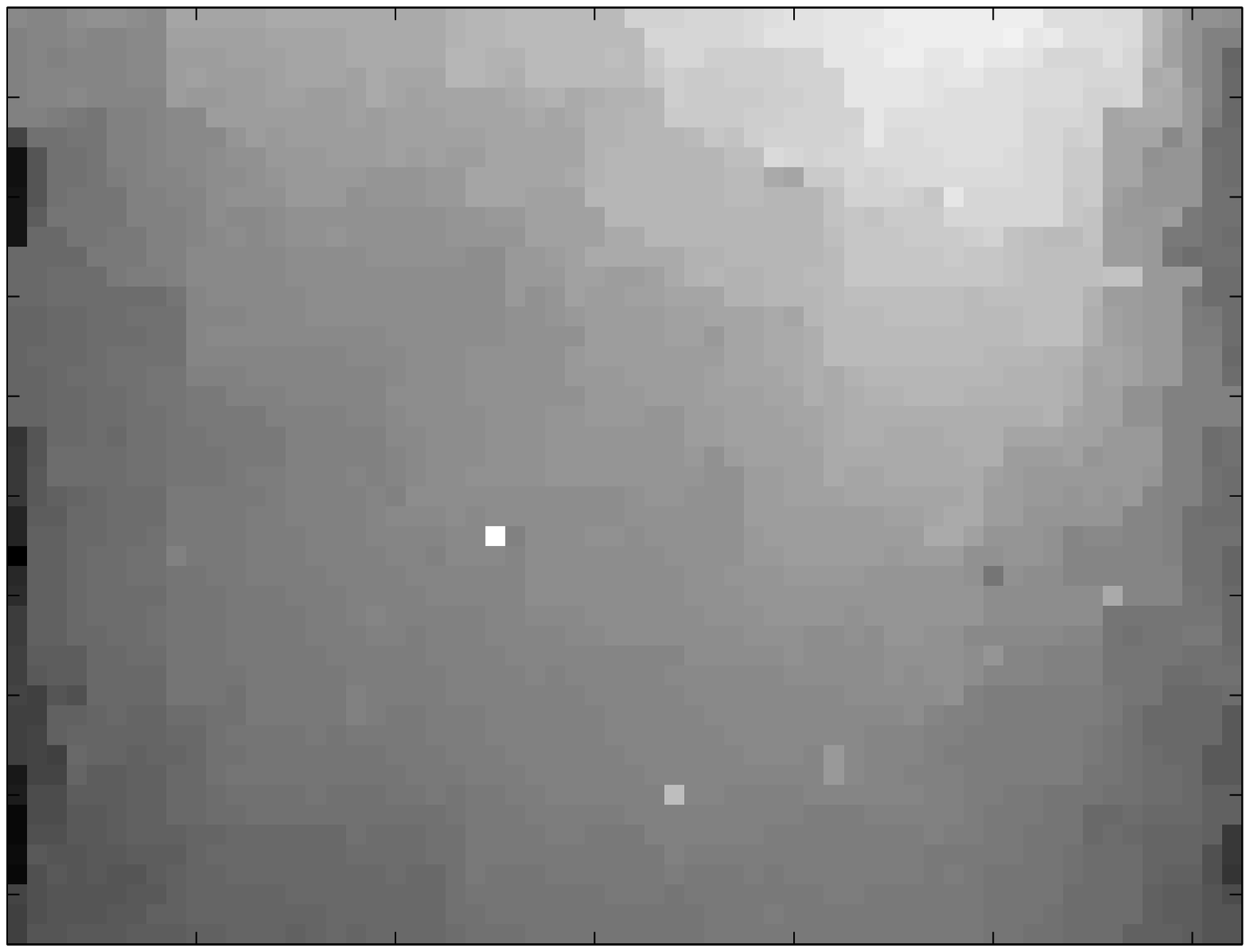}\hspace{-.05in}} &
{\hspace{-.05in}\includegraphics[width=0.125\textwidth]{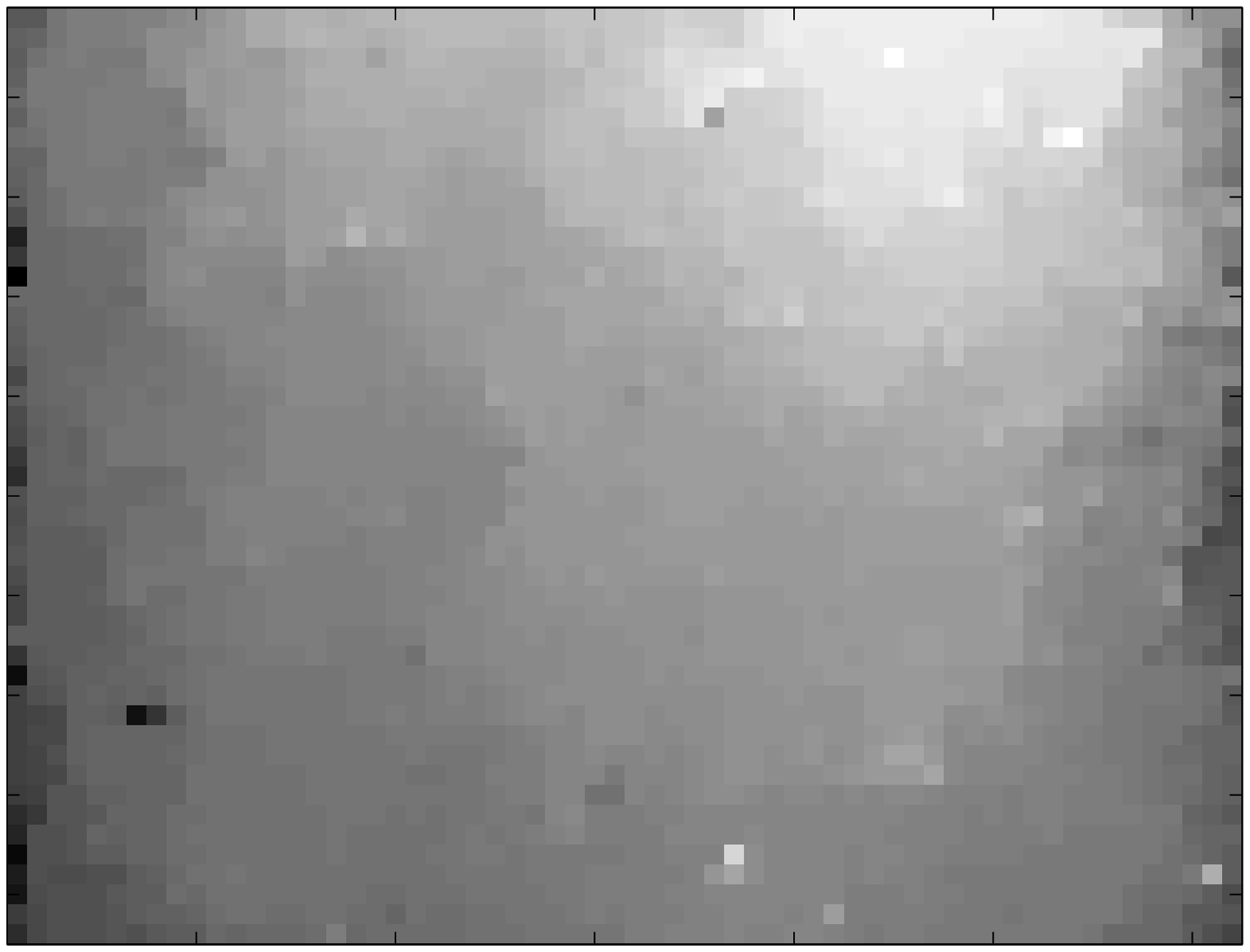}\hspace{-.05in}} &
{\hspace{-.05in}\includegraphics[width=0.125\textwidth]{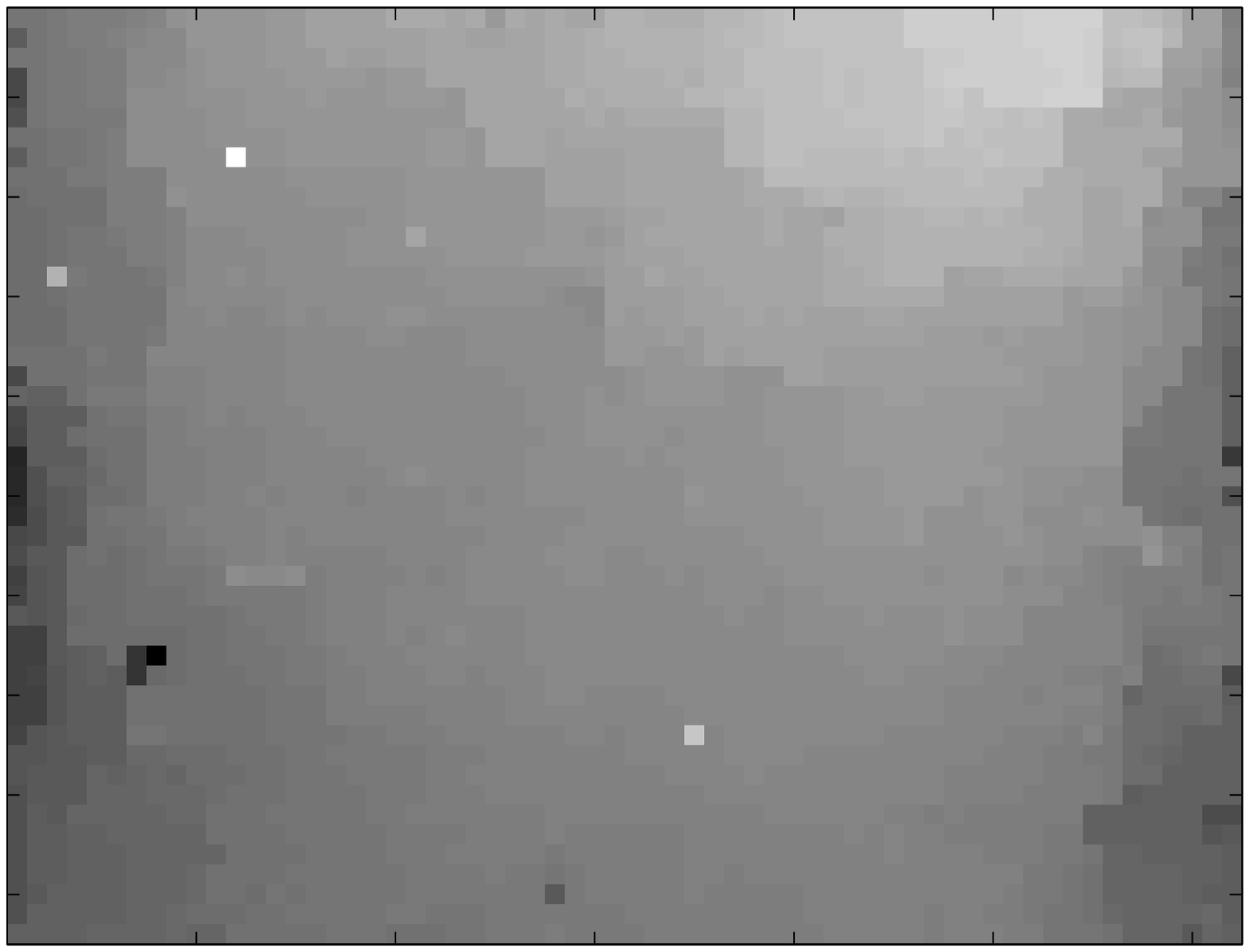}\hspace{-.05in}} \\
\hline 
{\hspace{-.05in}\begin{sideways} \hspace{.3in} $\frac{|{\bf\Omega}|}{mnN}=0.1$ \end{sideways}\hspace{-.05in}} &
{\hspace{-.05in}\begin{sideways} \hspace{.4in} Patch No.$1$ \end{sideways}\hspace{-.05in}} &
{\hspace{-.05in}\includegraphics[width=0.125\textwidth]{F17_Original_Patch1.eps}\hspace{-.05in}} &
{\hspace{-.05in}\includegraphics[width=0.125\textwidth]{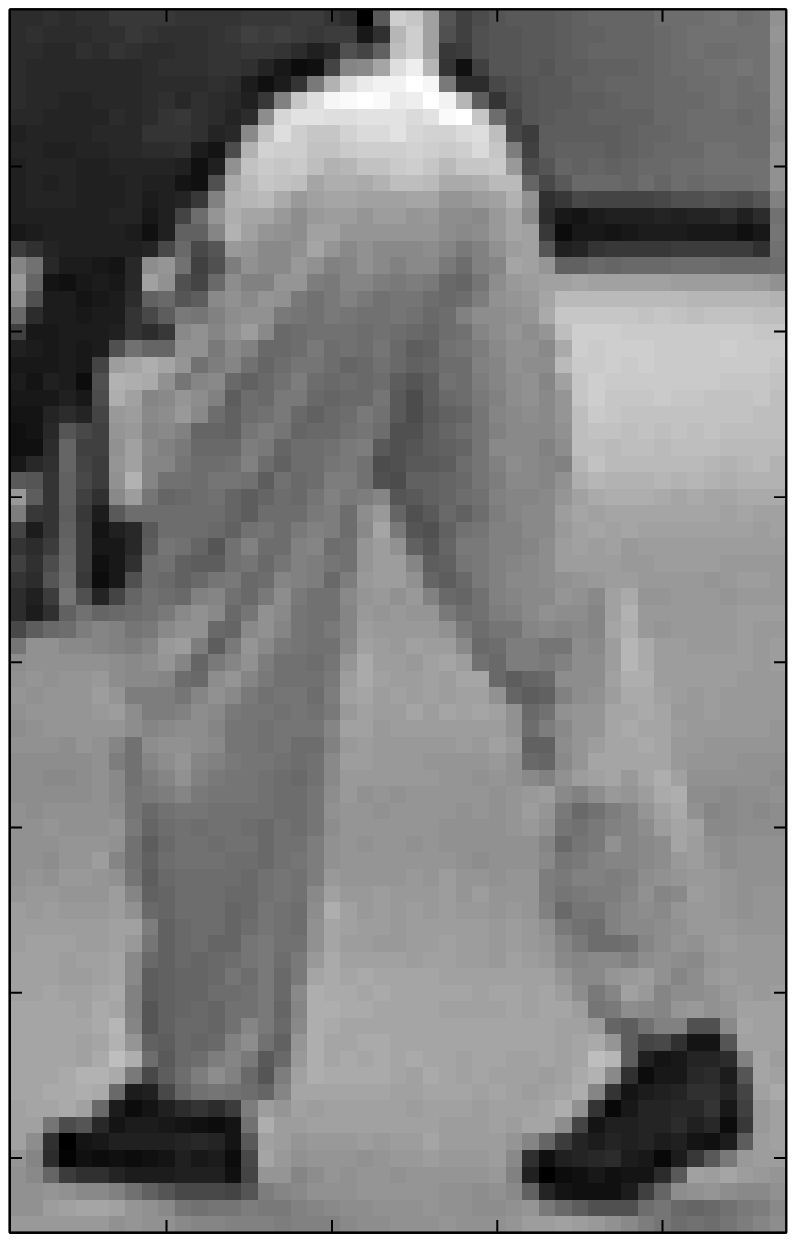}\hspace{-.05in}} &
{\hspace{-.05in}\includegraphics[width=0.125\textwidth]{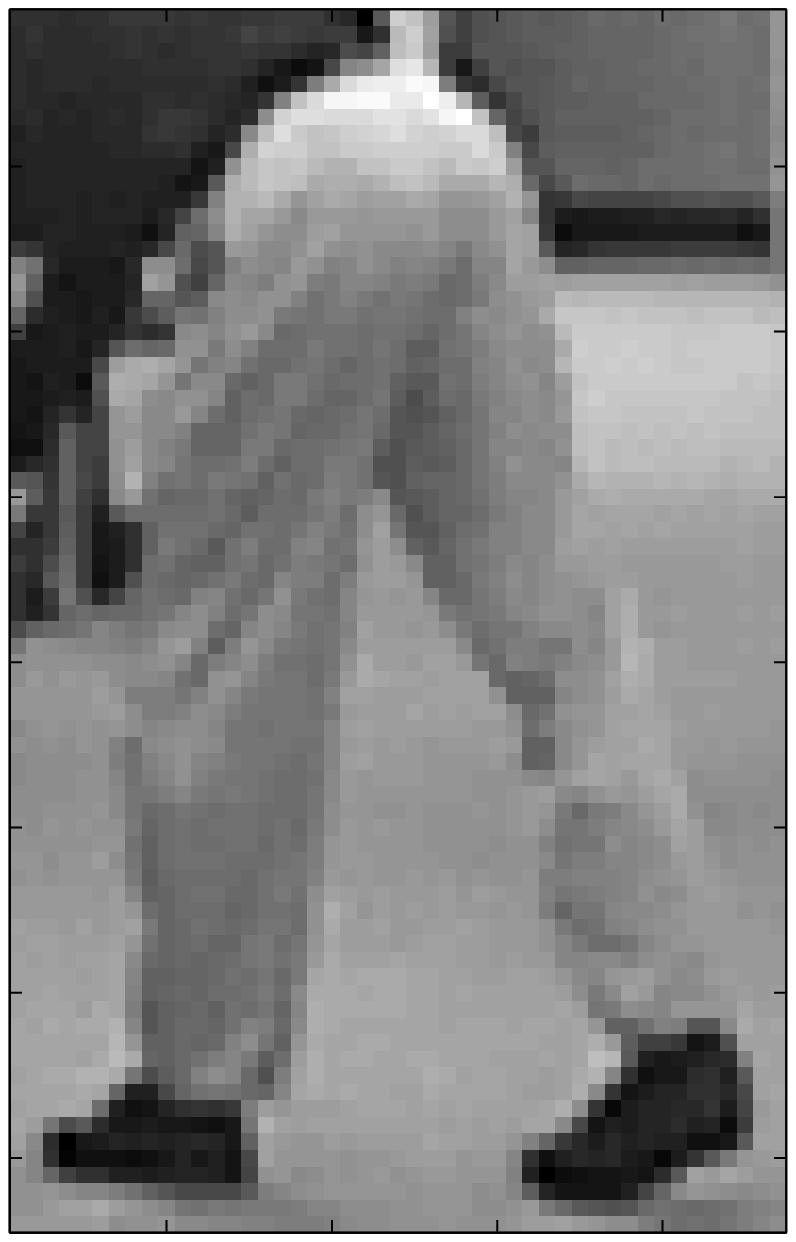}\hspace{-.05in}} &
{\hspace{-.05in}\includegraphics[width=0.125\textwidth]{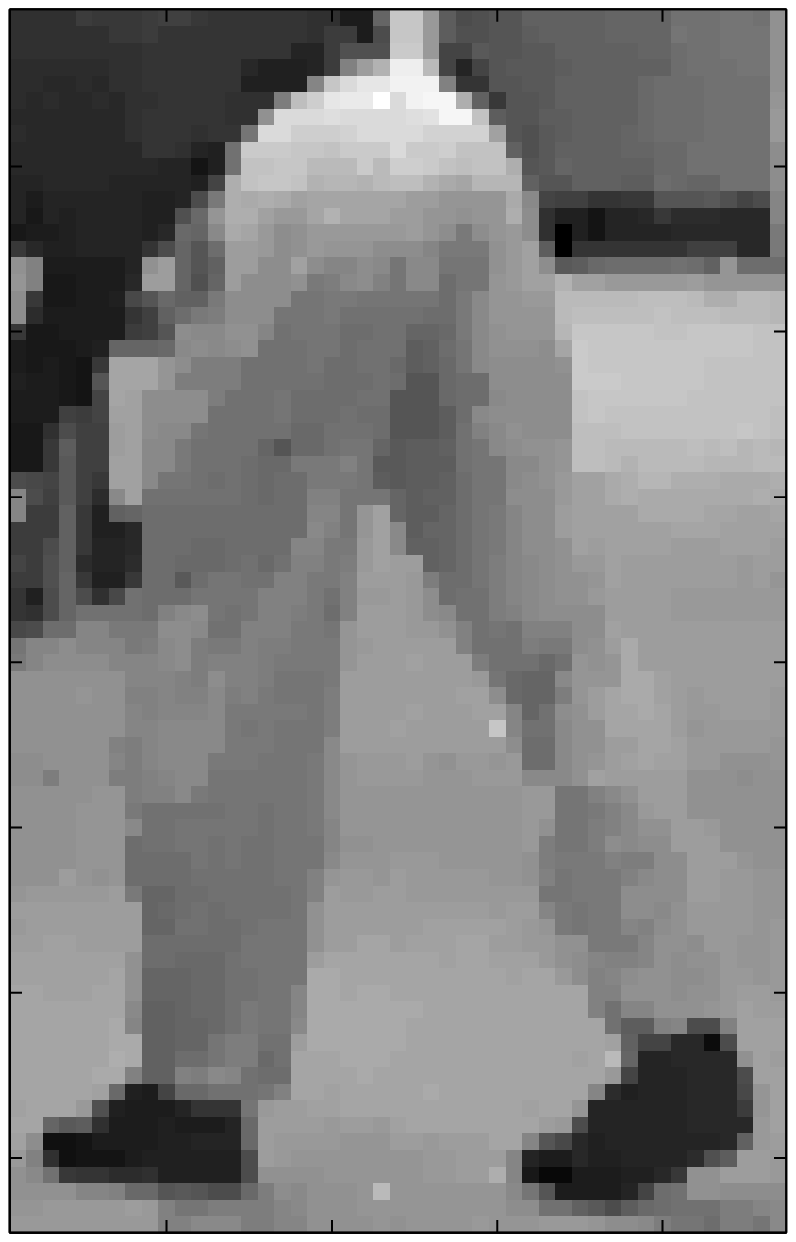}\hspace{-.05in}} &
{\hspace{-.05in}\includegraphics[width=0.125\textwidth]{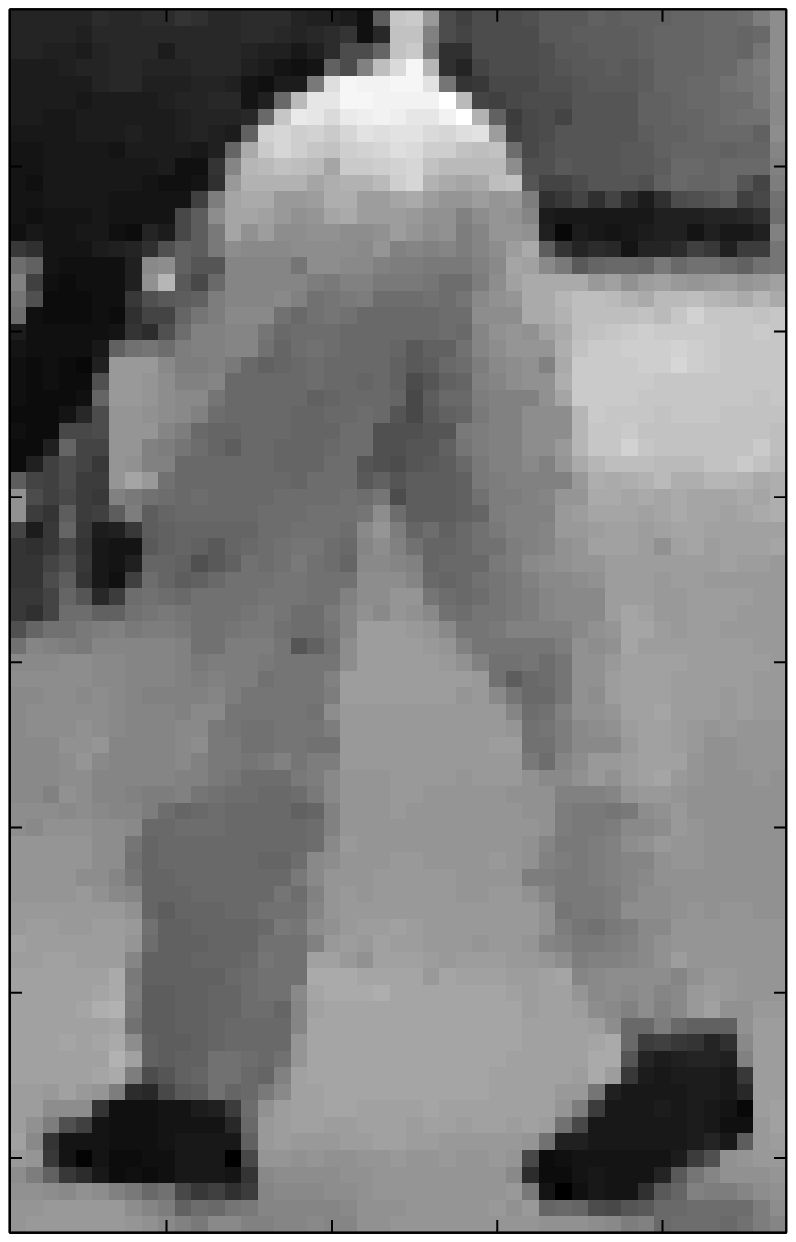}\hspace{-.05in}} &
{\hspace{-.05in}\includegraphics[width=0.125\textwidth]{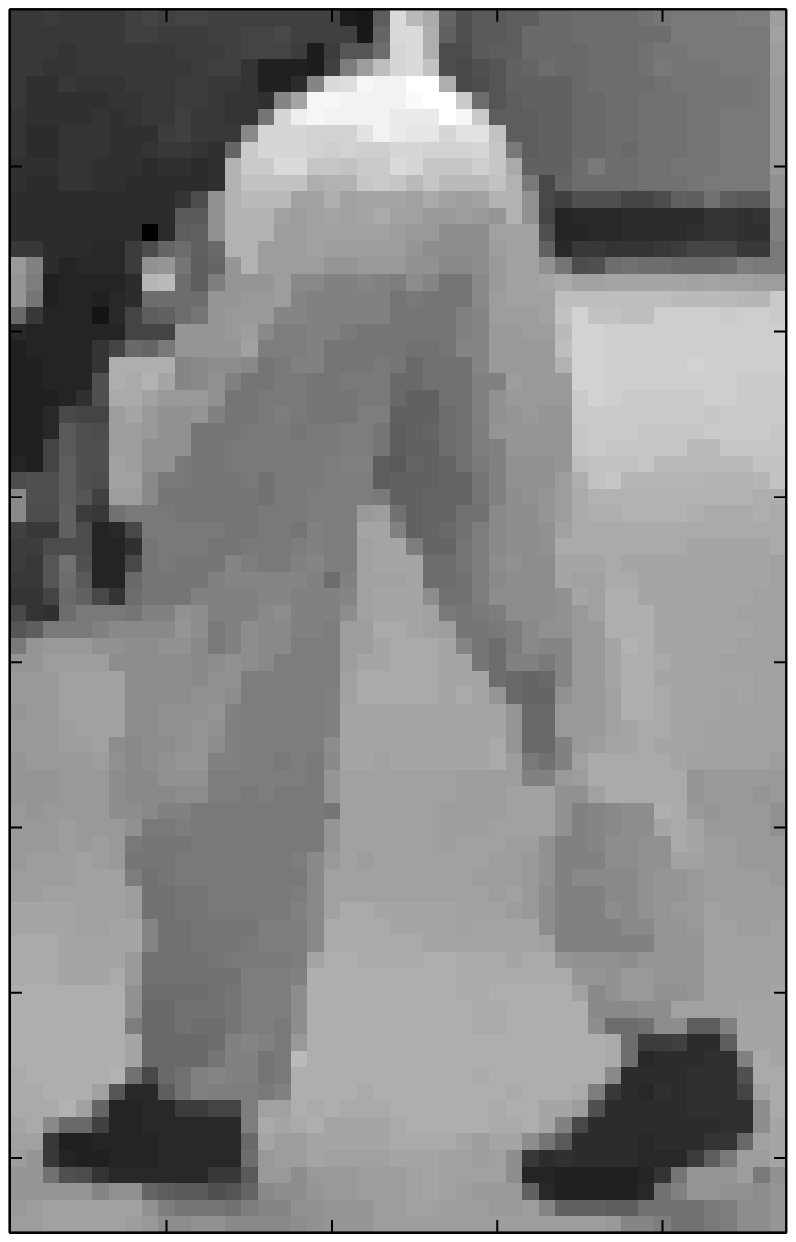}\hspace{-.05in}} \\
\hline 
{\hspace{-.05in}\begin{sideways} \hspace{.0in} $\frac{|{\bf\Omega}|}{mnN}=0.1$ \end{sideways}\hspace{-.05in}} &
{\hspace{-.05in}\begin{sideways} \hspace{.075in}Patch No.$2$ \end{sideways}\hspace{-.05in}} &
{\hspace{-.05in}\includegraphics[width=0.125\textwidth]{F17_Original_Patch2.eps}\hspace{-.05in}} &
{\hspace{-.05in}\includegraphics[width=0.125\textwidth]{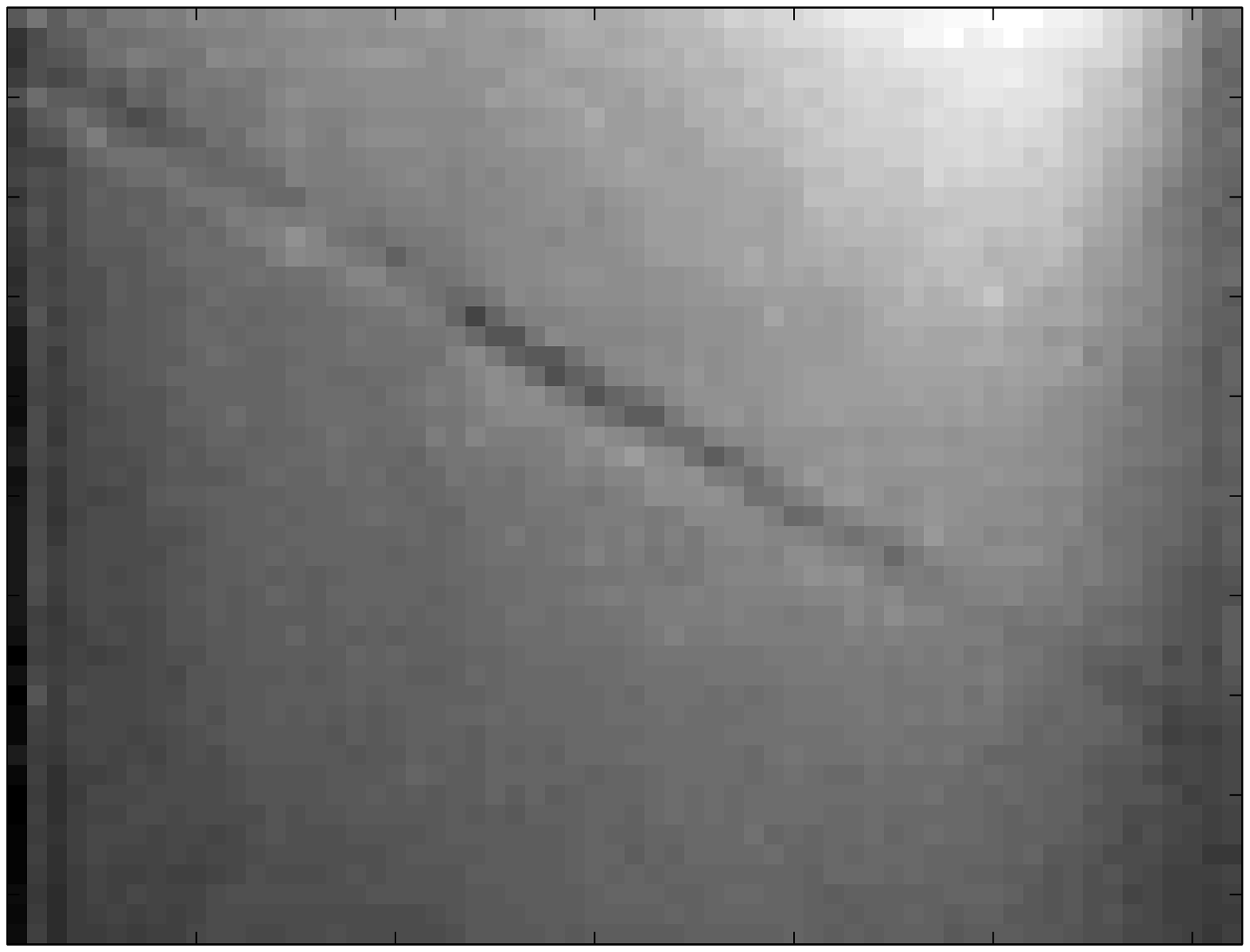}\hspace{-.05in}} &
{\hspace{-.05in}\includegraphics[width=0.125\textwidth]{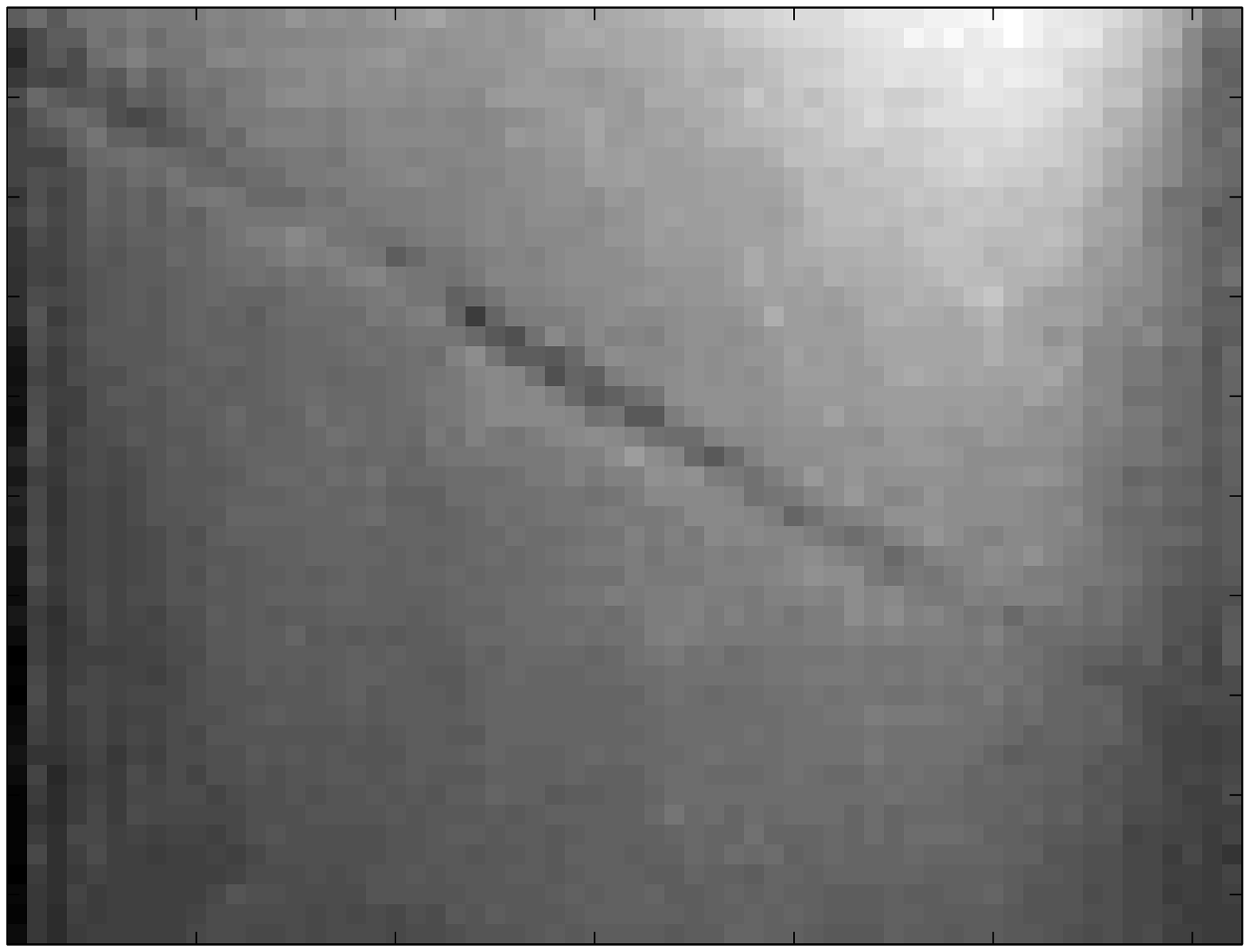}\hspace{-.05in}} &
{\hspace{-.05in}\includegraphics[width=0.125\textwidth]{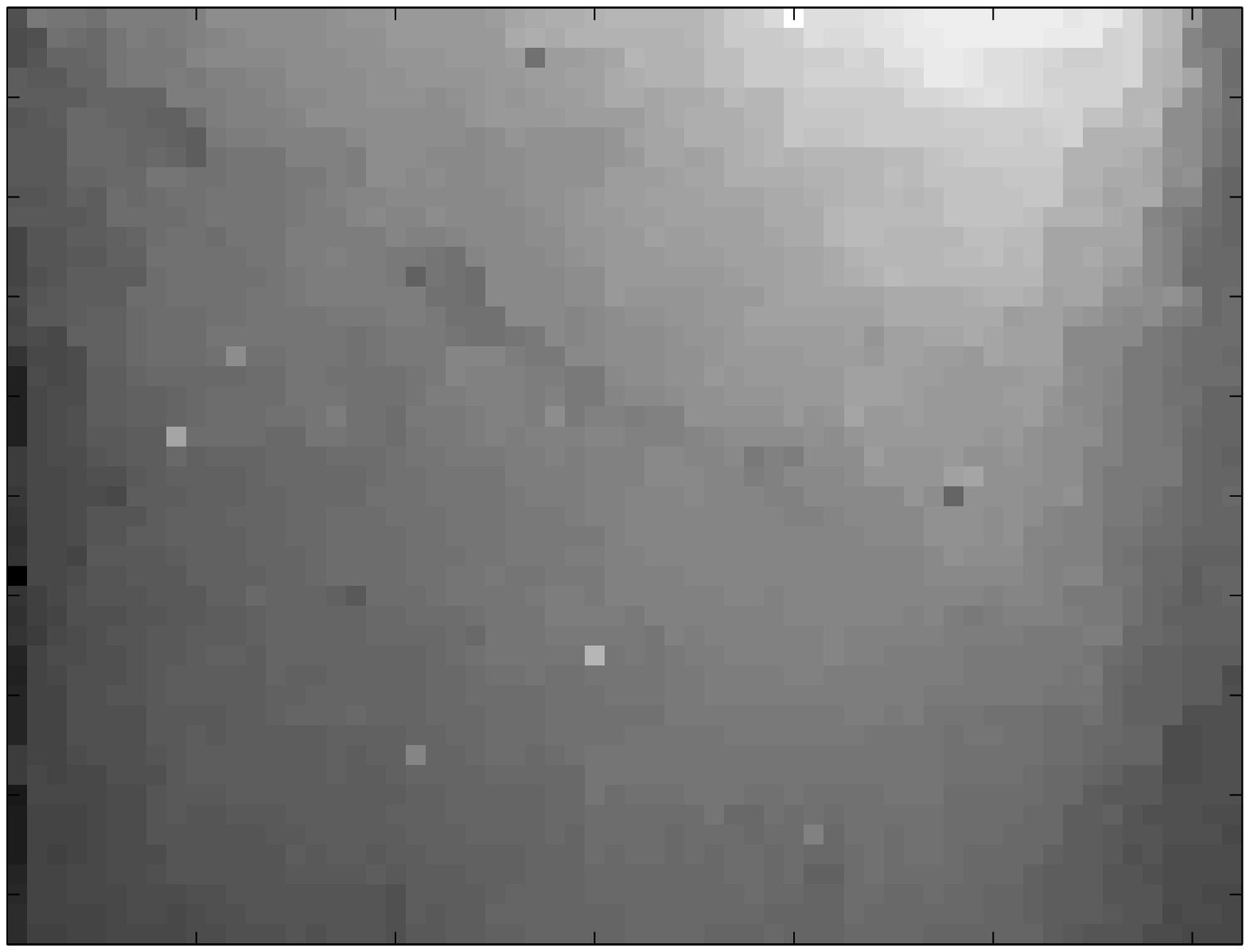}\hspace{-.05in}} &
{\hspace{-.05in}\includegraphics[width=0.125\textwidth]{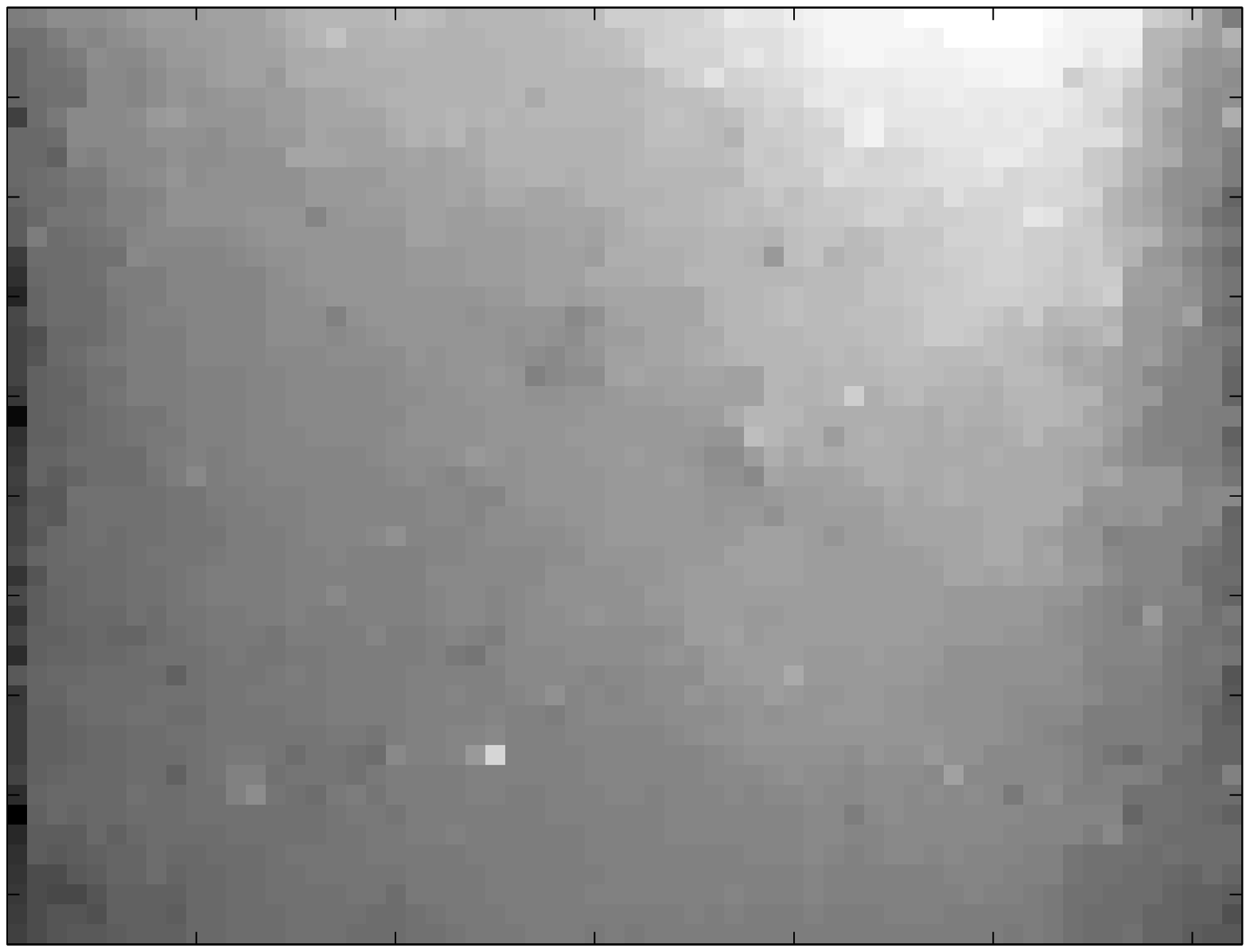}\hspace{-.05in}} &
{\hspace{-.05in}\includegraphics[width=0.125\textwidth]{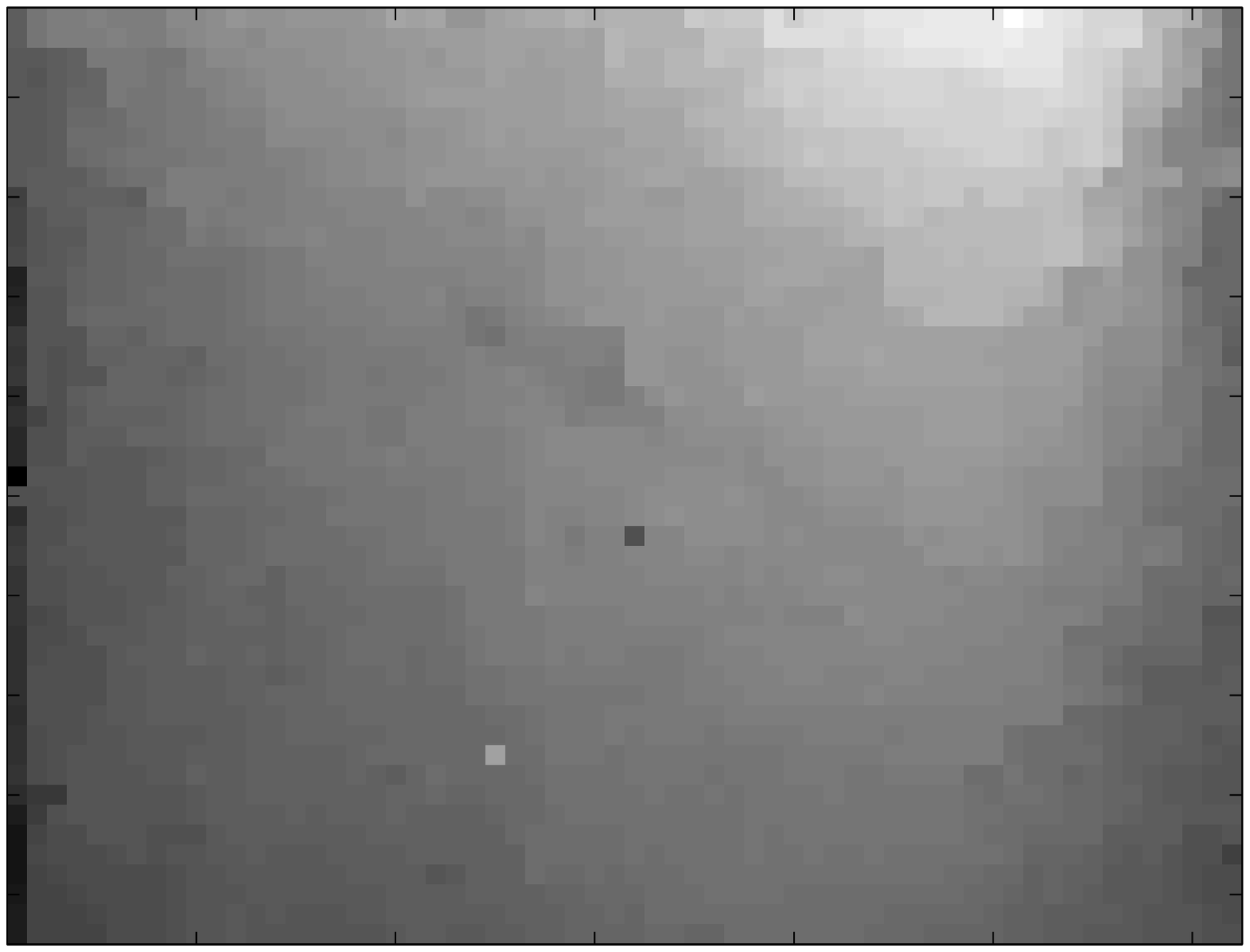}\hspace{-.05in}} \\
\hline
\end{tabular}
\end{center}
\end{table*}

\vspace{-.1in}
{
\subsection{Run Time}
Our final analysis is performed to study the computational processing time taken by each algorithm to recover \textit{Office Environment} video data.  Inverse operation ${\bf O}^{-1}$ in line \ref{algo:f-sub} of Algorithm \ref{algo_BC_ADMM} is the computational bottleneck of the BC-ADMM. We used eigenvalue decomposition \texttt{eig} in \textsc{Matlab} to decompose ${\bf D}{\bf D}^T$ in (\ref{eq18}) to acquire corresponding  eigenvalues ${\bf\Lambda}$ and eigenvectors ${\bf Q}$ to perform the least square minimization problem in line \ref{algo:f-sub}. The speed of our algorithm is not optimized, where fast implementation can be done by FFTs and FSTs in \cite{NgChanTang:1999,Capizzano:2003, Ng:2004} to handle different BCs. Also, multicore processing can be used to assign parallel minimization steps in BC-ADMM. In all competing algorithms: TFOCS, NESTA, and TVAL3 we ignored duplicating the data in order to have a fair comparison platform with BC-ADMM. All algorithms are run on a Windows 7 Dell Optiplex 780 desktop computer with an Intel(R) Core(TM)2 Quad CPU with 3.0 GHz processors and 4 GB of RAM. Here, only one core is used by \textsc{Matlab} for processing. The number of iterations as well as elapsed time to complete the minimization is demonstrated in Table \ref{table_convergence}. The values are averages over $10$ Monte-Carlo trials on each experiment. The ranking order in terms of average time elapsed to complete one iteration in each algorithm are: (1) NESTA $=0.43$, (2) TVAL3$=0.47$, (3) TFOCS  $=0.85$, and (4) BC-ADMM: $=1.02$ second. Our algorithm takes almost double time to complete each iteration in BC-ADMM. This is because the Toeplitz matrix ${\bf D}$ is constructed by $L=27$ tab filter, where in 3D tensorial calculation the order of complexity increases by $\mathcal{O}(L^3)$ due to curse of dimensionality. Nevertheless, the number of iterations in BC-ADMM to complete an experiment remains as low as TVAL3, where both use ADMM for minimization.}

\begin{table*}[htp]
\caption{Iterations and elapsed time (second) completed by each algorithm on different sampling rates}\vspace{-.1in}
\label{table_convergence}
\centering
\begin{tabular}{ |c|c|c|c|c|c|c|c|c|c|c| } 
\hline
\multirow{2}{*}{$\frac{|{\bf\Omega}|}{mnN}$} & \multicolumn{2}{ |c| }{AR-BC-ADMM} & \multicolumn{2}{ |c| }{P-BC-ADMM} & \multicolumn{2}{ |c| }{TFOCS} & \multicolumn{2}{ |c| }{NESTA} & \multicolumn{2}{ |c| }{TVAL3} \\ 
\cline{2-11}
& Iterations & CPU Time & Iterations & CPU Time & Iterations & CPU Time & Iterations & CPU Time & Iterations & CPU Time\\ \hline
0.01 &  215 & 217.98 & 266 & 271.72 & 400 & 334.08 & 543 & 223.29& 269 & 118.93 \\ \hline
0.02 &  180 & 180.96 & 210 & 217.24 & 400 & 338.10 & 344 & 142.74& 211 &  97.47 \\ \hline
0.03 &  165 & 170.66 & 192 & 197.58 & 400 & 336.50 & 307 & 128.44& 183 &  83.32 \\ \hline
0.05 &  146 & 152.56 & 171 & 175.42 & 400 & 338.13 & 268 & 114.41& 143 &  67.27 \\ \hline
0.10 &  127 & 134.35 & 144 & 148.66 & 400 & 339.83 & 226 & 97.43  & 113 &  53.91 \\ \hline
0.15 &  118 & 124.43 & 129 & 133.76 & 400 & 341.33 & 211 & 91.04  & 105 &  51.25 \\ \hline
0.25 &  108 & 111.98 & 114 & 116.32 & 400 & 351.30 & 188 & 83.58  &   86 &  41.86 \\ \hline        
\end{tabular}
\end{table*}\vspace{-.0in}



\vspace{-.1in}
{
\section{Concluding Remarks}\label{Sec:Conclusion}
We have proposed an alternative TV regularization model that utilizes HO accuracy differential FIR filter to encode high frequency components pertinent to rapid transitions/edge information in the signal. The proposed TV is extended to tensorial representation to encode video features by implicitly decomposing each space-time dimension. This is used as a new regularizer in CVS problem to recover frames from compressed samples. The solution is driven by ADMM to decouple the optimization steps and solve a unique quadratic minimization problem that is capable of handling different BCs. In particular, we suggest deploying AR-BC in BC-ADMM to avoid restoration artefacts. Extensive experimentation with commonly used video data sets suggests that the proposed TV regularizer is capable of accurately recovering video images and preserve much edge information in frames as well as discontinuous motion features in temporal domain. The proposed TV can be easily reduced to the problem of 2D imaging for any restoration problem such as compressed imaging, denoising, and deblurring.}\vspace{-.1in}

\appendices
%
\section{High-Order of Accuracy Differentiation}\label{Sec:NRD}
\begin{figure}[htp]
\centerline{
\subfigure[]{\includegraphics[width=0.2\textwidth,height=0.2\textwidth]{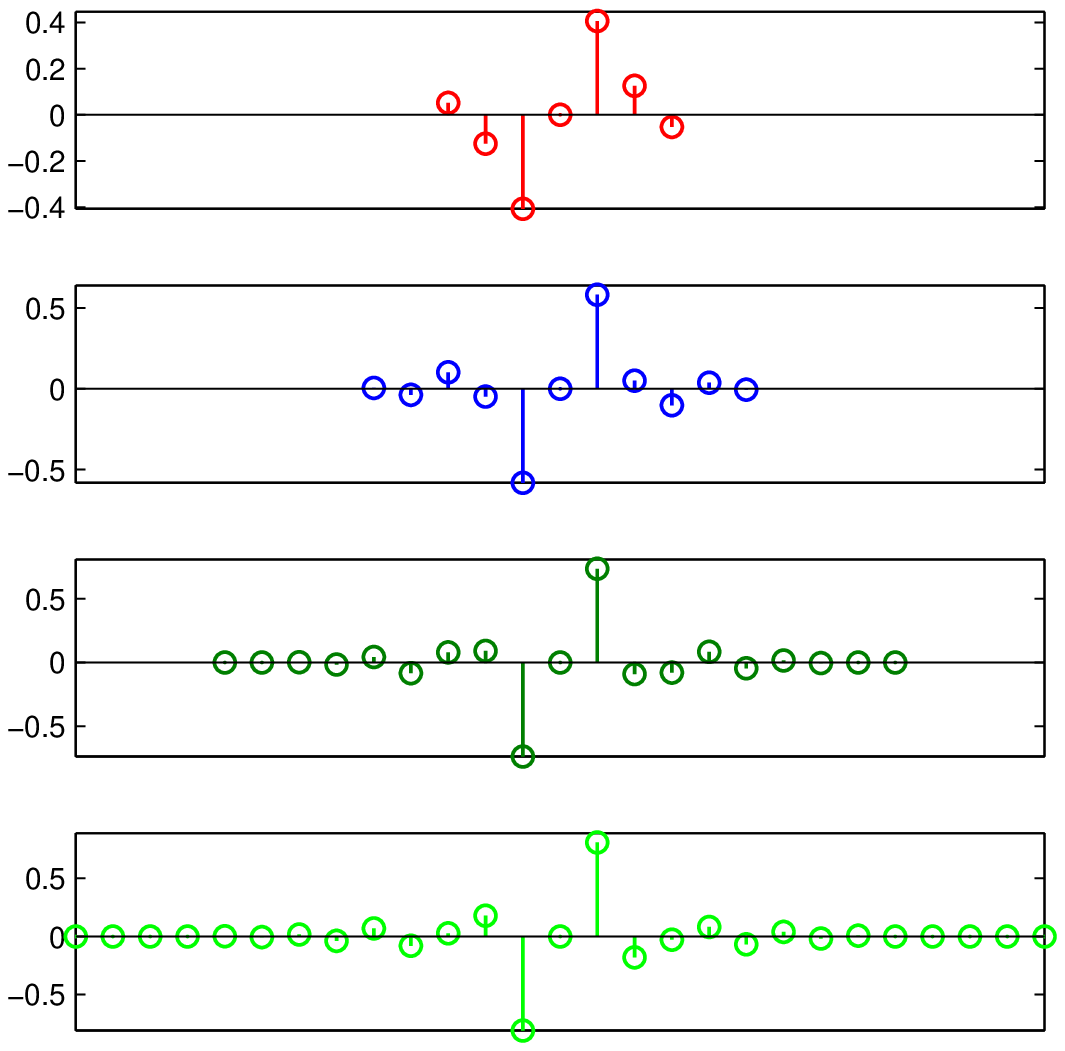}
\label{Fig:drvKernels}}
\subfigure[]{\includegraphics[width=0.185\textwidth]{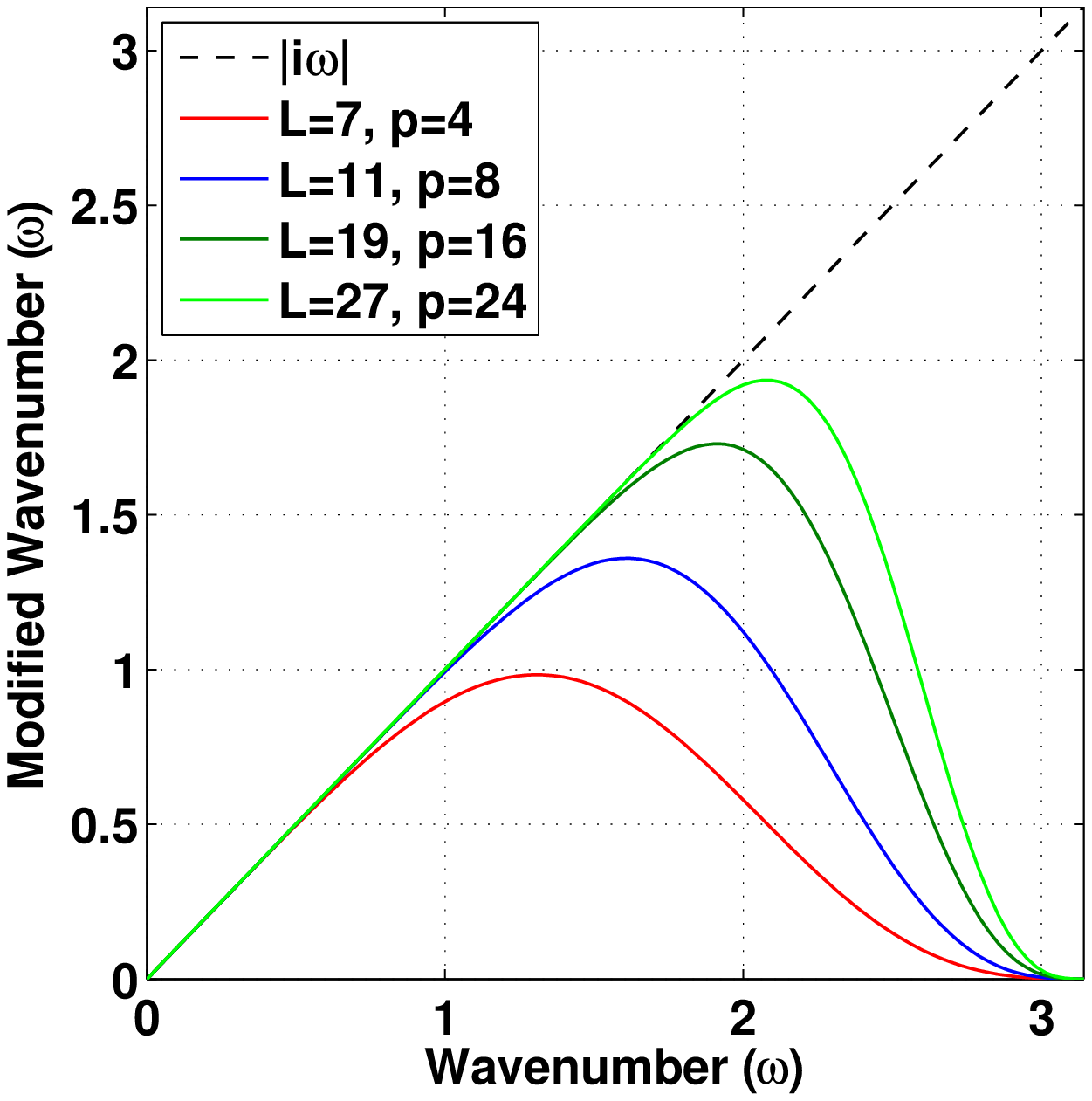}
\label{Fig:fitlerResRN_wavenumber}}}
\vspace{-.15in}
\caption{(a) Noise-robust derivative kernels generated at different level of accuracies, i.e. $L=\{7,11,19,27\}$; (b) Filter response}
\label{Fig:fitlerResRN}
\end{figure}

The general format of the numerical differentiators are designed to be an anti-symmetric digital FIR filter with $L-1$ accuracy order, i.e.
\begin{equation}\label{TD7}
v^{\prime}_{k} \approx \frac{1}{T}\sum^{(L-1)/2}_{\ell=1}{d_\ell\cdot\left(v_{k+\ell}-v_{k-\ell}\right)},
\end{equation}
where, $d_\ell$ is the associated filter coefficient. The corresponding transfer function is calculated by discrete Fourier transform (DFT) of (\ref{TD7}), i.e.
\begin{equation}\label{TD8}
H(\omega, T) \triangleq \frac{{\mathcal F}\left\{v^{\prime}_k\right\}}{{\mathcal F}\left\{v_k\right\}}=
\frac{2i}{T}\sum^{(L-1)/2}_{\ell=1}{d_\ell\sin{\ell T\omega}}.
\end{equation}
The derivative operator, in continuous  domain, carries the following transfer function $H_a(\omega)=i\omega$. Our goal, by deploying a digital derivative filter, is $H(\omega, T)$ to become closer to $H_a(\omega)$. This is the main assumption in polynomial interpolation in \cite{Jahne:1993:SIP:528863}. Holoborodko \cite{Holoborodko:2008} proposed a solution to approximate derivatives based on the concept of high frequency band suppression. The derivative kernel is expected to satisfy high precision on low frequencies and suppress higher frequency artifacts. Based on these two assumptions, the designed filter $H(\omega, T)$ is likely to behave as low-pass filter by taking the following in to the consideration:
\begin{equation}\label{TD9}
\frac{\partial^k H}{{\partial\omega}^k}(\omega, T){\Big{\vert}}_{\omega=0}= \frac{\partial^k H_a}{{\partial\omega}^k}(\omega){\Big{\vert}}_{\omega=0},~~~k\in\{0,\hdots,p\}
\end{equation}
to meet the first assumption, and
\begin{equation}\label{TD10}
\frac{\partial^k H}{{\partial\omega}^k}(\omega, T){\Big{\vert}}_{\omega=\pi} = 0,~~~k\in\{0,\hdots,q\}
\end{equation}
to behave as a low-pass filter. Substituting (\ref{TD8}) in  (\ref{TD9}-\ref{TD10}) and grouping the terms, yields
\begin{equation}\label{TD11}
\left\{
\begin{array}{lr}
\sum^{(L-1)/2}_{\ell=1}{\ell d_\ell}=1/2 & ,~k=0\\
\sum^{(L-1)/2}_{\ell=1}{\ell^{2k+1}d_\ell}=0 & ,~k=\{1,\hdots,\lfloor\frac{p-1}{2}\rfloor\}\\
\sum^{(L-1)/2}_{\ell=1}{(-1)^\ell \ell^{2k+1}d_\ell}=0 & ,~k=\{0,\hdots,\lfloor\frac{q-1}{2}\rfloor\}.
\end{array}
\right.
\end{equation}
The approximated coefficients $d_l$ identify the designed filter $H(\omega, T)$ in (\ref{TD8}). This filter will be exact at $\omega=0$ similar to ideal filter response $H_a(\omega)$ up to $p$ polynomial degree in frequency domain. $p$ is arbitrarily chosen up to $p<L-1$. In order to derive a unique solution via (\ref{TD11}), $q$ is chosen such that $\lfloor\frac{q-1}{2}\rfloor=\frac{L-1}{2}-\lfloor\frac{p-1}{2}\rfloor-2$. The value $q$ will define the order of tangency of  $H(\omega, T)$ at $\omega=\pi$. Accordingly, the corresponding system of matrix equations of (\ref{TD11}) is presented in (\ref{TD12}), where $n_o=\lfloor\frac{p-1}{2}\rfloor+1$, $m_o=\lfloor\frac{q-1}{2}\rfloor+1$, and $A=\left(\frac{L-1}{2}\right)$. The modified frequencies from the designed filter is shown in Figure \ref{Fig:fitlerResRN} for different accuracy levels, i.e. $L$ with exact polynomials $p=L-2$.

\begin{equation}\label{TD12}
\small{\left[\hspace{-.05in}\begin{array}{llll}
(1)^1 & (2)^1 & \hdots & A^1 \\
(1)^3 & (2)^3 & \hdots & A^3 \\
\vdots & \vdots& \ddots & \vdots \\
(1)^{n_o} & (2)^{n_o}  & \hdots & A^{n_o} \\
(-1)^1 (1)^1 & (-1)^2 (2)^1 & \hdots & (-1)^{A} A^1 \\
(-1)^1 (1)^3 & (-1)^2 (2)^3 & \hdots & (-1)^{A} A^3 \\
\vdots & \vdots& \ddots & \vdots \\
(-1)^1 (1)^{m_o} & (-1)^2 (2)^{m_o} & \hdots & (-1)^{A} A^{m_o} \\
\end{array}\hspace{-.05in}\right]
\left[\hspace{-.05in}\begin{array}{c}
d_1 \\ d_2 \\ \vdots \\ d_{A}
\end{array}\hspace{-.05in}\right]=
\left[\hspace{-.05in}\begin{array}{c}
\frac{1}{2}\\ 0 \\ \vdots \\ 0 \\ 0 \\ 0 \\ \vdots \\ 0
\end{array}\hspace{-.05in}\right],}
\end{equation}


As it shown, the modified modes are close to ideal responses around $\omega=0$. In contrast, the filter attenuates the higher frequencies and suppresses to zero close to $\omega=\pi$. The proposed filter in \cite{Holoborodko:2008} is capable of suppressing different frequency modes, which is dependent to tuning of the filter parameters. For a fixed length $L$, the suppression of different modes varies by the choices of $p$ which defines the exactness of low frequencies to be preserved.\vspace{-.1in}

\section{Proof of Proposition \ref{Proposition:Rank} }\label{Sec:Appendix2}
\begin{IEEEproof}
The operator $\nabla=[A; B; C]$ is a partitioned matrix with three block matrices, where each projects the input signal into directional derivative spaces.
We use the theorem from Marsaglia et al.\cite{MatsagliaStyan:1974} to find the rank of partitioned matrices which states the following equality for any two conformable matrices $A$ and $B$:
\begin{equation}\label{eq20}
r\left[\begin{array}{c}A\\ B\end{array}\right]=r(A)+r\left(B(I-A^{\dagger}A)\right),
\end{equation}
where $A^{\dagger}$ is Moore-Penrose pseudo inverse of $A$. First, we find the rank of two block matrices and extend the results to three sub-matrices. Let us substitute the first two block matrices of $\nabla$ into (\ref{eq20}) to find their rank, i.e. $A={\bf I}_N\otimes {\bf I}_n \otimes {\bf D}_m$ and $B={\bf I}_N\otimes {\bf D}_n \otimes {\bf I}_m$. We decompose ${\bf D}_m={\bf U}_m{\bf \Sigma}_m {\bf V}^{*}_m$ using singular value decomposition, where ${\bf U}_m,{\bf V}_m\in\mathbb{R}^{m\times m}$ are orthonormal matrices and ${\bf \Sigma}_m\in\mathbb{R}^{m\times m}$ is diagonal singular values which contains $m-p$ non-zeros. The Moore-Penrose pseudo inverse of ${\bf D}_m$ is defined by ${\bf D}^{\dagger}_m={\bf V}_m{\bf \Sigma}^{\dagger}_m {\bf U}^{*}_m$, where ${\bf \Sigma}^{\dagger}_m$ is reciprocal of non-zero diagonal values and leaving the zeros in place \cite{Golub:1996}. Substituting this decomposition in $A^{\dagger}A$ yields,
\begin{eqnarray}\label{eq21}
\lefteqn{I-A^{\dagger}A=} \\
& & {\bf I}_{nmN}-\left({\bf I}_N\otimes {\bf I}_n \otimes {\bf V}_m{\bf \Sigma}^{\dagger}_m {\bf U}^{*}_m\right)\left({\bf I}_N\otimes {\bf I}_n \otimes {\bf U}_m{\bf \Sigma}_m {\bf V}^{*}_m\right) \nonumber \\
& &={\bf I}_{nmN}-{\bf I}_N\otimes {\bf I}_n \otimes {\bf V}_m{{\bf \Sigma}^{\dagger}_m{\bf \Sigma}_m}{\bf V}^{*}_m \nonumber \\
& &={\bf I}_N\otimes {\bf I}_n \otimes \left({\bf I}_m-{\bf V}_m{{\bf \Sigma}^{\dagger}_m{\bf \Sigma}_m}{\bf V}^{*}_m\right). \nonumber
\end{eqnarray}
Note that ${\bf I}_m-{\bf V}_m{{\bf \Sigma}^{\dagger}_m{\bf \Sigma}_m}{\bf V}^{*}_m={\bf V}_m{({\bf I}_m-{\bf \Sigma}^{\dagger}_m{\bf \Sigma}_m)}{\bf V}^{*}_m$ and multiplying $B$ to (\ref{eq21}) from left yields,
\begin{eqnarray}\label{eq22}
\lefteqn{B(I-A^{\dagger}A) =} \\
& & \left({\bf I}_N\otimes {\bf D}_n \otimes {\bf I}_m\right)\left({\bf I}_N\otimes {\bf I}_n \otimes {\bf V}_m{({\bf I}_m-{\bf \Sigma}^{\dagger}_m{\bf \Sigma}_m)}{\bf V}^{*}_m\right) \nonumber \\
& &\hspace{.5in}={\bf I}_N\otimes {\bf D}_n \otimes {\bf V}_m{({\bf I}_m-{\bf \Sigma}^{\dagger}_m{\bf \Sigma}_m)}{\bf V}^{*}_m.\nonumber
\end{eqnarray}
To determine the rank in (\ref{eq20}), we use the corollary from \cite{laub:2005}, where the rank of Kronecker product of two matrices equals to the product of the ranks of two matrices, i.e. $r(A\otimes B)=r(A)r(B)$. By substituting (\ref{eq22}) in (\ref{eq20}), the rank is determined as follows, 
\begin{equation}\label{eq23}
r\left[\begin{array}{c}A\\ B\end{array}\right]=Nn(m-p)+N(n-p)p=N(nm-p^2).
\end{equation}
Here, $r\left[{\bf V}_m{({\bf I}_m-{\bf \Sigma}^{\dagger}_m{\bf \Sigma}_m)}{\bf V}^{*}_m\right]=p$, since ${\bf I}_m-{\bf \Sigma}^{\dagger}_m{\bf \Sigma}_m$ is a diagonal matrix containing $p$ non-zero values. Next, the third block matrix is considered to find the final rank of $\nabla$ as follows,
\begin{equation}\label{eq24}
r(\nabla)=r\left[\begin{array}{c}A\\ B\\ C\end{array}\right]=r\left[\begin{array}{c}A\\ B\end{array}\right]+r(C(I-\left[\begin{array}{c}A\\ B\end{array}\right]^{\dagger}\left[\begin{array}{c}A\\ B\end{array}\right])),
\end{equation}
where $C={\bf D}_N\otimes {\bf I}_n \otimes {\bf I}_m$. We use the SVD decomposition in $A={\bf I}_N\otimes {\bf I}_n \otimes {\bf U}_m{\bf \Sigma}_m {\bf V}^{*}_m$ and $B={\bf I}_N\otimes {\bf U}_n{\bf \Sigma}_n {\bf V}^{*}_n \otimes {\bf I}_m$ to simplify the following expression,
\begin{eqnarray}\label{eq25}
\lefteqn{C(I-\left[\begin{array}{c}A\\ B\end{array}\right]^{\dagger}\left[\begin{array}{c}A\\ B\end{array}\right])} \\
& & \hspace{-.2in}=({\bf D}_N\otimes {\bf I}_n \otimes {\bf I}_m)({\bf I}_{mnN}-{\bf I}_N \otimes {\bf I}_n \otimes {\bf V}_m{{\bf \Sigma}^{\dagger}_m{\bf \Sigma}_m}{\bf V}^{*}_m \nonumber \\
& &\hspace{1.2in} - {\bf I}_N \otimes {\bf V}_n{{\bf \Sigma}^{\dagger}_n{\bf \Sigma}_n}{\bf V}^{*}_n \otimes {\bf I}_m) \nonumber \\
& & \hspace{-.2in} = ({\bf D}_m\otimes {\bf I}_n \otimes {\bf I}_N)[{\bf I}_N \otimes {\bf V}_n{({\bf I}_n-{\bf \Sigma}^{\dagger}_n{\bf \Sigma}_n)}{\bf V}^{*}_n \otimes \nonumber\\
& & \hspace{1.2in} {\bf V}_m{({\bf I}_m-{\bf \Sigma}^{\dagger}_m{\bf \Sigma}_m)}{\bf V}^{*}_m]\hspace{.63in} \nonumber \\
& & \hspace{-.2in} = {\bf D}_N \otimes {\bf V}_n{({\bf I}_n-{\bf \Sigma}^{\dagger}_n{\bf \Sigma}_n)}{\bf V}^{*}_n \otimes {\bf V}_m{({\bf I}_m-{\bf \Sigma}^{\dagger}_m{\bf \Sigma}_m)}{\bf V}^{*}_m.\nonumber
\end{eqnarray}
By substituting (\ref{eq25}) and (\ref{eq23}) in (\ref{eq24}) the final rank is determined by
\begin{equation}\label{eq26}
r(\nabla)=N(nm-p^2)+(N-p)p^2=Nnm-p^3,
\end{equation}

Next we use the results from ((\ref{eq26})) to find the rank $r({\bf G})$. First, we find the rank of the following
\begin{equation}\label{eq27}
r\left[\begin{array}{c}{\bf \nabla}\\ \bar{\bf\Psi}\end{array}\right]=r({\bf \nabla})+r\left(\bar{\bf\Psi}({\bf I}_{mnN}-{\bf \nabla}^{\dagger}{\bf \nabla})\right).
\end{equation}
To simplify the ${\bf \nabla}^{\dagger}{\bf \nabla}$ in (\ref{eq27}), we use the SVD decomposition in all three derivative matrices deployed in ${\bf \nabla}$, i.e. ${\bf D}_m={\bf U}_m{\bf \Sigma}_m {\bf V}^{*}_m$, ${\bf D}_n={\bf U}_n{\bf \Sigma}_n {\bf V}^{*}_n$ and ${\bf D}_N={\bf U}_N{\bf \Sigma}_N {\bf V}^{*}_N$ and express as follows
\begin{equation}\label{eq28}
{\bf \nabla}^{\dagger}{\bf \nabla}={\bf V}_N{{\bf \Sigma}^{\dagger}_N{\bf \Sigma}_N}{\bf V}^{*}_N \otimes{\bf V}_n{{\bf \Sigma}^{\dagger}_n{\bf \Sigma}_n}{\bf V}^{*}_n \otimes {\bf V}_m{{\bf \Sigma}^{\dagger}_m{\bf \Sigma}_m}{\bf V}^{*}_m.
\end{equation}
Consequently, the following term will be simplified to,
\begin{eqnarray}\label{eq29}
\lefteqn{\bar{\bf\Psi}({\bf I}_{mnN}-{\bf \nabla}^{\dagger}{\bf \nabla})={\bf V}_N{({\bf I}_N-{\bf \Sigma}^{\dagger}_N{\bf \Sigma}_N)}{\bf V}^{*}_N \otimes} \\
& & {\bf\Psi}_n{\bf V}_n{({\bf I}_n-{\bf \Sigma}^{\dagger}_n{\bf \Sigma}_n)}{\bf V}^{*}_n \otimes {\bf\Psi}_m{\bf V}_m{({\bf I}_m-{\bf \Sigma}^{\dagger}_m{\bf \Sigma}_m)}{\bf V}^{*}_m. \nonumber
\end{eqnarray}
The rank of (\ref{eq29}) is equal to $r(\bar{\bf\Psi}({\bf I}_{mnN}-{\bf \nabla}^{\dagger}{\bf \nabla}))=p^3$. This is because the individual ranks from (\ref{eq29}) are found by $r\left[{\bf\Psi}_n{\bf V}_n{({\bf I}_n-{\bf \Sigma}^{\dagger}_n{\bf \Sigma}_n)}{\bf V}^{*}_n\right]=p$ since the orthonormal projection on ${\bf\Psi}_n$ is non-rank-deficient. Hence, by substituting the results from (\ref{eq26}) and (\ref{eq29}) into (\ref{eq27}), the final rank yields,
\begin{equation}\label{eq30}
r\left[\begin{array}{c}{\bf \nabla}\\ \bar{\bf\Psi}\end{array}\right]=mnN-p^3 + p^3=mnN,
\end{equation}
where, it is a full column rank partitioned matrix. Considering the third block matrix ${\bf\Phi}$ in ${\bf G}$ and using  similar analysis yields a full column rank matrix for ${\bf G}$, i.e. $r({\bf G})=mnN$.
\end{IEEEproof}

\vspace{-.1in}
\section*{Acknowledgment}
The authors would like to thank Dr. Stephan Becker for helpful suggestions to revise and improve TFOCS' TV toolbox to 3D problem. 
\vspace{-.1in}

\ifCLASSOPTIONcaptionsoff
  \newpage
\fi



\bibliographystyle{IEEEtran}
\bibliography{Bibliography,myref}
\end{document}